\input amstex
\documentstyle{amsppt}
\magnification=\magstep0
\define\cc{\Bbb C}
\define\z{\Bbb Z}
\define\r{\Bbb R}

\define\N{\Bbb N}
\define\Q{\Bbb Q}
\define\jj{\Bbb J}

\define\A{\Cal A}
\define\h{\Cal D}
\define\E{\Cal E}
\define\m{\Cal M}

\define\f{\Cal S}
\define\a{\Cal L}

\define\n{\Cal F}

\define\la{\lambda}
\define\om{\omega}

\define\e{\varepsilon}
\define\va{\varphi }
\define\CB#1{\Cal C_b(#1)}
\define\st{\subset }

\topmatter
 \title
 Equality of uniform  and Carleman spectra for  bounded measurable functions
  \endtitle
 \endkeywords
 \author
  Bolis Basit and  Alan J. Pryde
\endauthor
 \abstract
{In this paper we study various types of spectra of functions
$\phi:\jj\to X$, where $\jj\in\{\r_+,\r\}$ and $X$ is a complex
Banach space. We show that uniform spectrum defined in [15]
coincides with Carleman spectrum for $\phi\in L^{\infty}(\r,X)$.
This result holds true also for Laplace (half-line) spectrum for
$\phi\in L^{\infty}(\r_+,X)$. We also indicate a class of bounded
measurable functions for which Laplace spectrum and Carleman
spectrum are equal.}
\endabstract
\endtopmatter
\rightheadtext{ Equality of uniform  and Carleman spectra}
\leftheadtext{B.Basit,  A. J. Pryde}
  \TagsOnRight
\document
\pageno=1 \baselineskip=19pt

\head{\S 0. Introduction}\endhead

In this paper we study various types of spectra of functions
$\phi:\jj\to X$, where $\jj\in\{\r_+,\r\}$ and $X$ is a complex
Banach space. If $\phi\in BUC(\jj,X)$, the space of bounded
uniformly continuous functions, it is easily verified using
Bochner integration that the Laplace transform $ \Cal{L} \phi:
\cc_+\to X$, where $ \Cal{L} \phi(\la) =\int_0^{\infty}\, e^{-\la
\,t}\phi(t)\, dt$, has a holomorphic extension in a neighbourhood
of   some point $i\om\in i\r$ if and only if $ \Cal{L}_u \phi:
\cc_+\to BUC(\jj,X)$, where $ \Cal{L}_u \phi(\la)
(s)=\int_0^{\infty}\, e^{-\la \,t}\phi (t+s)\, dt$   for $s\in
\jj$, has   a holomorphic extension in a neighbourhood of $i\om$.
The Laplace spectrum  is $sp^{\Cal{L}} (\phi):=\{\om\in\r: i\om$
is a singular point for $ {\Cal{L}} \phi\}$ and the  uniform
Laplace spectrum is $sp^{\Cal{L}_u} (\phi):= \{\om\in\r: i\om$ is
a singular point for $ {\Cal{L}}_u \phi\}$.   The Carleman
transform (see (1.3)), spectrum (see (1.7)) and the uniform
spectrum are defined similarly. In [15], [24], the uniform
Carleman spectrum $ sp^{\Cal{C}_u}( \phi)$ for $\phi\in BC(\r,X)$
is introduced  and it is shown  that $ sp^{\Cal{C}} (\phi)\st
sp^{\Cal{C}_u} (\phi)$. Many properties of $ sp^{\Cal{C}} (\phi)$
are shown to hold true for $ sp^{\Cal{C}_u} (\phi)$ (see [15,
Proposition 2.3]); however, equality is not established.

In \S 1 we collect notation and  definitions and prove some
preliminaries. In particular, we prove in Proposition 1.5 a useful
necessary condition for a point $\om$ to be in the complement of
the  reduced Beurling spectrum $sp_{C_0(\jj,X)}(\phi)$.

In \S 2,  we  establish some tools which enable us to calculate
the Laplace and weak Laplace spectra and relate them to the
reduced Beurling  spectrum relative to a class $\A$ satisfying
(1.12). Using a new property (Theorem 2.4(i)) of the  weak Laplace
spectrum we give simple proofs of several tauberian results of
Ingham [20], [22] and their generalizations by Chill and
 others (see [2], [3],  [4, 4.10, p. 332],  [11] and references therein).

In \S 3, we study bounded $C_0$-semigroups $T(t)$, $t\in\jj$ with
generator $A$. We give  proofs  of the identities \,\,\,
 $
i\,sp^{\Cal{L}} (T(\cdot))= \,\sigma (A)\cap i\,\r$\,\,\ and, when
$\jj=\r$, of the equality of the  Arveson, Beurling and Carleman
spectra of both the orbits $T(\cdot)x$ and groups $T(\cdot)$ (see
(3.6), (3.7)); in particular, in Corollary 3.3 we show that
Arveson spectrum $sp^A(T(\cdot)x)$ is equal $\sigma_u(A,x)$, the
local unitary spectrum of $A$ at $x$ defined in [10, \S 3]. Since
isometric semigroups on $\r_+$ can be extended uniquely  to
isometric groups on $\r$, Corollary 3.3 extends [10, Theorem
2.2].

In \S 4, we prove  $ sp^{\Cal{L}_u} (\phi)= sp^{\Cal{L}} (\phi)$
and, when $\jj=\r$,
  $ sp^{\Cal{C}_u} (\phi)= sp^{\Cal{C}} (\phi)$ for $\phi\in L^{\infty}
(\jj,X)$. These results seem new for $BC(\r,X)$ and the proofs are
new even for  $BUC(\r,X)$.

Finally, in section \S5, we indicate a subclass of $L^{\infty}
(\r,X)$ for which the Laplace spectrum   coincides with the
Carleman spectrum. This class includes almost periodic, almost
automorphic, Levitan almost periodic and recurrent functions.

 \head{\S 1. Notation, Definitions and preliminaries}\endhead

In the following $\r_+=[0,\infty)$, $\jj\in\{\r_+,\r\}$, $\N= \{1,
2,\cdots \}$, $\N_0=\{0 \} \cup \N$,  $\cc_+
 =\{\la\in \cc: \text{Re \,}\la
>0\}$ and $\cc_-=\{\la\in \cc: \text{Re \,}\la < 0\}$.
 Denote by $X$ a   complex Banach space. If $Y$, $Z$ are locally
 convex topological spaces, $L(Y,Z)$ will denote the space of all bounded linear operators from $Y$ to
 $Z$ and   $ L(Y)=L(Y,Y)$.
   The Schwarz space of rapidly decreasing functions is denoted by $\f
   (\r)$  and
  $\f' (\r,X)=L(\f(\r),X)$ (see [31]) is the  space of $X$-valued tempered distributions
  on $\r$. The action of an element $ S\in \f' (\r,X)$ on
  $f\in\f(\r)$ is denoted  $<S,f>$.
 If $\phi$ is a  $ X$-valued function defined on  $\jj,$
then   $\phi_s$,  $\Delta _s\phi$ will stand for  functions
defined on $\jj$ by $\phi_{s}(t) = \phi(t+s)$, $\Delta _s\phi (t)=
\phi_{s}(t)-\phi (t) $ for all $s\in \jj $,
  $|\phi|$ will denote the function $|\phi|(t):= ||\phi(t)||$  for all $t\in \jj$ and $||\phi|| _{\infty} :=
  \text {  sup}_{t\in \jj} ||\phi(t)||$.
 If $\phi \in L_{loc}^1 (\jj, X)$, then $P\phi$ and $M_h\phi$ will denote  the
indefinite integral  and mollifier  of $\phi$ defined respectively
by $P\phi(t) = \int_{0}^{t} \phi (s)\,ds$ and $M_h \phi (t)=
(1/h)\int_0^h \phi (t+s)\,ds$ for $h >0$.
 For $f\in
   L^1(\r,\cc)$ and $\phi\in
   L^{\infty}(\r,X)$ or $f\in
   L^1(\r,X)$ and $\phi\in
   L^{\infty}(\r,\cc)$ the Fourier transform $\hat{f}$ and convolution $\phi*f$ are
   defined respectively by $\hat{f} (\om)=\int_{-\infty}^{\infty} \gamma_{-\om} (t)\, f(t)\,
   dt$ and $\phi*f (t)= \int_{-\infty}^{\infty}\phi(t-s) f(s)\,
   ds$, where  $\gamma_{\om} (t)= e^{i\, \om t} $.  The Fourier transform of $ S\in \f' (\r,X)$ is
  the tempered distribution
   $\widehat{S}$  defined by $<\widehat{S},f> =<S,\hat{f}>$ for all $f\in \f(\r)$. All integrals are
Lebesgue-Bochner integrals (see [4, pp. 6-15], [16, p. 318],
   [19, p. 76]).

 We recall [19, Definition 3.5.5,
p. 74] that an operator valued function $F: \jj\to L(X)$ is
strongly measurable (strongly integrable) if $F(\cdot)x$ is
measurable  (integrable) for each $x\in X$. We denote by
$L_s^{\infty} (\jj,L(X))$ the Banach space of all (essentially)
bounded strongly measurable operator-valued functions $F$ such
that $|F(\cdot)|\in L^{\infty} (\jj)$. Since $F(\cdot)x\in
L^{\infty} (\jj, X)$ for each $x\in X$  and $a\in \jj$ we may
define
 the $strong\,\, integral$ of $F$ by $(\int_0^{a}\, F(t)\,dt)x :=\int_0^{a}\, F (t)x\,dt$.  Similarly for $h>0$ and $t\in \jj$,
   the $strong\,\, mollifier$ $M_h F(t) \in L(X)$ is
defined by $(M_h F(t))x =M_h (F(t)x)=(1/h)\int_0^h F(t+s)x\, ds$.
  For example if
$T(t)\in L(X)$ for $t\in \jj$ is a bounded $C_0$-semigroup, then
by [16, p 616], $|T(\cdot)|\in L^{\infty} (\jj)$ and so
$T(\cdot)\in L_s^{\infty} (\jj, L(X))$. In particular, the
translation semi-groups  $S^{\jj}(t): BUC(\jj,X)\to BUC(\jj,X)$
  defined by $S^{\jj}(t)\phi=\phi_t$ for $t\in\jj$ are strongly continuous. Hence
$S^{\jj}(\cdot)\in L_s^{\infty} (\jj, L(BUC(\jj,X)))$. For
simplicity we write  sometimes $S(\cdot)=S^{\r}(\cdot)$ and
$S^{+}(\cdot)=S^{\r_+}(\cdot)$.

In this paper we  consider the space $L^1_{loc} (\r_+,X)$  as a
subspace of $L^1_{loc} (\r,X)$ by  identifying a
 function $\phi$ defined on $\r_+$ with its extension by $0$ to
$\r$.

If $\phi\in L^{1}_{loc}(\jj,X)\cap \f'(\r,X)$, then

\qquad  $f\phi|\,\r_+,\,\,  e_{\la}\, \phi|\,\r_+ \in L^1(\r_+,X)$
for $f\in \f(\r)$ and $e_{\la }(t)=e^{-\la \, t}$, \,\,    $\la
\in \cc_+$,

\noindent and the $Laplace\,\, transform$ $\Cal {L}\phi$  is
defined by

(1.1) \qquad $\Cal {L}\phi(\la)=\int_0^{\infty}\, e^{-\la \,t}\phi
(t)\, dt$ \,\,\,\, for\,\,\,\, $\la \in \cc_+ $.

\noindent Clearly,   $\Cal {L}\phi$ is holomorphic on
 $\cc_+$. In particular, if $F\in L_s^{\infty} (\jj,
L(X))$ then
 $\Cal {L} F(\cdot)\, x$ is holomorphic on
$\cc_+$ for each $x\in X$.  For $\la\in \cc_+$ we define  $\Cal
{L}F (\la)$  by

(1.2) \qquad  \,\,\,\, $\Cal {L} F(\la)x= \Cal {L}F(\cdot)\,
x(\la)$.

 Denote the $Carleman$ $transform$ of $\phi \in
L^{1}_{loc}(\r,X)\cap \f'(\r,X)$    by

(1.3) \qquad   $\Cal {C} {\phi} (\la)=
 \cases { \Cal{L^+}\phi (\la)= \int _0^{\infty}\, e^{-\la \, t} \phi (t)\,
dt,\qquad \,{\text{if\,\,} \la
 \in \cc_+}}\\
  { \Cal{L^-}\phi(\la) = - \int _0^{\infty}\, e^{\la\, t} \phi(-t)\, dt,\qquad {\text{\,
  if\,\,}
 \la \in \cc_-}.}\endcases$

\noindent Then $\Cal {C} {\phi}$ is an $X$-valued function which
is holomorphic on $\cc\setminus i\,\r $.

Similarly,  since the Carleman transform $\Cal{C} F(\cdot)x$ of
$F\in L_s^{\infty}(\r,L(X))$  is holomorphic on $\cc\setminus
i\,\r=\cc_+\cup \cc_-$ for each $x\in X$ we define  $\Cal {C}F
(\la)$ by

(1.4) \qquad  \,\,\, $\Cal{C} F(\la)x= \Cal{C} F(\cdot)x(\la)$.

\noindent It is easily verified that  $\Cal {L} F(\la)\in L(X)$
for $\la\in \cc_+$ and  $\Cal {C} F(\la)\in L(X)$ for $\la\in
\cc\setminus i\,\r$.  Since
  $\Cal {L} F(\cdot)\, x$ (respectively $\Cal {C} F(\cdot)(\la)x$) is holomorphic  on $\cc_+$ ($\cc\setminus i\,\r$) for
each $x\in X$, it follows that $\Cal {L}F$  (respectively $\Cal
{C}F$) is holomorphic on $\cc_+$ ($\cc\setminus i\,\r$), by  [19,
Theorem 3.10.1, p. 93].

If $\phi\in L^1(\jj,X)$, then  $\Cal {L} \phi$ has a continuous
extension  to $\overline{\cc_+}$  given also by (1.1). By the
Riemann-Lebesgue lemma,  if
  $g (\om)= \Cal {L} \phi(i\om)$, then $g\in C_0 (\r,X)$. Note that $g=\widehat{\phi|\r_+}$
   is the Fourier transform of
  $\phi|\r_+$ extended by $0$ to $\r$. Moreover, if $\phi\in \f
(\jj,X)$, then $g\in C^{\infty}_0 (\r,X)$. If $\phi \in
L^1(\r,X)$, then $\Cal {C} \phi $ has a continuous extension to
 $\cc$ and the Fourier transform $ \hat{\phi}= \Cal {C}\phi(i\cdot)$ satisfies $ \hat{\phi}\in C_0 (\r,X)$. If
$f\in \f (\r,X)$, then $\hat{f}\in \f (\r,X)$ by [32, p. 146] for
$X=\cc$. It follows that if $\phi\in L^{1}_{loc}(\r_+,X)\cap
\f'(\r,X)$, then
 $\widehat {\phi}\in  \f'(\r,X)$, by [32, p.
151] for $X=\cc$. So, if   $e_a (t)=e^{-at}$, then
$\a\phi(a+i\cdot)=\widehat{e_a\phi}  $ is  an $\f'(\r,X)$-valued
  function for all $a >0$.
Moreover, if $f\in \f(\r)$,
 $<\Cal {L} \phi (a+i\cdot),f>\,\,= <
\widehat{e_a\phi},f> \,\,=\,\, <{e_a\phi},\hat{f}>\,\,\to
<{\phi},\hat{f}>\,\,= <\widehat {\phi},f>$, where the limit exists
as $a\searrow 0$ by the Lebesgue dominating convergence theorem.
This means for  $\phi\in L^{1}_{loc}(\r_+,X)\cap \f'(\r,X)$,

(1.5)\qquad $\lim_{a\searrow 0}\a\phi (a+i\cdot)=\widehat {\phi} $
\qquad in \qquad $ \f'(\r,X)$.

For a  holomorphic function  $\zeta: \Sigma\to X$, where
$\Sigma=\cc_+$ or $\Sigma=\cc\setminus i\,\r$,  the point
$i\,\om\in i\,\r$ is called a $regular$ $point$ for $\zeta$ or
$\zeta$ is called $holomorphic$ at $i\,\om$, if $\zeta$ has an
extension $\overline{\zeta}$ which is holomorphic in a
neighbourhood $V\st \cc$ of $i\,\om$.

 Points $i\,\om$ which are $not$ $regular$ $points$ are
called $singular$ $points$.

 The $Laplace$ $spectrum$  of $\Phi\in L^1 _{loc}(\jj,X)\cap \f'(\r,X)$  or $\Phi\in  L_s^{\infty}(\jj,L(X))$
  is defined by

 (1.6) \qquad $sp^{\Cal{L}}
(\Phi) :=\{\om\in\r: i\, \om$ is a singular point for $\Cal
{L}\Phi \}$. See [4, p. 275].

\noindent The Laplace spectrum is called also the half-line
spectrum (see [4, p. 275]).

 The $Carleman \,\, spectrum$ of $\Phi\in  L^1 _{loc}(\r,X)\cap \f'(\r,X)$ or $\Phi\in  L_s^{\infty}(\r,L(X))$
 is defined by

(1.7) \qquad $sp^{\Cal{C}} (\Phi) :=\{\om\in\r: i\, \om$ is a
singular point for $\Cal {C} {\Phi}\}$. See [4, (4.26)].

We note that  for $\phi\in L^{\infty} (\r^+,X)$, $\om\not\in
sp^{\Cal{L}}(\phi)$ respectively $\phi\in L^{\infty} (\r,X)$,
$\om\not\in sp^{\Cal{C}}(\phi)$), one has

(1.8)\qquad $\lim _{\la\to 0}\Cal{L} (\gamma_{-\om }\phi)(\la)=
\overline{\Cal{L}} (\gamma_{-\om} \phi)(0)$  respectively

\qquad\qquad \,\,\,\, $\lim _{\la\to 0}\Cal{C} (\gamma_{-\om
}\phi)(\la)= \overline{\Cal{C}} (\gamma_{-\om} \phi)(0)$.

\noindent Indeed,  by the definitions that $0\not\in sp ^{\Cal{L}}
(\gamma_{-\om}\phi)$ or $0\not\in sp ^{\Cal{C}}
(\gamma_{-\om}\phi)$, it follows $\Cal{L} (\gamma_{-\om} \phi)$ or
$\Cal{C} (\gamma_{-\om} \phi)$ has a holomorphic extension $
\overline{\Cal{L}} (\gamma_{-\om} \phi)$ or $ \overline{\Cal{C}}
(\gamma_{-\om} \phi)$ in a neighbourhood of $0$.

For a  holomorphic function  $\zeta: \cc_+\to X$,  the point
$i\,\om\in i\,\r$ is called a  $weak$ $regular$ $point$ for
$\zeta$ if there exist $\e
> 0$ and  $h\in L^1 (\om-\e, \om+\e)$ such that

(1.9) \qquad $\lim_{a\searrow 0} \int _{-\infty}^{\infty}\zeta
(a+i\,s)\va (s)\, ds= \int _{\,\om-\e}^{\,\om+\e}h (s)\va (s)\,
ds$

\qquad \qquad\qquad \qquad\qquad for all $\va\in \h(\r) $ with
supp$ \,\va \st ]\om-\e, \om+\e[$.

 Points $i\,\om$ which are not weak regular points are
called $weak$ $singular$ $points$.

 The $weak\,\, Laplace\,\, spectrum$ of $\Phi\in
  L^1 _{loc}(\jj,X)\cap \f'(\r,X)$ or  $\Phi\in
L_s^{\infty}(\jj,L(X))$ is defined by (see [4, p. 324])

(1.10) \qquad $sp^{w\Cal{L}} (\Phi) :=\{\om\in\r: i\, \om$ is not
a weak regular point for $\Cal {L} {\Phi}\}$.

The definitions imply

(1.11) \qquad $sp^{w\Cal{L}} (\Phi)\st sp^{\Cal{L}} (\Phi)\st
sp^{\Cal{C}} (\Phi)$

\noindent  and   $sp^{w\Cal{L}} (\Phi)=\emptyset$ if $\Phi\in
L^1(\r,X)$.

 We conclude this section by recalling (see [7, p. 117]) that a subset
$\A\st L^{\infty}(\jj,X)$ is called $BUC$-$invariant$ if for each
$\psi\in BUC(\r,X)$ with $\psi|\,\jj\in \A$, one has
$\psi_a|\,\jj\in \A$ for each $a\in \r$.

In the following we assume

(1.12) \qquad $\A$ is  a $BUC$-$invariant$ (closed) subspace of
$L^{\infty}(\jj,X)$.

\noindent Examples of such $\A$ include

\qquad \qquad\qquad  $\{0\}$,  $C_0 (\jj,X)$, $AP(\r,X)$,
$AAP(\jj,X)$,

\noindent  the spaces consisting  respectively of the zero
function (when $\jj=\r$), all continuous vanishing at infinity,
almost periodic (when $\jj=\r$) and asymptotically almost periodic
functions. Such an $\A$ is called a $\Lambda$-class if
additionally $\A \st BUC(\jj,X)$, contains all constants and is
also closed under multiplication by characters.

 A point $\om\in \r$ is called $\A$-$regular$
for $ \phi\in L^1 _{loc}(\jj,X)\cap \f'(\r,X)$, where  $\A$ is a
class satisfying (1.12) if there is $f\in \f(\r)$ such that
$\hat{f}(\om)\not =0$ and $\phi*f|\,\jj \in \A$. The $reduced$
$Beurling$ $spectrum$ of $\phi$ with respect to $\A$ is defined by
(see [5, (4.1.1)], [6, (2.9)], [14] and references therein)

(1.13) \qquad $sp_{\A} (\phi):= \{\om\in \r: \om$ is not an
$\A$-regular point for $\phi \}$.

We note that if $ \phi\in L^{\infty}(\jj,X)$,  then a point
$\om\in \r$ is  $\A$-regular for $ \phi$   if and only if there is
$f\in L^1(\r)$ such that $\hat{f}(\om)\not =0$ and $\phi*f|\,\jj
\in \A$.

If $\A=\{0\}$, then $sp_{\{0\}} (\phi)$ is just the Beurling
spectrum defined in (1.17) below. We recall the following
property of $sp_{\A} (\phi)$. For $\phi\in
L^1_{loc}(\r,X)\cap\f'(\r,X)$, $f\in \f (\r)$,

(1.14) \qquad $sp_{\A} (\phi*f) \st sp_{\A} (\phi)\cap $ supp $
\hat{f} $.

For the  proof, see  the references in [8, Proposition 1.1 (ii)]
 and, for the case $\A=\{0\}$, [26, Proposition 06(i), p. 25].

We recall (see  [7], [9], [29], [30]) that a function $\phi\in
L^1_{loc} (\jj,X)$ is called $ergodic$  if there is a constant
$m(\phi)\in X$ such that

\qquad \qquad sup$_{x\in \jj}|| \frac{1}{T}\int_0^T \phi(t+s)\,ds
-m(\phi)||\to 0 $ as $T\to \infty$.

\noindent The limit $m(\phi)$ is called the $mean$ of $\phi$. The
set of all such ergodic functions will be denoted by $\E(\jj,X)$.
We set  $\E_b(\jj,X) =\E(\jj,X)\cap L^{\infty} (\jj,X)$,
$\E_0(\jj,X)=\{\phi\in \E_b(\jj,X): m(\phi)=0\} $, $\E_{ub}(\jj,X)
=\E(\jj,X)\cap BUC(\jj,X)$ and $\E_{u,0}(\jj,X)
=\E_{ub}(\jj,X)\cap \E_0 (\jj,X)$.

 Simple calculations show that for $\phi\in L^{\infty} (\jj,X)$,
 $t\in \jj$, $h, T >0$, one has

$\frac{1}{T}\int_0^T [\phi(t+s)-\phi(t+s+h)]
\,ds=\frac{1}{T}\int_0^h [\phi(t+s)-\phi(t+s+T)]\,ds$,

$\frac{1}{T}\int_0^T [M_h\phi(t+s)-\phi(t+s)]
\,ds=\frac{1}{hT}\int_0^h
\int_0^u[\phi(t+s)-\phi(t+s+T)]\,ds\,du$,

\noindent which implies

(1.15)\qquad  $\phi-\phi_h$ and $\phi-M_h\phi\in \E_0(\jj,X)$.

In the sequel we need the following analogue of Wiener's theorem
on  Fourier series. See [12, Proposition 1.1.5 (b), p. 22].

\proclaim{Proposition 1.1} Let $f\in L^1 (\r)$ and $\hat {f}
(\om)\not = 0$ for all $\om\in K$ a compact set. Then there is
$g\in L^1 (\r)$ such that $\hat {f} (\om) \cdot \hat {g} (\om)=1$
for all $\om\in K$.
\endproclaim

For $\phi \in L^{\infty}(\r,X)$, the corresponding  Arveson
spectrum,  Beurling spectrum and   maximal ideal (see [14], [27,
p. 184]) are defined respectively by

(1.16)\qquad $sp^A(\phi) =\{\la \in \r: $ for each
 $\e
>0$

 \qquad \qquad\qquad \qquad there exists $f\in L^1 (\r)$,  supp $
\hat {f} \st $
 $
]\la-\e,\la+\e[$ and $\, f* \phi\not = 0\}$

(1.17)\qquad $sp^B (\phi) = \{\la \in \r: $  if  $f\in L^1 (\r)$,
$ \hat {f}(\la) =1$ then $ f* \phi\not = 0\}$,

(1.18)\qquad $\,\,I(\phi)\,\,\,\, =\{f\in L^1(\r):\,  \phi*f
=0\}$.

The following result is well known and relates (1.17) to [4, p.
321], [5, Definition 4.1.2]. We include a proof for the benefit of
the reader.

 \proclaim{Proposition 1.2} For $\phi \in
L^{\infty}(\r,X)$, one has

(1.19)\qquad  $sp^B (\phi)= \{\la\in \r:
 \hat {f} (\la)=0 $ for all  $f\in I(\phi)\}= sp^A(\phi)$.
\endproclaim
\demo{Proof} The first equality  of (1.19) is easily verified.
Now,
 let $\la \not \in sp^B (\phi)$. There are $f_0\in L^1 (\r)$
with $\hat{f_0}(\la)=1$ but $f_0*\phi=0$; $\delta
>0$ such that $|\hat{f}_0(\mu)| \ge 1/2$ for $\mu\in
[\la-\delta,\la+\delta]$; and by Proposition 1.1,
 $g_0\in L^1 (\r)$ such that $\hat{g_0}(\mu)\cdot
\hat{f_0}(\mu)=1$ for $\mu\in [\la-\delta,\la+\delta]$. Let $f\in
L^1 (\r)$ and supp $\hat{f} \st [\la-\delta,\la+\delta]$. One has
$f*\phi =(f_0* g_0)*f*\phi =0$. This implies $\la \not \in sp^A
(\phi)$ and gives $sp^A (\phi)\st sp^B (\phi)$.  The converse is
proved similarly. \P
\enddemo

 \proclaim{Remark 1.3}  Assume $\A \st L^{\infty}(\jj,X)$
is a $BUC$-invariant subspace (see (1.12)).

 (i)  $C_0 (\jj,X) \st \A$ if $\jj=\r_+$ but this is not necessarily true if  $\jj=\r$.

(ii) If   $\jj=\r_+$, $\tilde{\phi}\in  L^1_{loc} (\r,X)$,
 $\tilde{\phi} |\r_+
 =\phi\in L^1_{loc} (\r_+,X)\cap \f'(\r,X)$ and  $\tilde{\phi}|\,(-\infty,0]\in L^{\infty} ((-\infty, 0],X)$,
  then $sp_{\A} (\tilde{\phi}) =sp_{\A} (\phi)$.

(iii) If $\A \st BUC(\jj,X)$  and $\phi\in L^{\infty} (\r,X)$ then
definition (1.13) is equivalent to  Definition 4.1.2 of [5]. In
particular, $sp_{\{0\}}(\phi)=sp^B (\phi)$.

(iv) If $\A= C_0(\jj,X)$ and $ \phi\in L^1(\r,X)$, then $sp_{\A}
(\phi)=\emptyset$.
\endproclaim

\demo{Proof} (i)  By the assumption if $\phi\in BUC(\r,X) $ and
$\phi$ has compact support from $(-\infty,0]$, then $\phi_t|\r_+
\in \A$ for all $t\in \r$.  It follows $C_c (\r_+,X)\st \A$ and so
$C_0 (\r_+,X)\st \A$  (see also the proof of Theorem 2.2.4 in [5,
p. 13]). A counter example for the case $\jj=\r$ is $\A=
AP(\r,X)$.

(ii) For $f\in \f(\r)$, one has $\tilde{\phi} * f |\r_+ (t)=
\int_0^{\infty} \phi (s) f(t-s)\, ds + \int_{-\infty}^0
\tilde{\phi} (s) f(t-s)\, ds =\int_0^{\infty} \phi (s) f(t-s)\, ds
+ \xi (t)$, where $\xi\in C_0 (\r_+,X)$.  This proves (ii).

(iii) We prove that  a point $\om_0\in \r$ is  $\A$-regular for $
\phi$
 if  there is $f_0\in L^1(\r)$ such
that $\widehat{f_0}(\om_0)\not =0$ and $\phi*f_0|\,\jj \in \A$,
since the converse is obvious. Choose
 $\delta
>0$ such that $\hat{f}_0(\om)\not = 0$ for $\om\in
[\om_0-\delta,\om_0+\delta]$ and by Proposition 1.1,
  $g_0\in L^1 (\r)$ such that $\hat{g_0}(\om)\cdot \hat{f_0}(\om)=1
$ for $\om\in [\om_0-\delta,\om_0+\delta]$.
 Let $f\in \f (\r)$, $\hat{f}(\om_0)\not =0$ and
supp $\hat{f} \st [\om_0-\delta,\om_0+\delta]$. Using Bochner
integration we have $ (\phi*f_0)*g_0|\,\jj\in \A$ and so
$\phi*f|\,\jj=\phi*(f_0*g_0*f)|\,\jj=(\phi*f_0*g_0)*f|\,\jj\in
\A$. This proves the  first part and the second part follows
taking $\A=\{0\}$.

 (iv) Since  $ \phi\in L^1(\r,X)$, it follows $\phi*f \in C_0 (\r,X)$
for all  $f\in C_0 (\r)$.  Since $\f(\r)\st C_0 (\r)$ the
definitions imply
 $sp_{C_0(\r,X)} (\phi)=\emptyset$ and $sp_{C_0(\r_+,X)}
 (\phi)=\emptyset$. \P
\enddemo

In Propositions 1.4, 1.5    $\psi$ will denote an element  of $ \f
(\r)$ with the following properties:

(1.20)\qquad $\widehat {\psi}$  has compact support,  $\widehat
{\psi} (0)=1$ and $\psi, \widehat {\psi} $ are non-negative.

\noindent An example of such $\psi$ is given by
$\psi=\widehat{\va}^2$, where $\va(t)=a\, e^{\frac{1}{t^2-1}}$,
$|t|\le 1$,  $\va=0$ elsewhere on $\r$ with $a$ some suitable
constant.

\proclaim{Proposition 1.4}
 The sequence $f_n(t)= n\,\psi (n\,t)$ is an approximate identity
for the space of uniformly continuous functions $UC(\r,X)$, that
is $||\phi*f_n-\phi||_{\infty}\to 0$ as $n\to \infty$.
\endproclaim

\demo {Proof}  Given  $\phi\in UC(\r,X)$ and $\e >0$ there exists
$k >0$ such that  $||\phi (t+s)-\phi(t)|| \le k  |s|+\e$ for all
$t,s\in \r$. But
 $\phi*f_n (t)- \phi(t)= \int_{-\infty}^{\infty} [\phi (t-\frac{s}{n})-\phi(t)]
 \psi (s)\,ds$ which gives $||\phi*f_n - \phi||_{\infty}\le (k/n)\int_{-\infty}^{\infty}
 |s|\psi (s)\,ds+ \e \int_{-\infty}^{\infty}
 \psi (s)\,ds$.  The result follows. \P
\enddemo

\proclaim{Proposition 1.5} Let $\phi\in L^{\infty}(\jj,X)$ or
$\phi\in UC(\jj,X)$.  Assume $0\not\in sp_{C_0 (\jj,X)}(\phi)$.
Then  $\phi$ is bounded and ergodic with mean $0$, in other words
$\phi\in \E_0 (\jj,X)$.
\endproclaim

\demo{Proof} By the Definition above (1.13),  $0\not\in sp_{C_0
(\jj,X)}(\phi)$ implies there is $f_0\in \f(\r)$ such that
$\widehat {f_0} (0)\not =0$ and ${\phi}* f_0|\,\jj \in C_0
(\jj,X)$.  Choose $\e
>0$  such that  $\widehat {f_0} (\om)\not = 0$ for all $\om \in [-2\e,2\e]$ and
$ sp_{C_0 (\jj,X)}(\phi)\cap  [-2\e,2\e]=\emptyset $. By
Proposition 1.1,  there is $g_0\in L^1 (\r)$ such that supp
$\widehat{g_0}\st [-2\e,2\e]$ and $\widehat{f_0}
(\om)\widehat{g_0}(\om)=1$ for $\om \in [-\e,\e]$. Set $f=
{f_0}*{g_0}$. If $\jj=\r_+$, let $\tilde{\phi}=\phi $ on $\r_+$
and $\tilde{\phi}=\phi(0) $ elsewhere on $\r$. It follows as in
the proof of Remark 1.1 (ii) that $\tilde{\phi}* f_0|\,\jj
 \in C_0 (\jj,X)$  and using Bochner integration  $\tilde{\phi}*
f|\,\jj  = (\tilde{\phi}* f_0)*g_0|\,\jj \in C_0 (\jj,X)$. Hence
$\tilde{\phi}* f|\,\jj \in \E_{u,0} (\jj,X)$. Set
$F=\tilde{\phi}-\tilde{\phi}*f $. Then $0\not\in sp_{\{0\}} (F)$.
It follows $F$ is bounded, by [8, Theorem 4.2]. Hence $PF$ is
bounded  by [8, Corollary 4.4] and so $F\in \E_0 (\r,X)$. This
implies $\phi = [F+ \tilde{\phi}* f]\,|\jj \in \E_{0}(\jj,X) $. \P
\enddemo

\proclaim{Remark 1.6} Proposition 1.5 is  true for $\phi \in L^1
_{loc}(\jj,X)$ with $\Delta_h^n \phi\in BUC(\jj,X)$ for all $h>0$
and  some $n\in \N$. Here $\Delta_h^n \phi
=\Delta_h(\Delta_h^{n-1} \phi)$ for $n >1$. However, it is not
valid for arbitrary $\phi\in L^1 _{loc}(\jj,X)\cap \f'(\r,X)$ as
the following example shows. The function $\phi (t)= te^{it}$ has
Beurling spectrum $sp^ {B}(\phi)=\{1\}$ and a
 direct calculation shows that $\phi$ is not ergodic.
\endproclaim

 \head{\S 2. Laplace and weak Laplace spectra}\endhead

In this section we  establish some tools which enable us to
calculate Laplace and weak Laplace spectra and relate them to the
reduced Beurling  spectrum relative to a class $\A$ satisfying
(1.12). We prove  new properties (Theorems 2.3(i), 2.4(i)) of the
weak Laplace spectrum which enable us to give simple proofs of
several tauberian results of Ingham [20], [22] and their
generalizations by Chill and
 others (see  [4, 4.10, p. 332] and references therein).

\proclaim{Proposition 2.1}  If $\Phi \in L^ 1_{loc}(\jj,X) \cap
\f'(\r,X)$ or $\Phi\in L_s^{\infty}(\jj,L(X))\}$, then

(i) $sp^{\Cal{L}} (\Phi) =sp^{\Cal{L}} (\Phi_a)$ for each $a\in
\jj$.

(ii) $sp^{\Cal{L}} (\Phi)=\cup_{h>0} sp^{\Cal{L}} (M_h \Phi)$.

(iii)  $sp^{\Cal{L}} (\gamma_{\om}\Phi)= \om + sp^{\Cal{L}}
(\Phi)$.

(iv) The statements (i), (ii), (iii) hold true for $sp^{w\Cal{L}}
$ and, when $\jj=\r$, for $sp^{\Cal{C}}$.
\endproclaim

\demo{Proof}  (i) A simple calculation shows for $\la\in
\cc^{\pm}$ and $\jj=\r$

(2.1)\qquad $\Cal {L^{\pm}}\Phi_a (\la) = e^{\la \,a}\Cal
{L^{\pm}}\Phi (\la)- e^{\la \,a} \int_0^{a}\, e^{-\la \,t }\Phi
(t)\,dt$.

\noindent  It follows $\Cal {L^+}\Phi$  (respectively $\Cal
{C}\Phi$) is holomorphic at $i\,\om$ if and only if $\Cal
{L^+}\Phi_a $ (respectively $\Cal {C}\Phi_a$) is holomorphic at
$i\,\om$. This proves (i)   for $sp^{\Cal{L}}$ and $sp^{\Cal{C}}$.

(ii) Another  calculation shows for  $\la \in \cc^{\pm}$ and
$\jj=\r$

(2.2)\qquad $\Cal {L^{\pm}}(M_h \Phi) (\la) =$
  $g(\la\, h) \Cal {L^{\pm}}\Phi (\la) - (1/h)\int_0^h (e^{\la \, v} \int_0 ^u e^{-\la \, t}\Phi(t+v)
 dt)\,dv$,

\noindent where  $g$ is the entire function given by
$g(\la)=\frac{e^{\la}-1}{\la}$ for $\la\not = 0$. Let $i\,\om\in
i\, \r$ be a regular point for $\Cal {L^+}\Phi$ and let $\overline
{\Cal {L^+}}\Phi: V \to X$ be a holomorphic extension of $\Cal
{L^+}\Phi$ to a neighbourhood $V \st \cc$ of $i\,\om$. Then
$\overline {\Cal {L^+}}(M_h \Phi) (\la) =$
  $g(\la\,h)\overline  {\Cal {L^+}}\Phi (\la) -
  (1/h)\int_0^h (e^{\la \, v} \int_0 ^u e^{-\la \, t}\Phi(t+v)\,
 dt) \,dv$, $\la \in V$  is a holomorphic extension of
$\Cal {L^+}(M_h \Phi)$. So $i\,\om$ is a regular point for $\Cal
{L^+}( M_h \Phi)$.
 If $\om\in \r$ there is  $h_0 >0$ such
that $g(i\om\, h_0)\not =0$. Similarly as above, if $i\,\om\in i\,
\r$ is a regular point for $\Cal {L^+} (M_{h_0}\Phi)$, then
$i\,\om$ is a regular point for $\Cal {L^+} \Phi$. This proves the
first part of (ii). The second part follows similarly noting that
(2.2) implies $\Cal {C}(M_h \Phi) (\la) =$
  $g(\la\, h) \Cal {C}\Phi (\la) - (1/h)\int_0^h (e^{\la \, v} \int_0 ^u e^{-\la \, t}\Phi(t+v)
 dt)\,dv$. This proves (ii)   for $sp^{\Cal{L}}$ and $sp^{\Cal{C}}$.

(iii) This follows easily from the definitions noting that $\Cal
{L}(\gamma_{\om}\Phi) (\la)= \Cal {L}\Phi (\la -i\om)$  and $\Cal
{C}(\gamma_{\om}\Phi) (\la)= \Cal {C}\Phi (\la -i\om)$. This
proves (iii) for $sp^{\Cal{L}}$, $sp^{\Cal{C}}$ and
 $sp^{w\Cal{L}}$.
 \P

(iv) The proofs of (i) and (ii) for  $sp^{w\Cal{L}}$ follow
similarly as in the case $sp^{\Cal{L}}$ using [28, Theorem 6.18,
p.146] as in the proof of Proposition 2.4(i) below.

\enddemo

\proclaim{Example 2.2} Let $\phi (t)=e^{it^2}$, $t\in \r$. Then
$sp^{\Cal {L}}(\phi) =\emptyset$ and $sp^{\Cal {C}}(\phi) =\r$.
Moreover, $M_h \phi \in  C_0 (\r,\cc)$ and $sp^{\Cal {L}}(M_h\phi)
=\emptyset$ for $h>0$.
\endproclaim
\demo{Proof} By Proposition 1.1 (i), (iii), it is readily verified
 that $sp^{\Cal {L}}(\phi_a)=sp^{\Cal {L}}(\phi)=2a+ sp^{\Cal
 {L}}(\phi)$ for each  $a\in \r$. This implies that either $sp^{\Cal {L}}(\phi)
 =\emptyset$ or $sp^{\Cal {L}}(\phi) =\r$. We claim
 that $ sp^{\Cal {L}}(\phi)=\emptyset $. Indeed, let $ y (\la)= \Cal
 {L}\phi (\la)$ for $\la\in \cc_+$. Then $y' (\la)+ (\la/2i)y (\la)= 1/2i$.  Solving this equation, one gets
  $\psi (\la,a) = e^{-(\la ^2/4i)} (a+ (1/2i) P e^{(\la
 ^2/4i)})$ is a general solution, where $a\in \cc$. One has
 $\Cal{L}\phi (0)= ((1+i)\pi^{1/2})/2^{3/2}=a_0$ and so $ \Cal {L}\phi (\la)= \psi (\la,a_0)$ is a particular solution.
  Since $ e^{(\la^2/4i)}$ is an entire function, $\psi (\la,a_0)$ is an entire extension  of  $\Cal
 {L}\phi$, implying
  $sp^{\Cal {L}}(\phi)
 =\emptyset$.  A similar argument  shows that either $sp^{\Cal {C}}(\phi)
 =\emptyset$ or $sp^{\Cal {C}}(\phi) =\r$. Since  $sp^{\Cal {B}}(\phi)=sp^{\Cal {C}}(\phi)
 \not=\emptyset$, one gets  $sp^{\Cal {C}}(\phi)
 =\r$. Since $\int_0^{\infty} e^{i\,t^2}\, dt$ is
an improper Riemann integral, it follows $P\phi (T)=\int_0^T
e^{i\,t^2}\, dt \to a_0$ as $T\to\infty$. Since  $M_h \phi(t)=
P\phi (t+h)-P\phi (t)$, $M_h\phi\in C_0 (\r,\cc)$ for $h>0$.
Finally, by Proposition 2.1(ii) and   $sp^{\Cal {L}}(\phi)
 =\emptyset$, one concludes $sp^{\Cal {L}}(M_h\phi)
=\emptyset$ for $h>0$.
  \P
\enddemo

\proclaim{Theorem 2.3 (Ingham)} Let $\phi \in
L^{1}_{loc}(\r_+,X)\cap \f'(\r,X)$ and
  $sp^{w\Cal{L}}(\phi)=\emptyset$.

(i)\,\,\, $\phi*g \in C_0(\r,X)$ for all $g\in \f(\r)$ with
$\hat{g}\in \h(\r)$.

(ii)\,\,\,If $\phi \in L^{\infty}(\r_+,X)$, then $\phi*g \in
C_0(\r,X)$ for all $g\in L^1(\r)$.

(iii) \,\,\,If $\phi \in BUC(\r_+,X)$ or more generally  $\phi$ is
slowly oscillating, then $\phi \in C_0(\r_+,X)$.

(iv) \,\,\,If $sp^{\Cal{L}}(\phi)=\emptyset$, then
$(P\gamma_{-\om}\phi-\overline{\Cal{L}}\phi(i\om))\in C_0(\r_+,X)$
 for all $\om\in \r$.
\endproclaim
\demo{Proof} (i) As  $\phi \in \f' (\r,X)$,   the Fourier
transform $\widehat{\phi} \in \f' (\r,X)$ is given by

(2.3) \qquad $<\widehat{\phi},f>= <\phi, \hat
{f})>=\int_0^{\infty} \phi (s)\hat{ f}(s)\, ds$, \,\,  $f\in \f
(\r)$.

The condition $sp^{w\Cal{L}}(\phi)=\emptyset$ means that
$\widehat{\phi}|\, \h(\r) =h$, where $h$ is locally integrable on
$\r$ (see for example [4, Lemma 4.9.3]) and so

(2.4)\qquad $\int_0^{\infty} \phi (s)\hat{ f}(s)\, ds=
\int_{-\infty}^{\infty} h(\eta) f(\eta)\, d\eta$,\,\,\, $f\in
\h(\r)$.

 Take $g\in \f(\r)$ with $\, \hat{g}\in \h(\r)$.
 Using (2.3) and (2.4), one gets

(2.5)\qquad $\phi*g(t)= \int_0^{\infty} \phi (s) g(t-s)\, ds=
(1/2\pi)\int_{-\infty}^{\infty} h(\eta) e^{it\eta}\hat {g}(\eta)\,
d\eta$.

Now, $h\hat{g} \in L^1 (\r,X)$ and so by the Riemann-Lebesgue
lemma (see [4, p. 45] or [21, p. 123]), $\phi*g\in C_0(\r,X)$.
This proves (i).

(ii) Given  $g\in L^1(\r)$ choose $(g_n) \st \f(\r)$ with
$(\widehat{g_n})\st \h(\r)\}$ and $||g_n-g||_{L^1}\to 0 $ as
$n\to\infty$. Since   $\phi \in L^{\infty}(\r_+,X)$, we conclude
 $||\phi*
g_n- \phi*g||_{{\infty}}\to 0 $ as $n\to\infty$ and so by (i),
$\phi* g\in C_0 (\r,X)$). This proves (ii).

(iii) Consider first the case $\phi\in BUC(\r_+,X)$ and let
$\om\in\r$, $\e>0$. Take $g\in \f(\r)$ with supp ${\hat{g}} \st
[\om-\e,\om+\e]$ and ${\hat{g}}(\om)\not =0$. By (i), $\phi*g \in
C_0(\r,X)$.  By Remark 1.1 (ii),  we conclude $\om \not\in
sp_{C_0(\r_+,X)}(\tilde{\phi})$, where $\tilde{\phi}\in BUC(\r,X)$
is any extension of $\phi$. This implies
$sp_{C_0(\r_+,X)}(\tilde{\phi})=\emptyset$ and so
$\phi=\tilde{\phi}|\r_+\in C_0(\r_+,X)$ by [5, Theorem 4.2.1]. For
the  general case, it follows $sp_{C_0(\r_+,X)}(M_h
\phi)=\emptyset$ for all $h>0$ and hence $M_h \phi\in C_0(\r_+,X)$
for all $h>0 $. This completes the proof of (iii), by [4, Theorem
3.2.3, p. 250].

(iv) Replacing $\phi$ by $\gamma_{-\om}\phi$, we may assume
$\om=0$. Set $\Psi (t) = P\phi (t)-\overline{\Cal{L}}\phi(0)$ for
$t\in\r_+$ and $\Psi= -\overline{\Cal{L}}\phi(0)$ elsewhere on
$\r$. Simple calculation shows
$\Cal{L}\Psi(\la)=(\Cal{L}\phi(\la)-\overline{\Cal{L}}\phi(0))/\la$
for  $\la\in \cc_+$. Since $\overline{\Cal{L}}\phi$ is holomorphic
on $\overline{\cc_+}$, it follows that $\Cal{L}\Psi$ has a
holomorphic extension to an open neighbourhood of
$\overline{\cc_+}$. This implies $ sp^{\Cal{L}} (\Psi)=\emptyset$.
So, by  (1.11) and part (i), one gets $\Psi*g|\,\r_+ \in
C_0(\r_+,X)$ for all $g\in \f(\r)$ with
 $\, \hat {g}\in \h(\r)$.
  Consider the sequence $f_n (t)$ of Proposition 1.3.
 Since $\Psi\in UC(\r,X)$  and $\Psi*f_n|\r_+\in C_0 (\r_+,X)$ for each $n\in \N$, one concludes
 $\Psi|\,\r_+\in C_0 (\r_+,X)$ by Proposition 1.3. This finishes the proof.
  \P
\enddemo

 \proclaim{Theorem 2.4}  Let $\phi \in L^1_{loc}
 (\r_+,X)\cap \f'(\r,X)$. Then

(i) $sp^{w\Cal{L}}(\phi *f) \st sp^{w\Cal{L}}(\phi)\cap\text
{\,supp\,} \hat {f}$ \qquad   for all $f\in \f (\r)$.

(ii) $sp_{C_0(\r_+,X)}(\phi) \st sp^{w\Cal{L}}(\phi)$.

(iii) If $\phi \in L^{\infty}(\r_+,X)$ and  $\om_0\not \in
sp^{w\Cal{L}}(\phi)$, then $\gamma_{-\om_0}\phi \in \E_0(\r_+,X)$.

(iv) If $\phi \in L^{\infty}(\r_+,X)$ and  $\om_0\not \in
sp^{\Cal{L}}(\phi)$, then $ (P\gamma_{-\om_0}\phi-\overline{\Cal
{L}}\phi(i\om_0))\in C_0(\r_+,X)$.
\endproclaim
 \demo{Proof} (i) Direct calculations show that for $\la\in
\cc$

\qquad $\Cal{L}{\phi}*f (\la)=\Cal{L}\phi(\la)\cdot \Cal{L}f
(\la)+ \zeta (\la)$ with $\zeta (\la)=\int_0^{\infty}\Cal{L}\phi_s
(\la) f(-s)\, ds$.

Since $f\in \f (\r)$, $\lim_{a\searrow 0}\Cal{L}f
(a+i\,\eta)=\int_0^{\infty} e^{-i\,\eta\, t} {f}(t)\,dt :=g(\eta)$
for each $\eta \in \r$.  Since  $g\in C^{\infty} (\r)$ we conclude
$sp^{w\Cal{L}}(f)=\emptyset$. Now, assume $i\,\om$ is a weak
regular point for $\Cal{L}\phi$. So there exists  $\e >0$ and $h
\in L^1 (\om-\e,\om+\e)$ satisfying $\lim_{a\searrow 0} \int
_{-\infty}^{\infty}\Cal{L}\phi (a+i\,\eta)\va (\eta)\, d\eta= \int
_{i\,\om-\e}^{i\,\om+\e}h (s)\va (\eta)\, d\eta$
 for all $\va\in \h(\r) $ with supp$ \,\va \st
]\om-\e, \om+\e[$. Then by [28, Theorem 6.18, p. 146],
$\lim_{a\searrow 0} \int _{-\infty}^{\infty}\Cal{L}\phi
(a+i\,\eta)\Cal{L} f (a+i\,\eta)\,\va (\eta) d\eta= \int
_{\om-\e}^{\om+\e}h (\eta) g (\eta)\,\va (\eta) d\eta$
 for all $\va\in \h(\r) $ with supp$ \,\va \st ]\om-\e, \om+\e[$.
 It follows $i\,\om$ is a weak regular point for
$\Cal{L}\phi\cdot \Cal{L}f$. Since  $\Cal{L}\phi_s(a+i\eta)=
e^{(a+i\eta)s}[\Cal{L}\phi (a+i\eta)- \int_0^s
e^{-(a+i\eta)t}\phi(t)\, dt] $,  the same argument shows $\,\om$
is a weak regular point for $\zeta$. It follows $\,\om$ is a weak
regular point for $\Cal{L} ({\phi}*f)$ showing
$sp^{w\Cal{L}}({\phi}*f) \st sp^{w\Cal{L}}(\phi)$. By (1.11),
(1.14), one gets $sp^{w\Cal{L}}({\phi}*f)\st
sp^{\Cal{C}}({\phi}*f) \st \text {supp}\hat{f}$. This finishes the
proof of (i).

(ii) Let $\om \not \in  sp^{w\Cal{L}}(\phi) $. Choose $\e >0$ and
$f\in\f(\r)$ such that $sp^{w\Cal{L}}(\phi)\cap
[\om-\e,\om+\e]=\emptyset$, $\hat{f}(\om)=1$ and supp $\hat{f}\st
[\om-\e,\om+\e]$. By part (i), it follows $
sp^{w\Cal{L}}(\phi*f)=\emptyset$ and so by Theorem 2.3 (ii),
$\phi*f|\r_+\in C_0 (\r_+,X)$. This implies $\om\not\in
sp_{C_0(\r_+,X)}(\phi)$ and proves (ii).

(iii)  Replacing $\phi$ by $\gamma_{-\om_0}\phi$,  we may assume
$\om_0 =0$. By part (ii), one concludes $0\not\in
sp_{C_0(\r_+,X)}(\phi)$. The statement follows by Proposition 1.4.

(iv)  As in the proof of Theorem 2.3 (iv),  we may assume  $\om_0
=0$ and   conclude $0\not\in sp^{\Cal{L}}(P\phi-\overline{\Cal
{L}}\phi(0))= sp^{\Cal{L}}(\phi)$. Choose $\e
>0$  and $f\in \f(\r)$ such that $ sp^{\Cal{L}}(\phi) \cap [-\e,\e]
=\emptyset$, supp $\hat{f} \st [-\e,\e]$ and $\hat{f} (\om)=1$ if
$|\om|\le \e/2$. One has $(P\phi*f- \overline {\Cal{L}} \phi(0))=
(P\phi- \overline {\Cal{L}}\phi (0))*f$ and so  by Part (i), it
follows $sp ^{w\Cal {L}} (P\phi*f- \overline {\Cal{L}}
\phi(0))=\emptyset$. Choose $g\in \f (\r)$  with $\hat {g}\in \h
(\r)$ and $\hat {g}(\om)=1$ if $|\om|\le \e$.  Then  $(P\phi *f-
\overline {\Cal{L}}\phi (0))= (P\phi *f- \overline {\Cal{L}}\phi
(0))*g \in C_0 (\r_+,X)$, by Theorem 2.3(i). It is easily verified
 that $0\not\in sp^{\Cal{C}}(\phi-\phi*f)$ and hence $P(\phi-\phi*f)\in
BUC(\r,X)$, by [8, Corollary 4.4]. It follows $P\phi=
(P(\phi-\phi*f)+ P\phi*f)|\, \r_+\in BUC(\r_+,X)$. \P
\enddemo
We conclude this section giving a short proof for [4, Theorem
4.9.7].

\proclaim{Theorem 2.5 } Let $\phi \in BUC(\r_+,X)$ or more
generally $\phi \in L^{\infty}(\r_+,X)$ with $\phi$  slowly
oscillating. If $sp^{w\Cal{L}}(\phi)$ is countable and
$\gamma_{-\om}\phi\in \E(\r_+,X)$ for all $\om\in
sp^{w\Cal{L}}(\phi)$, then $\phi \in AAP(\r_+,X)$.
\endproclaim
\demo{Proof} Consider first the case $\phi \in BUC(\r_+,X)$. Let
$\tilde{\phi}\in BUC(\r,X)$ be an extension of $\phi$. Using $C_0
(\r_+,X)\st AAP(\r_+,X)$, (1.11) and Theorem 2.4 (ii), one gets
$sp_{AAP(\r_+,X)}(\tilde{\phi})\st sp_{C_0
(\r_+,X)}(\tilde{\phi})\st
sp^{w\Cal{L}}(\tilde{\phi})=sp^{w\Cal{L}}({\phi})$. It follows,
$sp_{AAP(\r_+,X)}(\tilde{\phi})$ is countable. This and the
assumptions imply  $\phi =\tilde{\phi}|\r_+ \in AAP(\r_+,X)$, by
[5,Theorem 4.2.6]. The case $\phi$ is slowly oscillating can be
proved similarly to  Theorem 2.3(ii) noting that direct
calculation shows that if $\gamma_{-\om}\phi\in \E(\r_+,X)$, then
$\gamma_{-\om}M_h\phi\in \E(\r_+,X)$ for all $h>0$. \P
\enddemo

\head{\S 3. Bounded $C_0$}-semigroups (groups)\endhead

For a bounded $C_0$-group $T(t)\in L(X)$, $t\in \r $, we make the
following definitions recalling that $T(\cdot)$ is strongly
measurable  but not necessarily measurable:

(3.1)\qquad $sp^B (T(\cdot)):=$
  $\{\la \in \r: $ there is
 $y\in X$ such that if $f\in L^1
(\r)$

\qquad\qquad\qquad\qquad \qquad\qquad\qquad\qquad  and
 $
\hat {f}(\la) =1$ then  $ f* T(\cdot)y\not = 0\}$.

The Arveson spectrum of $T(\cdot)$ is defined ([2, p. 365], [14,
Definition 4]) by

(3.2)\qquad $sp^A (T(\cdot))=\{\la \in \r: $ for each $\e >0$
there is $y\in X$ and

\qquad\qquad \qquad\qquad \qquad$f\in L^1 (\r)$ with supp $ \hat
{f} \st ]\la-\e,\la+\e[$ and $\, f* T(\cdot)y\not = 0\}$.

Using (1.17),  it is easily verified that

(3.3) \qquad $sp^B (T(\cdot))= \cup _{x\in X} sp^B
 (T(\cdot)x)$.

\proclaim{Proposition 3.1}  Let $ T(t)\in L(X)$, $t\in\r_+$ be a
bounded $C_0$-semigroup with  generator $A$. Then

(3.4) \qquad $ i\,sp^{\Cal{L}} (T(\cdot))=  \,\sigma (A)\cap
i\,\r=i\, \cup_{x\in X} \, sp^{\Cal{L}} (T(\cdot)x)$.
\endproclaim
\demo{Proof}  By [25, (3.5), (3.7), p. 9], one has
$(\la-A)^{-1}x=\int_0^{\infty }\, e^{-\la \,t} T(t)x\,dt$ for $\la
\in \cc_+$, $x\in X$. It follows by (1.1) and (1.2) that  $
(\la-A)^{-1}=\Cal{L}T(\cdot)$.  Since $R(\la):=(\la-A)^{-1} $ is
holomorphic on $\rho (A)=\cc\setminus \sigma (A)$ (see [25,
(5.21), p. 20]), it follows that if $i\,\om \in \rho (A)\cap
i\,\r$ then $\om\not\in sp^{\Cal{L}} (T(\cdot))$ proving\,\,\,\,
 $i\,sp^{\Cal{L}} (T(\cdot))\st  \sigma (A)\cap i\,\r$.
 The definitions imply
 $\cup_{x\in X} sp^{\Cal{L}} (T(\cdot)x)\st
\,sp^{\Cal{L}} (T(\cdot))$. So,
 it remains to show $ \sigma (A)\cap i\,\r \st i\,\cup_{x\in X} sp^{\Cal{L}}
(T(\cdot)x)$.
  Assume  $\om  \in \r\setminus \cup_{x\in X} sp^{\Cal{L}}
(T(\cdot)x)$. Then for each $x\in X$, there is an open disk
$V_x\st \cc$ with center $i\,\om$  and  a holomorphic function
$F_x: V_x\to X$  such that $ F_x(\la)=\Cal {L}T(\cdot)x (\la)=
R(\la)x$ for $\la \in V_x \cap \cc_+$. It follows  $(\la
I-A)F_x(\la)=x$,\,\, for \,\, \,\,  $\la \in V_x \cap \cc_+$, \,\,
$x\in X$ and $\la F_x (\la)-F_{ Ax} (\la)=x$,\,\, for \,\,$\la \in
V_x \cap V_{Ax} \cap \cc_+$, \,\,$x\in D(A)$.
 Since $F_x$ is continuous at $i\,\om$ for each $x\in X$ and $A$
is a closed operator, the identities remain valid for
$\la=i\,\om$. The identity $(i\,\om I-A)F_x(i\,\om)=x$ implies
$i\,\om I-A$ is onto. Now, assume $x_0 \in D(A)$ with $(i\,\om
I-A)x_0=0$. Then $i\,\om F_{x_0} (i\,\om ) - F_{Ax_0} (i\, \om)=
x_0=F_{0} (i\, \om)=0$. This gives $x_0=0$ and proves $i\,\om I-A$
is one to one. By the closed graph theorem $i\,\om \in \rho (A)$.
This  completes the proof. \P
\enddemo

  We note that if $g\in L^1 (\r)$ with $\hat
{g}(\la) =1$,
 supp $\hat {g}\st  [\la-\delta,\la+\delta]$, $ T(t)\in L(X)$, $t\in\r$ is a bounded $C_0$-group,
   $x\in X$ and  $y=\int_{\infty}^{\infty} g(s) T(-s)x\, ds $,
 then

(3.5) \qquad $T(\cdot)y=g*T(\cdot)x$ and $sp^B (T(\cdot)y) \st
[\la-\delta,\la+\delta]$.

For $y\in X$ set $X_y := \overline{\text{span}}\{z\in X: z=T(t)y,
t\in \r\}$, $T_y (t)z= T(t)z$ for all $z\in X_y$ and $A_y$ the
generator of $T_y (\cdot)$.

\proclaim{Proposition 3.2}  Let $ T(t)\in L(X)$, $t\in\r$ be a
bounded $C_0$-group and  $x\in X$.

(3.6)\qquad
 $i sp^B (T(\cdot)x) =i sp^{\Cal{C}}
(T(\cdot)x)= i sp^{\Cal{L}} (T(\cdot)x)=\sigma (A_x) $,

(3.7)\qquad  $ i sp^A (T(\cdot))= i sp^B (T(\cdot))= i
sp^{\Cal{C}} (T(\cdot))= i sp^{\Cal{L}} (T(\cdot))= \sigma (A)$.

\noindent Furthermore each of these sets is  closed.
\endproclaim
\demo{Proof} By the first equality of(1.19) for $\phi= T(\cdot)x$,
one gets $sp^B (T(\cdot))x$ is closed for each $x\in X$.
 We prove $sp^B (T(\cdot))$ is closed. Let $(\la_n)\st sp^B
(T(\cdot))$  and $\la_n \to \la$ as $n\to \infty$. We restrict
ourself to the case $\la_n >  \la_{n+1}$ for $n\in \N$. Choose
$x_n \in X$ such that $||x_n||=1$ and $f* T(\cdot)x_n\not =0$ for
each $f\in L^1 (\r)$ with $\hat{f}(\la_n)=1$.  By (3.5),  one can
replace $(x_n)$ by $(y_n)$ with $sp^B (T(\cdot)y_n)\st I_n
=[\la_n-\delta_n,\la_n+\delta_n]$, where $ 0 <\delta_n  < (\la_n
-\la_{n+1})/2$ and  $\delta_n > \delta_{n+1}$ for $n\in \N$. Set
$y=\sum_{k=1}^{\infty} y_k/2^k$.   Because the intervals $I_n$ are
compact and disjoint, it follows that choosing $h_n\in L^1 (\r)$
with  $\hat{h}_n=1$ on  $I_n$ and  $\hat{h}_n=0$ on a
neighbourhood of  $\cup_{k\not =n} I_k$ one has $h_n *T(\cdot)y =
(T(\cdot)y_n)/2^n$. Hence $sp^B (T(\cdot)y_n)\st sp^B (T(\cdot)y)$
for each $n\in \N$. Let $f\in L^1 (\r)$ and $\hat{f}(\la)=1$.
Since $\la_{n}\to \la$ as $n\to \infty$ and $\hat{f}$ is
continuous, there is $n (f)$ such that $\hat{f}(\la_{n})\not =0$
for $ n \ge n (f)$. It follows $f*T(\cdot)y_{n}\not = 0$ for $ n
\ge n (f)$. This implies $f*T(\cdot)y\not = 0$ and  proves $\la
\in sp^B (T(\cdot))$ and so $sp^B (T(\cdot))$ is closed.

By  (3.3)  and (1.19) we get

 (3.8) \qquad $ sp^B (T(\cdot))= \cup _{x\in X} sp^B (T(\cdot)x)= \cup _{x\in X}
sp^A (T(\cdot)x)$.

We claim $sp^A (T(\cdot))=  \cup _{x\in X} sp ^A (T(\cdot)x)$.
 Indeed, the definition gives $ sp^A (T(\cdot)x) $

 \noindent $\st sp^A
(T(\cdot))$ for each $x\in X$ and so $ \cup _{x\in X} sp^A
(T(\cdot)x)\st sp^A (T(\cdot))$. As $sp^B (T(\cdot))=\cup _{x\in
X} sp^A (T(\cdot)x)$ is closed,  for each $\la \not \in \cup
_{x\in X} sp^A (T(\cdot)x)$ there is $\delta
>0$ such that $[\la-\delta,\la+\delta]\cap  (\cup _{x\in X} sp^A
(T(\cdot)x))=\emptyset$. Let  $f\in L^1 (\r)$ and supp $\hat{f}\st
[\la-\delta,\la+\delta]$. Then $sp^B (f*T(\cdot)x)=\emptyset$  so
$f*T(\cdot)x=0$ for each $x\in X$. This implies $\la \not \in sp^A
(T(\cdot))$, $sp^A (T(\cdot))= \cup _{x\in X} sp ^A(T(\cdot)x)$
and
  proves
the first  identity of (3.7).

As in the proof of Proposition 3.1, $T(\cdot)\in L_s^{\infty} (\r,
L(X))$ and $\Cal {C}T(\cdot) (\la)=(\la-A)^{-1}$ for $\la \in
\cc\setminus i\,\r$. Moreover,  if $i\,\om$ is a regular point for
$\Cal{L}T(\cdot)$ then  by Proposition 3.1 $i\,\om\in \rho (A)$.
So, $i\,\om$ is a regular point for $(\la-A)^{-1}$. It follows
$i\, \om$ is a singular point for $\Cal {C} T(\cdot)$ if and only
if $i\, \om$ is a singular point for $\Cal {L}T (\cdot) $. This
and Proposition 3.1 imply $ sp^{\Cal{C}} (T(\cdot))= sp^{\Cal{L}}
(T(\cdot))= \sigma(A)\cap i\r$. Noting that $sp^{\Cal{L}}
(T(\cdot)z)\st sp^{\Cal{L}} (T(\cdot)x)$ for each $z\in X_x$ and
using Proposition 3.1 one gets $ \sigma(A_x)\cap i\r=i
sp^{\Cal{L}} (T_x(\cdot))=i sp^{\Cal{L}} (T(\cdot)x) $. Since
$\cup_{y\in X} sp^{\Cal{C}} (T(\cdot)y)\st sp^{\Cal{C}}
(T(\cdot))=sp^{\Cal{L}} (T(\cdot))= \cup_{y\in X} sp^{\Cal{L}}
(T(\cdot)y)\st  \cup_{y\in X} sp^{\Cal{C}} (T(\cdot)y)$, one gets
$sp^{\Cal{C}} (T(\cdot))= \cup_{y\in X} sp^{\Cal{C}} (T(\cdot)y)
$. But  $sp^{\Cal{C}} (T(\cdot)z)\st sp^{\Cal{C}} (T(\cdot)x)$ for
each $z\in X_x$. So, by the above, we conclude $ sp^{\Cal{C}}
(T(\cdot)x)= sp^{\Cal{C}} (T_x(\cdot))= sp^{\Cal{L}} (T_x(\cdot))=
sp^{\Cal{L}} (T(\cdot)x)$ for each $x\in X$.
 By
(1.19), Proposition 1.2 and [4, Theorem 4.8.4] $sp^B (T(\cdot)x)=
sp^{\Cal{C}} (T(\cdot)x)$ for each $x\in X$. Since $T(\cdot)$ is
bounded, one gets $\sigma (A)\st i\r$ and  $\sigma (A_x)\st i\r$
for each $x\in X$.
 So,
 (3.6) and (3.7) follow by the above, (3.4) and (3.8). \P
\enddemo

\proclaim {Remark 3.3} (i)   Definition  (1.16) of $ sp^A
(T(\cdot)x)$ is easily seen to be equivalent  to the the
definition of the  Arveson spectrum of $x$ with respect to
$T(\cdot)$ in [2, p. 365], [3],  [14, Definition 4].

(ii) In [14] a proof that $sp^A (T(\cdot))= \overline{\cup _{x\in
X} sp^A (T(\cdot)x)}$ is outlined. So, the result $sp^A
(T(\cdot))=\cup _{x\in X} sp^A (T(\cdot)x)$ seems new.

(iii) The Arveson spectrum defined in [17, p. 285], [18] is easily
seen to be $
 i\, \{\la: \hat{f}(\la)=0$ for all $f \cap_{x\in X} I
(T(\cdot)x)\}=   \sigma (A)= i\,sp^A (T(\cdot))$, by [17, p. 285]
or [18] and  Propositions 3.1, 3.2.  This  gives a proof of the
remark above Theorem 5 in [14].
\endproclaim

 \head{\S 4. Uniform  Laplace and Carleman spectra}\endhead

In this section we recall the definition  of uniform spectrum
$sp^{\Cal{L}_u} (\phi)$ for functions from $ BUC(\jj,X)$ and
extend it to functions from $L^{\infty}(\jj,X)$.  For $ \phi \in
BUC(\jj,X)$, the Carleman spectrum coincides with Arveson spectrum
of the generator $A_{\phi}$ of the group of translations
$S(\cdot)$ restricted to the subspace $L(\phi)=\overline
{span}\{\phi_t: t\in \r\}$. In the case $\phi\in BC(\r,X)$ it is
proved  in [15, Proposition 2.3 (iii)] that $sp^{\Cal{L}}
(\phi)\st sp^{\Cal{L}_u} (\phi)$ and the inclusion may be strict.
Here, we prove that $sp^{\Cal{L}_u} (\phi)=sp^{\Cal{L}} (\phi)$
for $\phi\in L^{\infty}(\jj,X)$ (Theorem 4.2).

This result is new and our proof seems new even for the case
$\phi\in BUC(\jj,X)$.

Using (1.1)  for $\phi \in L^{\infty}(\jj,X)$   and  (1.3) when
$\jj=\r$, define

(4.1) \qquad $ \Cal{L}_u \phi(\la) (s)=\Cal {L} {\phi_s} (\la)$
\qquad for $s\in \jj$, $\la\in \cc_+$;

(4.2) \qquad   $ \Cal{C}_u\phi (\la)(s)=\Cal {C} {\phi_s} (\la)$
\qquad for $s\in \r$, $\la \in \cc\setminus i\,\r$.

\proclaim{Lemma 4.1} (i) For  $ \Cal{L}_u \phi$, one has $
\Cal{L}_u \phi(\la) \in BUC(\jj,X)$ and $ \Cal{L}_u \phi(\cdot)$
is holomorphic  on $\cc_+$. If $\phi \in BUC(\jj,X)$, then $
\Cal{L}_u \phi = \Cal {L} S^{\jj}(\cdot)\phi$.

(ii) For $ \Cal{C}_u \phi$, one has $ \Cal{C}_u \phi(\la)\in
BUC(\r,X)$ and $ \Cal{C}_u \phi(\cdot)$ is holomorphic on
$\cc\setminus i\r$. If $\phi \in BUC(\r,X)$, then $ \Cal{C}_u \phi
= \Cal {C} S^{\r}(\cdot)\phi$.

(iii) If $\jj=\r^+$, $\om_0\not \in sp^{\Cal{L}} (\phi)$ (or
$\jj=\r$, $\om_0\not \in sp^{\Cal{C}} (\phi)$)  and ${\Cal{L}}
\phi : V\to X$ (respectively ${\Cal{C}} \phi : V\to X$) is a
holomorphic extension of ${\Cal{L}} \phi$ (respectively ${\Cal{C}}
\phi$) to an open neighbourhood $V\st \cc$ of $i\,\om_0$, then

(4.3)\qquad $\overline {\Cal{L}}_u \phi (\la) (s)= e^{\la \,s}
(\overline {\Cal{L}} \phi(\la)-\int_{0}^{s} e^{-\la \,t}\phi (t)
\, dt)$,

\qquad\qquad \,\,\, $\overline {\Cal{C}}_u \phi (\la) (s)= e^{\la
\,s} (\overline{\Cal{C}} (\phi)(\la)-\int_{0}^{s} e^{-\la \,t}\phi
(t) \, dt)$

\noindent  are  extensions of  ${\Cal{L}}_u \phi$ and $
{\Cal{C}}_u
 \phi$ respectively; moreover,   $\overline {\Cal{L}}_u \phi(\la) \in BUC(\r_+,X)$
 and
 $\overline {\Cal{C}}_u \phi(\la) \in BUC(\r,X)$
   for each $\la \in
 V$.
\endproclaim
\demo{Proof} (i),(ii).  Set $f_{\la}(t) = \cases{ e^{-\la t},
\text{\,\, if\,\,} t \ge 0}\\ { 0,\text {\quad \,\,\, if \,} t<
 0}\endcases$ and $g_{\la}(t)= -t\, f_{\la}(t)$.  Then  $f_{\la}$, $ g_{\la} \in L^1(\r)$ for all $\la\in \cc_+$.
  The proof of the case $\jj=\r_+$ can be
reduced to the case $\jj=\r$  noting that $\Cal{L}_u \phi
(\la)(s)= \check{f_{\la}}*{\phi}(s)$ for $s\in\r^+$.  So, we
consider the case $\jj=\r$. For Re $\la
 >0$, we have
  $\Cal{L}_u \phi (\la)(s)=\Cal{L} \phi_s (\la)=\Cal{L^+} \phi_s (\la)$ and
 $\Cal{L}_u \phi (\la)(s)= \Check{f_{\la}}*\phi(s)$.
 This implies
$\Cal{L}_u \phi(\la)\in BUC(\r,X)$ (see [4, Proposition 1.3.2
(c)]). Now, let $(\la_n)\st \cc$, $\la_n \to \la$ as $n\to
\infty$, $\la_n \not = \la$ and Re $\la_n \ge $ Re $\la/2$. One
can show that for $t \ge 0$, one has $|\frac{ e^{-\la_n\, t} -
\,e^{-\la\,t}}{\la_n -\la}|\le  t
 e^{-(1/2)\text{Re\,}\la\,t}$,
  and so  by the Lebesgue
 dominating convergence theorem $ || \frac{f_{\la_n} -  f_{\la}}{\la_n -\la}- g_{\la}||_{L^1}
\to 0 $ as  $n\to \infty$. This implies $||\frac{\Cal{L}_u \phi
({\la_n}) - \Cal{L}_u \phi({\la})}{\la_n -\la} - \Check
{g_{\la}}*\phi||_{\infty}\to 0$ as $n\to \infty$. Hence $\Cal{L}_u
\phi$  is differentiable at $\la$  and $\frac {d\Cal{L}_u \phi
(\la)}{d\la}= \check {g}_{\la}*\phi$. So, $\Cal{L}_u \phi$ is
holomorphic on $\cc_+$. If $\phi \in BUC(\r,X)$, then $S(\cdot)$
is strongly continuous, and so  $\Cal{L^+} \phi_s\, (\la)=\int
_0^{\infty}\, e^{-\la \, t} \phi (s +t)\, dt= \int _0^{\infty}\,
e^{-\la \, t} S(t)\phi (s)\, dt= (\int _0^{\infty}\, e^{-\la \, t}
S(t)\phi \, dt)(s)$. The last equality follows  by [32, Corollary
2, p. 134] since evaluation at $s$ is a bounded linear operator.
This proves (i) and (ii) on $\cc_+$ because $\Cal{C}_u \phi
(\la)(s)=\Cal{L^+} \phi_s (\la)$.
  The  case $\text{Re}\la < 0$ being
similar, one gets (ii).

(iii) If $\la \in V\cap \cc_+$, then $\overline {\Cal{L}}_u \phi
(\la) = \overline {\Cal{L}}_u \phi (\la) \in BUC(\r_+,X)$ by part
(i). If $\la=i\,\om \in V\cap i\,\r$, then
  $\overline {\Cal{L}}_u \phi (i\,\om) (s)= e^{i\,\om
\,s} (\overline{\Cal{L}}\phi(i\,\om)-\int_{0}^{s} e^{-i\,\om
\,t}\phi (t) \, dt)$ is  bounded uniformly continuous  by  Theorem
2.4 (iv). This implies $\overline {\Cal{L}}_u \phi (i\,\om)\in
BUC(\r_+,X)$. If $\la \in V\cap \cc_-$, then $\overline
{\Cal{L}}_u \phi (\la) (s)= e^{\la \,s}
 \overline {\Cal{L}} \phi(\la)-\int_{0}^{s} e^{\la \,(s-t)}\phi (t) \, dt)= e^{\la \,s}
 \overline {\Cal{L}} \phi(\la)- h_{\la}*\overline {\phi} (s)$, where $h_{\la} (t)=
 e^{\la\,t}$ if $t\ge 0$ and $h_{\la} (t)=0$ if $t <0$ and $\overline {\phi}$ is as in the proof
 above. Since $\text {Re}\la < 0$, $h_{\la}\in L^1 (\r)$ and
 $e^{\la\,s}$ is uniformly continuous and bounded on $\r_+$. This implies $\overline {\Cal{L}}_u \phi
(\la) \in BUC(\r_+,X)$ for each $\la\in V$.

If $\la \in V\cap \cc\setminus i\r$, then $\overline {\Cal{C}}_u
\phi (\la) = \Cal{C}_u \phi (\la) \in BUC(\r,X)$ by part (ii). If
$\la=i\,\om \in V\cap i\,\r$, then
  $\overline {\Cal{C}}_u \phi (i\,\om) (s)= e^{i\,\om
\,s}(\overline{\Cal{C}}\phi(i\,\om)-\int_{0}^{s} e^{-i\,\om
\,t}\phi (t) \, dt)$ is uniformly continuous and  bounded by [8,
Corollary 4.4 ]. This proves $\overline {\Cal{C}}_u \phi (\la) \in
BUC(\r,X)$ for each $\la\in V$.
 \P
\enddemo

 If $\phi \in L^{\infty}(\jj,X)$ then $\om\in \r$ is said to be
L-$uniformly\,\, regular$ (respectively C-$uniformly\,\, regular$)
for $\phi$  if $i\,\om$ is regular for $\Cal{L}_u \phi: \cc_+\to
BUC(\r_+,X)$ (respectively $\Cal{C}_u \phi:  \cc\setminus i\,\r\to
BUC(\r,X)$). The corresponding $uniform\,\, spectra$ are the sets
 $sp^{ \Cal{L}_u} (\phi)$ \, and \,  $sp^{ \Cal{C}_u}
(\phi)$ \,\, of real numbers which are  not L- and C-uniformly
regular respectively.

Since  $\Cal{L}_u \phi (\la) (0)=\Cal{L} \phi (\la) $, it follows

(4.4)\qquad $ sp^{\Cal{L}} (\phi)  \st sp^{\Cal{L}_u} \phi  \,\,$
and $\,\, sp^{\Cal{C}} (\phi)  \st sp^{\Cal{C}_u} \phi  $.

\proclaim{Proposition 4.2} Let $\phi \in L^{\infty}(\jj,X)$. Then
$sp^{\Cal{L}}(\phi)= sp^{\Cal{L}_u} (\phi)$ for $\jj=\r_+$ and
$sp^{\Cal{C}}(\phi)= sp^{\Cal{C}_u} (\phi)$ for $\jj=\r$.
\endproclaim
\demo{Proof} By (4.4) we need to prove $\om_0 \not \in
sp^{\Cal{L}}(\phi)$ (respectively $\om_0 \not \in
sp^{\Cal{C}}(\phi)$) implies $\om_0 \not \in sp^{\Cal{L}_u}
(\phi)$ (respectively $\om_0 \not \in sp^{\Cal{C}_u}(\phi)$) . By
Lemma  4.1 (iii),   $\overline {\Cal{L}}_u \phi (\la)\in
BUC(\r_+,X)$ (respectively $\overline {\Cal{C}}_u \phi (\la)\in
BUC(\r,X)$) for each $\la \in V$. Moreover, $\overline {\Cal{L}}_u
\phi (\cdot)(s)$ (respectively $\overline {\Cal{C}}_u \phi
(\cdot)(s)$ is holomorphic on $V$ for each $s\in \r^+$ (
respectively $s\in \r$). By (4.3), this implies $x^*\circ
(\overline {\Cal{L}}_u \phi (\cdot)(s))=\overline {\Cal{L}}
(x^*\circ \phi_s) $ (respectively $x^*\circ (\overline {\Cal{C}}_u
\phi (\cdot)(s))=\overline {\Cal{C}} (x^*\circ \phi_s) $) is
holomorphic on $V$ for each $s\in \r^+$ (respectively $s\in \r$)
and each $x^*\in X^*$. This implies $\overline {\Cal{L}}_u \phi$
($\overline {\Cal{C}}_u \phi $) is holomorphic on $V$ by [19,
Definition 3.10.1, Theorem 3.10.1] since the set of functionals
$\{\phi\to x^*\circ \phi(s): x^*\in \ X^*, s\in \jj\}$ is a total
subspace
  [16, p. 418]  of $(BUC(\jj,X))^*$. This proves $\om_0 \not \in sp^{\Cal{L}_u}
(\phi)$  ($\om_0 \not \in sp^{\Cal{C}_u} (\phi)$).
 \P
\enddemo

\proclaim{Corollary 4.3} Let $\phi\in BUC(\jj,X)$. Then

\qquad $sp^{\Cal{L}}(\phi)= sp^{\Cal{L}} (S^{\r_+}(\cdot)\phi)\,$
for $\jj=\r_+$ and $\,sp^{\Cal{C}}(\phi)$$= sp^{\Cal{C}}
(S^{\r}(\cdot)\phi)$ for $\jj=\r$.
\endproclaim
\demo{Proof} By Lemma 4.1, $\Cal{L}_u \phi
=\Cal{L}(S^{\r_+}(\cdot)\phi)$ and $\Cal{C}_u \phi
=\Cal{C}(S^{\r}(\cdot)\phi)$. This implies $sp^{\Cal{L}_u}(\phi)=
sp^{\Cal{L}} (S^{\r_+}(\cdot)\phi)$ and $sp^{\Cal{C}_u}(\phi)=
sp^{\Cal{C}} (S^{\r}(\cdot)\phi)$. By Proposition 4.2, one gets
$sp^{\Cal{L}}(\phi)= sp^{\Cal{L}} (S^{\r_+}(\cdot)\phi)$ and
$sp^{\Cal{C}}(\phi)= sp^{\Cal{C}} (S^{\r}(\cdot)\phi)$.
 \P
\enddemo

\head {\S 5 Conditions for $sp^{\Cal{L}}=sp^{\Cal{C}}$}\endhead

 In the following we indicate a subclass of $L^{\infty}
(\r,X)$ for which the half-line spectrum  $sp^{\Cal{L}}$ coincides
with Carleman spectrum $sp^{\Cal{C}}$. This class includes almost
periodic, almost automorphic, Levitan almost periodic and
recurrent functions (see [1], [5], [7], [13], [23]).  For that let
$\phi \in L^{\infty}(\r,X)$. Set

(5.1)\qquad $L(\phi)=\overline{\text{\, span}}\{\phi_t:
t\in\r\}$\,\,\, and\,\,\, $L^+(\phi)=L(\phi)|\,\r_+$,

(5.2)\qquad $LC(\phi)=\overline{\text{\, span}}\{\phi*f:
 f\in L^1(\r)\}$,\, \, $LC^+(\phi)=LC(\phi)|\,\r^+$,

(5.3)\qquad $m: LC(\phi)\to LC^+(\phi)$,\,\,\, where \,\,\,
 $m(\psi)=\psi|\r_+$.

\noindent Note that $L(\phi)$, $LC(\phi)$ are closed translation
invariant subspaces of $L^{\infty}(\r,X)$, $BUC(\r,X)$
respectively. Moreover, using Bochner integration if $\phi\in
BUC(\r,X)$, then $L(\phi)=LC(\phi)$ (see for example [5, Lemma
1.2.1]).

\proclaim{Theorem 5.1} Let $\phi \in L^{\infty}(\r,X)$. Assume the
restriction mapping $m: LC(\phi)\to LC^+(\phi)$   is an isometric
linear bijection. Then $sp^{\Cal{L}} (\phi)=sp^{\Cal{C}} (\phi)$.
\endproclaim

\demo{Proof} First, we prove the case  when $\phi \in BUC(\r,X)$.
Let $i\om_0$ be a regular point of $\Cal{L} S(\cdot)\phi$ and
$\overline{\Cal{L}} S(\cdot)\phi$ be a holomorphic extension to
$\cc \cup V$, where $V$ is an open disk  with center $i\om_0$.
Since $S(\cdot)\phi$ is Bochner integrable in $BUC(\r,X)$ it
follows that $\Cal{L} S(\cdot)\phi (\la)\in L(\phi)$ for each
$\la\in \cc_+$  and hence $\overline{\Cal{L}} S(\cdot)\phi
(i\om)\in L(\phi)$ for each $i\om \in V\cap i\r$. Using a Taylor
expansion for $\overline{\Cal{L}} S(\cdot)\phi$ about the point
$i\om_0$, one can conclude $\overline{\Cal{L}} S(\cdot)\phi
(\la)\in L(\phi)$ for each $\la\in V$. Since $m$ is an isometric
linear bijection, we conclude $m\circ\overline{\Cal {L}}
S(\cdot)\phi $ is a holomorphic extension to $V$ of $\Cal{L}
S^+(\cdot)\phi|\,\r_+$ and hence $i\om_0$ is a regular point of
$\Cal{L} S^+(\cdot)\phi|\, \r_+$. This implies $sp^{\Cal{L}}
(S^+(\cdot)\phi|\,\r_+) \st sp^{\Cal{L}} (S(\cdot)\phi)$. The
converse can be proved similarly since $m^{-1}: LC^+(\phi)\to
LC(\phi)$ is also a linear   isometric mapping. So, $sp^{\Cal{L}}
(S^+(\cdot)\phi|\,\r_+) = sp^{\Cal{L}} (S(\cdot)\phi)$. Hence by
 Corollary 4.3 and (3.6), we conclude $sp^{\Cal{L}} (\phi)= sp^{\Cal{L}}
(S^+(\cdot)\phi|\,\r_+) = sp^{\Cal{L}} (S(\cdot)\phi)=
sp^{\Cal{C}} (S(\cdot)\phi)=sp^{\Cal{C}} (\phi)$.

 Now, assume  $\phi \in
L^{\infty}(\r,X)$. Then $M_h\phi \in LC(\phi)\st BUC(\r,X)$ for
each $h>0$. It follows $LC(M_h \phi)$ is a closed translation
invariant subspace of $LC(\phi)$. The assumptions imply $m: LC(M_h
\phi)\to LC^+(M_h \phi)$  is an isometric linear bijection. So, by
the above $sp^{\Cal{L}} (M_h\phi)=sp^{\Cal{C}} (M_h\phi)$ and by
Proposition 1.1 (ii), we conclude $sp^{\Cal{L}}
(\phi)=\cup_{h>0}\,\, sp^{\Cal{L}} (M_h\phi)=\cup_{h>0}\,\,
sp^{\Cal{C}} (M_h\phi)=sp^{\Cal{C}} (\phi)$. \P
\enddemo

In the following $AP(\r,X),AA(\r,X), LAP_b(\r,X), RC_b(\r,X)$ will
denote respectively the class of almost periodic, almost
automorphic,  bounded Levitan-almost periodic and continuous
bounded recurrent functions.

  \proclaim{Corollary 5.2} Let
$\phi\in \A\in \{AP(\r,X),\,\,AA(\r,X),  LAP_b(\r,X),\,\,
RC_b(\r,X) \}$. Then
 $sp^{\Cal{L}} (\phi)$$=sp^{\Cal{C}}
(\phi)$.
\endproclaim
\demo{Proof} This follows from Theorem 5.1, since $LC(\phi)\st \A$
and $m: LC(\phi) \to LC^+(\phi)$ is a linear isometric bijection,
by  [5, Theorem 2.1.9].
\enddemo

  \proclaim{Remark 5.3} (i) Let
$\phi (t)= e^{it^2}$. Then $L(\phi)= \phi\cdot AP(\r,\cc)$ and $m:
 L(\phi)\to  L^+(\phi)$ is a linear isometric bijection. But by Example
 2.2, $sp^{\Cal{L}}
(\phi)=\emptyset$ and $ sp^{\Cal{C}} (\phi)=\r$. Also, note that
$LC(\phi)\st C_0(\r,X)$ and so $L(\phi)\cap LC(\phi)=\{0\}$.

(ii) If $\phi\in D(\r,X)$, the Banach space of  distal functions
[13, p. 177] or $\phi\in \m\A \cap L^{\infty} (\r,X)$ with $\A$ as
in Corollary 4.6, then also, $sp^{\Cal{L}} (\phi)=sp^{\Cal{C}}
(\phi)$.
\endproclaim

\smallskip

\smallskip

\indent School of Math. Sci., P.O. Box No. 28M, Monash University,
 Vic. 3800.

\indent E-mails "bolis.basit\@sci.monash.edu.au",\qquad
"alan.pryde\@sci.monash.edu.au".

\enddocument

************************************************

Let $ T(t)\in L(X)$, $t\in\r$ be a bounded $C_0$-group and $x\in
X$,

$Y=\overline{\text{span}\{y: y= T(t)x, \, t\i \r\}} $ and $T_x
(\cdot)=T (\cdot)| Y$.

\noindent  Denote by $A_x$ the generator of $T_x (\cdot)$.

\proclaim{Corollary 2.5}

(2.9)$ \qquad   sp^{\Cal{L}} (T(\cdot)x)=  sp (T(\cdot)x)$.
\endproclaim
\demo{Proof} This follows from (2.5) because $sp (T_x(\cdot))=sp
(T(\cdot)x)=sp (T(\cdot)y)$ for each  non zero element $y\in Y$.
\P
 \enddemo

For $\phi \in BUC(\r,X)$ let $\Lambda(\phi)= \overline{
span\{\phi_t :\,\, t\in\r\}}$.  Then $\Lambda(\phi)$ is a
translation invariant subspace of $BUC(\r,X)$ and $\phi*f\in
\Lambda(\phi)$ for each $f\in L^1(\r)$. Set $S_{\phi}(t)\psi=
\psi_t$, $t\in \r$, $\psi\in \Lambda(\phi)$. Then
 $S_{\phi}(t)\in L(\Lambda(\phi))$, $t\in \r$
is a $C_0$-group.

For $\zeta\in L^{\infty}(\r,X)$ define

 $I(\zeta)=\{f\in L^1(\r):\, \zeta *f =0\}$\,\,\,\, and \,\,\,\, $s$-$I (S_{\phi}(\cdot)):=\cap_{\psi\in \Lambda(\phi)} I(S_{\phi}(\cdot)\psi)$.

\noindent In the special case when $S_{\phi}(\cdot)\in
L^{\infty}(\r, L(\Lambda(\phi)))$, one gets $s$-$I (S(\cdot))=I
(S(\cdot))$.

\noindent Simple calculations show that for $\phi \in BUC(\r,X)$,
one has

$I(S_{\phi}(\cdot)\psi)=I(\psi)$ and $I(\phi)\st I(\psi)$ for each
$\psi\in \Lambda (\phi)$.

\proclaim{Corollary 1.6} Let $\phi \in AP(\r,X)$. Then
$sp^{\Cal{L}} (\phi)=sp^{\Cal{C}} (\phi)=-i \sigma (D_{\phi})$.
\endproclaim
\demo{Proof} \P
\enddemo

\proclaim{Proposition 1.7} Let $\phi^+ \in BUC(\r_+,X)$. Then
$sp_u (\phi^+)=sp^{\Cal{L}} (\phi^+)$.
\endproclaim
\demo{Proof} Denote by $\phi\in BUC(\r,X)$ an extension of
 $\phi^+$ and consider $\Lambda(\phi)$. As in [Basit, SF, Lemma
 6.3(i)], $\sigma (D_{\phi})= P\sigma (D_{\phi})$, the point
 spectrum of $D_{\phi}$. Assuming $\om \in sp_u (\phi^+)$ but $\om \not\in sp^{\Cal{L}}
 (\phi^+)$, there exists  an open  neighbourhood $V\st \cc$ of
 $i\om$ and a unique holomorphic $g: V\to X$ such that $g (\la)= \Cal {C}
 {\phi}(\la)$, Re $\la >0$, $\la \in V$. Define  $g (\la, s): =e^{\la \, s}(g (\la)-\int_0^s e^{-\la \, s}\phi (t)\,
 dt)$. One can verify  $g(\la,\cdot)\in C_u(\r_+, X)$, $\la \in V$ and
 $g(\cdot,s)$ is holomorphic on $V$ for each $s\in \r_+$. Also,
 one gets
 $(\la I-D_{\phi}) g(\la,s)= \phi (s)$, $s\in\r_+$, $\la \in V$ and means that $(\la I-D_{\phi})^{-1}\phi=g(\la,\cdot) $.
  This implies that $i\om $ is not in $P\sigma (D_{\phi})$ by [Hille-ph, p
  54]. But by Theorem 1.4, $\la \in sp_u \phi^+ \st sp \phi = i\sigma (D_{\phi})= i P\sigma (D_{\phi})$. This is a contradiction which proves our claim. \P
\enddemo

\proclaim{Proposition 1.9} Let $\phi \in L^{\infty}(\r_+,X)$. Then
$sp^+_u (\phi)=sp^+ (\phi)$.
\endproclaim

\demo{Proof} Similarly as in the Proof of Proposition 1.4.\P
\enddemo

\proclaim{Corollary 1.10 (Ingham)} Let $\phi \in L^{\infty}(\r,X)$
and $i\om \in i\,\r \setminus sp^+ (\phi)$. Then sup$_{s \ge 0}
||\int_0^ s e^{-i\om t}\phi(t)\, dt|| <\infty$.
\endproclaim
\demo{Proof} Let $\overline{\psi}: V\to BUC(\r_+,X)$ be the
holomorphic extension of  $\psi $ to a neighbourhood of $i\om \in
i\,\r \setminus sp^+ (\phi)$. Then

 $\overline{\psi}(\la)(s)= e^{\la \,t}(
\overline{\psi}(\la)(0)- \int _0^s e^{-\la \,t}\phi (t)\, dt)$.

\noindent In particular, $\overline{\psi}(i\om)(s)= e^{i\om \,t}(
\overline{\psi}(i\om)(0)- \int _0^s e^{-i\om \,t}\phi (t)\, dt)$
and hence

$||\int _0^s e^{-i\om \,t}\phi (t)\, dt||\le 2 ||\overline{\psi}
(i\om)||_{\infty}$. \P
\enddemo

\proclaim{Proposition 1.11} Assume $\phi\in BUC(\r,X)$. Then $\om
\in sp_{\A} \phi$ if and only if for each $f\in L^1 (\r)$, $\hat
{f} (\om)\not = 0$ implies $\phi*f |\jj \not \in \A$.
\endproclaim

 \demo{Proof}  Assume for each $f\in L^1
(\r)$, $\hat {f} (\om)\not = 0$ implies $\phi*f |\jj \not \in \A$.
If $\om \not \in sp_{\A} \phi$, then by [bb, Proposition 2.5 (e)],
  there is $g \in  L^1 (\r)$ such that $\hat {g} (\om)\not = 0$  and $\phi*g |\jj \in \A$.
This is a contradiction which proves $\om \in sp_{\A} \phi$.

Conversely,  assume $\om \in sp_{\A} \phi$ and there is $f\in L^1
(\r)$ such that $\hat {f} (\om)\not = 0$ and $\phi*f |\jj  \in
\A$. By [bb, Proposition 2.5 (e)], one gets $\om \not\in sp_{\A}
\phi$, a contradiction which proves for each $f\in L^1 (\r)$,
$\hat {f} (\om)\not = 0$ implies $\phi*f |\jj \not \in \A$. \P
\enddemo

In the following we use the following space

(1.3)\qquad $\A ^{p,\infty}= L^{p}(\r,X)\cap BUC(\r,X)=
 L^{p}(\r,X)\cap C_0(\r,X)$, $1\le p< \infty$

\noindent which is Banach space endowed with the norm
$||\phi||_{p,\infty}= ||\phi||_{p}+ ||\phi||_{\infty} $.
Obviously, $\A^{p,\infty}$ is translation invariant and closed
under multiplication by characters.

As is customary   we define

(1.4)\qquad $f*\phi  (t)=\int _{-\infty}^{\infty} f(s)\phi (t-s)\,
ds$ for $1\le p \le \infty$,

\qquad \qquad \qquad  $f \in L^1 (\r,\cc)$, $\phi \in L^p (\r,X)$
or $f \in L^1 (\r,X)$, $\phi \in L^p (\r,\cc)$

\noindent the integral exists as a Bochner integral for a.e. $t\in
\r$. Moreover,
 $ \phi*f= f*\phi \in L^p
(\r,X)$ and $||\phi*f||_p \le ||\phi||_p \,||f||_1$.
 (See [abhn,
Propositions 1.3.1, 1.3.2]).

\proclaim{Lemma 1.11} Let $\phi \in L^p (\r,X)$, $f \in L^1
(\r,\cc)$. Then

$ \phi *f =  \int _{-\infty}^{\infty} f(s)\phi _{-s}\, ds$, as a
Bochner integral in $L^p (\r,X)$, $1\le p < \infty$.
\endproclaim
 \demo{Proof} First we prove the special case  where $\phi \in \A^{p,\infty}$ of (1.3).
 Indeed by the dominated convergence theorem and H\"{o}lder's inequality, $ \phi *f (\cdot) = \int _{-\infty}^{\infty} f(s)\phi
(\cdot -s)\, ds\in   L^p (\r,X)\cap C_0 (\r,X)$. Set  $ \psi (s)=
\phi_{-s}$. By [Rudin, 1.1.5, for the case $X=\cc$ but also valid
for general $X$],  $\psi: \r\to \A^{p,\infty}$ is a bounded
uniformly continuous function. It follows $f\psi:\r\to
\A^{p,\infty}$ is Bochner integrable in each of  $\A^{p,\infty}$,
$L^p (\r,X)$ and $C_0 (\r,X)$ . Moreover, the three integrals are
equal. This implies $\int _{-\infty}^{\infty} f(s)\phi _{-s}\,
ds\in  L^p (\r,X)\cap C_0 (\r,X)$ . Consider the evaluation
operator $\epsilon_t \zeta =\zeta (t)$ where $\zeta \in C_0
(\r,X)$. Then $\epsilon_t\in L(C_0 (\r,X), X) $ and one concludes
$\epsilon_t \int _{-\infty}^{\infty} f(s)\phi _{-s}\, ds =\int
_{-\infty}^{\infty} f(s)\phi (t-s)\, ds =\phi*f (t)$.  This proves
the special case. If $\phi \in L^p (\r,X)$, $1\le p < \infty$,
there is $(\phi_n) \st C_0 (\r,X)\cap L^p (\r,X)$, $||\phi_n -
\phi|| _{L^p }\to 0$ as $n\to \infty $.  Set $\psi_n (s)
=(\phi_n)_{-s}$. By  the dominated convergence theorem  one gets
$lim_{n\to \infty} \int _{-\infty}^{\infty} f(s)\psi_n (s)\, ds=
\int _{-\infty}^{\infty} f(s)\phi_{-s}\, ds $. Also, using the
inequality $||\phi*f||_{L^p} \le ||\phi||_{L^p}||f||_{L^1}$, one
gets $||\phi_n*f -\phi*f||_{L^p}\to 0$ as $n\to \infty$. This and
the special case proves our claim. \P
\enddemo

{\bf{\text{Remark}}}. The proof is not valid for the case
$p=\infty$ because $s\to \phi_s$ is not necessary  measurable.

\

Let $A$ be a  linear operator with domain $D(A)\st X$. Define the
composition operator $A_p\phi (t)=A\circ\phi(t)= A \phi(t) $,
$t\in \r$ with domain

$D(A_p)=\{\phi \in L^p (\r, X): \phi(t)\in D(A)$ for almost all
$t\in\r$, $A\circ\phi \in  L^p (\r, X)\}$.

\noindent We note that if $\phi \in D(A_p)$, then we can and often
will assume that $\phi (t)\in D(A)$ for all $t\in \r$.

\proclaim{Proposition 1.12} If $A$ is a closed operator in $X$
with domain $D(A)\st X$, then  $A_p$ is a closed operator,  $1\le
p\le \infty$.
\endproclaim

\demo{Proof} Let $(\phi_n)\st D(A_p)$ and $||\phi_n- \phi||_{p}\to
0$, $||A_p\phi_n- \psi||_{p}\to 0$ as $n\to \infty$. This implies
that there is a subsequence $(\phi_{n_k})\st (\phi_n)$ which
converges to $\phi$ a.e on $\r$ (see the proof of Theorem 10.42
[Rudin, Principles of Mathematical Analysis, p. 216] for the case
$X=\cc$, $p=2$ valid also for any $X$, $1\le p <\infty$. Since $A$
is closed, one concludes $\phi (t)\in D(A)$ and $A\phi (t)=\psi
(t)$ for almost all $t\in \r$. It follows $\phi \in D(A_p)$ and
$\psi=A_p\phi$. \P

\enddemo

\proclaim{Proposition 1.13} Let $A$ be a closed linear operator
with domain $D(A)\st X$. Assume $\phi\in L^{p}(\r,X)$ for some
$1\le p \le \infty$, $\phi \in D(A_p)$. Then $sp(A_p\phi) \st
sp(\phi)$.
\endproclaim
 \demo{Proof} Denote by $I(\phi)=\{f\in L^1(\r,\cc): \phi*f =0\}$.
By [abhn, Proposition 1.1.7] and (1.4), one has $(A_p\phi)*f (t)
  =\int_{-\infty}^{\infty} f(s) A\phi (t-s)\, ds =A \int_{-\infty}^{\infty} f(s) \phi(t-s)\,
ds=A_p (\phi*f)(t)$ a.e. $t\in \r$. This implies $I(\phi)\st
I(A\circ\phi)$ and hence $sp(A_p\phi) \st
 sp(\phi)$.
  \P
\enddemo

\proclaim{Proposition 1.14} Let $\phi\in L^{p}(\r,X)$, $1\le p\le
\infty$ and $0\not \in sp (\phi)$. Then  $P\phi  \in \E(\r,X)\cap
BUC (\r,X)$ and $(P\phi - MP\phi)\in  L^{p}(\r,X)\cap BUC(\r,X)$.
If $\phi (\r)$ is relatively compact (respectively  weakly
relatively compact), then $P\phi(\r)$ is relatively compact
(respectively  weakly relatively compact).
\endproclaim
 \demo{Proof} By H\"{o}lder's inequality, one can verify that $P\phi$ is uniformly continuous.  By [bp, 95, Theorem 3.6],  $P\phi$ is weakly
 bounded and therefore bounded by  [Dsc, 20 Theorem, p. 66]. So $P\phi\in
 BUC(\r,X)$.
 Choose  $\delta >0$ such that $[-\delta,\delta]\cap sp (\phi)=\emptyset$. There is $f\in
 L^1 (\r,\cc)$ such that $\hat {f} (0)=1$ and supp $\hat {f} \st
 [-\delta,\delta]$.  But $(P\phi*f)' =
 \phi*f=0$,
 so there is $x\in X$ such
 that $(P\phi-x)* f=0$. This implies $0\not \in sp (P\phi -x)$.
 Using the same argument as above, it follows $P(P\phi-x)\in
 BUC(\r,X)$ implying
 $(P\phi -x) \in \E(\r,X)$ with mean
 $ M P\phi- x=0$. For each $s\in \r$ we have $\Delta_s P\phi
=\phi*\chi_{[-s,0]}\in L^p (\r,X)$ and so $\Delta_s P\phi\in
\A^{p,\infty}$. Moreover, if $\phi (\r)$ is relatively compact
(respectively  weakly relatively compact) then $\Delta_s P\phi
(\r)= \phi * \chi_{[-s,0]}\st \overline{co}\phi (\r)$, so
$\Delta_s P\phi (\r)$ is relatively compact  (respectively  weakly
relatively compact) for each $s\in \r$. Therefore,  by [BB, 96,
Theorem 2.4 (ii)], $P\phi -MP\phi\in \A^{p,\infty}$, $1\le p
<\infty $ and
 if $\phi (\r)$ is relatively compact  (respectively  weakly relatively compact) then $P\phi (\r)$ is
also relatively compact (respectively  weakly relatively compact)
by [bb95, Theorem 3.1.2].
 \P
\enddemo

Let $\Lambda $ be a closed subset of $\r$.  Set

 $\Lambda^{p}=\{ \phi
\in L^{p}(\r,X): \,\, sp (\phi)\st \Lambda\}$ , \qquad $1\le p\le
\infty$,

$\Lambda_{rc}=\{\phi \in \Lambda^{\infty}: \phi (\r)$ is
relatively compact $\}$.

It is easy to check that $\Lambda^{p}$ and $\Lambda_{rc}$ are
closed translation invariant subspaces of $L^{p}(\r,X)$
respectively $L^{\infty}(\r,X)$. Denote by $B_{\Lambda_p}$,
$B_{\Lambda_{rc}}$ the differential operators acting respectively
on $\Lambda^{p}$, $\Lambda_{rc}$. Then one can verify that
$B_{\Lambda_p}$, $B_{\Lambda_{rc}}$ are closed operators in
$L^{p}(\r,X) $ with domains

\qquad $D(B_{\Lambda_p})=\{\phi\in \Lambda_{p}: \phi'\in
\Lambda_{p}\}$, \qquad $1\le p \le \infty$,

\qquad $D(B_{\Lambda_{rc}})=\{\phi\in \Lambda_{rc}: \phi'\in
\Lambda_{rc}\}$.

\proclaim{Lemma 1.15} Let $\phi\in \Lambda_{rc}$ respectively
$\phi\in \Lambda^{p}$ and $\beta  \in \r\setminus \Lambda$. Then

(1.5)\qquad $\psi'-i\beta \psi = \phi$

\noindent has a unique solution $\psi =
\gamma_{\beta}[P(\gamma_{-\beta}\phi)- M P(\gamma_{-\beta}\phi)]$
where $ \psi\in \Lambda_{rc}$ respectively $\psi \in \Lambda^{p}$.
\endproclaim
 \demo{Proof} Directly, one can check that $\psi$ satisfies (1.5) and by Proposition
 1.12, one can verify $\psi \in  \Lambda_{rc}$ respectively $\psi
\in \Lambda^{p}$. Uniqueness follows because if $\psi'-i\beta \psi
= 0$  then  $\psi = e^{i\beta t}x$ which gives then $sp
(\psi)=\{\beta\}$ or $x=0$.
\enddemo

\proclaim{Proposition 1.16}
$\sigma(B_{\Lambda_{rc}})=\sigma(B_{\Lambda^{p}})=i\Lambda$, $1\le
p\le \infty$.
\endproclaim
 \demo{Proof} Let $x\in X\setminus 0$, $\om \in \Lambda$.  If $p=\infty$, then
 $x\gamma_{\om}$, $(x\gamma_{\om})'= i\om x\gamma_{\om}\in \Lambda_{rc} \st
 \Lambda^{\infty}$. This proves $
 i\Lambda \st\sigma(B_{\Lambda^{\infty}})$, $
 i\Lambda \st \sigma(B_{\Lambda_{rc}})$.  If $1\le
p < \infty$, one can show that for each $n\in \N$ there is
$f_{n,p} \in C^1 (\r)$  with supp $f_{n,p}\st [-n,n]$,
$||f_{n,p}||_p =1$ and $||f'_{n,p}||_p < 1/n$. If  $\psi_n =
\gamma_{\om} f_{n,p}\, x$, $x\in X\setminus \{0\}$, then
$||(B_{\Lambda_p}-i\om)\psi_n||_{L^p}\to 0 $ as $n\to \infty$,
$1\le p  < \infty$. This implies that $i\Lambda \st \sigma
_{ap}(B_{\Lambda_{p}}) \st \sigma (B_{\Lambda_{p}}) $. Here
$\sigma _{ap}(B_{\Lambda_{p}})$ denotes the approximate point
spectrum (see [abhn, p. 463]).

 Conversely, let $\beta  \in
\r\setminus \Lambda$ and  $\phi\in \Lambda_{rc}$ respectively
$\phi\in \Lambda_{p}$.   By Lemma 1.15, equation (1.5) has a
unique solution. So, by the closed graph theorem it follows
$i\r\setminus i\Lambda\st \rho (B_{\Lambda_{rc}})$ respectively
$i\r\setminus i\Lambda\st \rho (B_{\Lambda^{p}})$. This  completes
the proof of the statement. \P
\enddemo

Now let $A$ be the generator of an analytic $C_0$-semigroup
$T(t)$, $t\ge 0,$ of linear bounded operators on $X$. We consider
the composition operators  $A_{\Lambda_{rc}}$, $A_{\Lambda_p}$:

(1.6)\qquad $D(A_{\Lambda_{rc}})=\{\phi \in \Lambda_{rc}:
\phi(t)\in D(A), t\in \r,\,\,\, A\circ\phi\in \Lambda_{rc}\}$,

(1.7)\qquad $A_{\Lambda_{rc}}\phi (t) =A\circ\phi (t)= A\phi (t)$
for all $\phi\in D(A_{\Lambda_{rc}})$.

(1.8)\qquad $D(A_{\Lambda_p})=\{\phi \in \Lambda_p: \phi(t)\in
D(A), t\in \r,\,\,\, A\circ\phi\in \Lambda_p\}$,

(1.9)\qquad $A_{\Lambda_p}\phi (t) =A\circ\phi (t)= A\phi (t)$ for
all $\phi\in D(A_{\Lambda_p})$.

 \proclaim{Proposition 1.17} The operator $A_{\Lambda_{rc}}$ defined
 by (1.6), (1.7) is the generator of an analytic $C_0$-semigroup on
 $\Lambda_{rc}$.
\endproclaim
 \demo{Proof} Consider the semigroup $T_{\Lambda_{rc}}(t)$ of
 composition by $T(t)$ on $\Lambda_{rc}$. We show that the map
$t\mapsto T_{\Lambda_{rc}}(t)\phi$ is continuous for each $\phi
\in \Lambda_{rc}$. Indeed, the mapping  $t\mapsto T(t)x$ is
continuous  at $0$ uniformly with respect to $x\in \overline
{\phi(\r)}$. So, sup$_{s\in \r} ||T(t)\phi (s)- \phi(s)||\to 0$ as
$t\to 0$. This implies that $T_{\Lambda_{rc}}(\cdot)$ is strongly
continuous. Let $G$ be the generator of $T_{\Lambda_{rc}}(\cdot)$.
We show that $ G=A_{\Lambda_{rc}}$.
 Indeed, let $\phi \in
D(A_{\Lambda_{rc}})$. Then by [25, Theorem 2.4, p. 4], $\frac
{T(t)\phi (s)-\phi (s)}{t}= \frac{1}{t}\int_{0}^{t} T(u) A\phi
(s)\, du $,\,\,\, $t>0$, $s\in \r$. Since
 sup$_{s\in \r}\,||\frac{1}{t}\int_{0}^{t} T(u)
A\phi (s)\, du-A\phi(s)|| \le  \frac{1}{t}\int_{0}^{t}$ sup
$_{s\in \r}|| T(u) A\phi (s) -A\phi(s)||\, du$ $ \to 0 $
 as  $0< t\to 0$,
  one concludes that $\phi \in D(G)$ and
$ G\phi= A_{\Lambda_{rc}}\phi$. Conversely, if $\phi \in D(G)$,
then $\frac{T(t)\phi(s)-\phi(s)}{t}\to (G\phi)(s)$ for a.e. $s\in
\r$. So, $\phi(s)\in D(A)$, $A\phi (s)=(G\phi) (s) $ and $\phi \in
D(A_{\Lambda_{rc}})$ proving $G=A_{\Lambda_{rc}}$. Now, we show
that $\sigma (A_{\Lambda_{rc}})\st \sigma (A)$. Indeed, let $\mu
\in \rho (A)$. We show that that for each $\phi \in \Lambda_{rc}$
the equation $\mu\psi - A_{{\Lambda_{rc}}}\psi= \phi$ has a unique
solution $\psi \in D(A_{{\Lambda_{rc}}})$. Indeed, for a.e. $s\in
\r$ the equation $\mu \,x - A x= \phi (s)$ has a unique solution
$\psi (s)= (\mu-A)^{-1}\phi (s)$. This implies
$(\mu-A_{{\Lambda_{rc}}})^{-1}$ is defined on $\Lambda_{rc}$ and
since it is closed, it is bounded by the closed graph theorem [32,
p. 79]. This proves $\mu \in \rho (A_{{\Lambda_{rc}}})$ and hence
$\sigma (A_{\Lambda_{rc}})\st \sigma (A)$. Since
$||(\mu-A_{\Lambda_{rc}})^{-1}||\le ||(\mu-A)^{-1}||$ for $\mu\in
\rho (A)$ and $A$ is the generator of an analytic $C_0$-semigroup,
we conclude $A_{\Lambda_{rc}}$ is the generator of an analytic
$C_0$-semigroup, by [25, Theorem 5.2, p. 61].
\enddemo

\proclaim{Proposition 1.18} The operator $A_{\Lambda _p}$ defined
 by (1.8), (1.9) is the generator of an analytic $C_0$-semigroup on
 $\Lambda_p$ for each $1\le  p< \infty$.
\endproclaim
\demo {Proof}  Consider the semigroup
 $T_{L^p}(t)$ of
 composition by $T(t)$ on $L^p (\r,X)$.
As in the proof of Proposition 1.15,the map  $t\mapsto
T_{L^p}(t)h$ is continuous
  for each
   $h\in {\Cal H_A}=\{ g: \,$ step function, $g: \r\to D(A) \}$.
Since $D(A)$ is dense in X,  ${\Cal H_A}$ is dense in $L^p
(\r,X)$, $1\le  p< \infty$. It follows that  $T_{L^p}(t)$ is
strongly continuous and so  the semigroup $T_{\Lambda_p}(t)$ of
 composition by $T(t)$ on $\Lambda_p$ is  strongly continuous on $\Lambda_p$.
 One can proceed as in
 the proof of Proposition 1.15. \P
\enddemo

\proclaim{Proposition 1.19} The operators $A_{\Lambda_{rc}}$
$B_{\Lambda_{rc}}$ are commuting and satisfy condition (P).
Moreover,   $\overline{B_{\Lambda_{rc}}- A_{\Lambda_{rc}}}^A
=\overline{B_{\Lambda_{rc}}- A_{\Lambda_{rc}}}$.
\endproclaim
 \demo{Proof} Condition (P) follows by Propositions 1.16 and 1.17.
$\overline{B_{\Lambda_{rc}}- A_{\Lambda_{rc}}}^A
=\overline{B_{\Lambda_{rc}}- A_{\Lambda_{rc}}}$.
\enddemo

\proclaim{Proposition 1.19} Let $A$ be a linear operator with
domain $D(A)$ in $X$. Assume that $(A-\la I)^{-1} \in L(X)$ for
some $\la \in \cc$. Then $A$ is closed.
\endproclaim
\demo {Proof} Let $(x_n)\st D(A)$, $x_n \to x$, $Ax_n\to y$ as
$n\to \infty$. We claim that $x\in D(A)$ and $y=Ax$. Indeed, the
assumptions imply

$(A-\la I)x_n \to y- \la x$,  $x_n= ((A-\la I)^{-1} (A-\la I)x_n
\to (A-\la I)^{-1}(y- \la x)= x$.

\noindent So $(A-\la I)x = y-\la x$ proving our claim. \P
\enddemo

\proclaim{Proposition 1.20} Let $(\phi_n),\, (\phi'_n),
\,\phi,\,\psi  \in L^p (\r,X)$, $||\phi_n- \phi||_{L^p}\to 0$,
$||\phi'_n- \psi||_{L^p}\to 0$ as $n\to \infty$, $1\le p\le
\infty$. Then $\phi'=\psi$.
\endproclaim
\demo {Proof} Let $\zeta\in  L^p (\r,X)$. By H\"{o}lder inequality
with  $1/p+1/q =1$, one has
 $||\int_0^t \zeta (s)\, ds||\le
|t|^{1/q}(\int_0^t ||\zeta (s)||^p\, ds)^{1/p}$ and so $||P
\phi'_n (t)-P\psi(t)||\to 0$ as $n\to \infty$ for each $t\in \r$.
But $P \phi'_n (t)= \phi_n (t)-\phi_n (0)$. Passing to a
subsequence if needed, the assumptions and again by [Rudin, p.
216] as in the proof of  Theorem 10.42 imply $ \phi_n (t)\to \phi
(t)$ as $n\to \infty$ for each $t\in \r$. This implies
$\phi(t)-\phi (0)= \int_0^t \psi (s)\, ds$ and so $\phi'=\psi$, by
[abhn, Proposition 1.2.2 (a), p. 16]. \P
\enddemo

\smallskip

\smallskip

\indent School of Math. Sci., P.O. Box No. 28M, Monash University,
 Vic. 3800.

\indent E-mail "bolis.basit\@sci.monash.edu.au".

\smallskip

\indent E-mail "alan.pryde\@sci.monash.edu.au".

\enddocument

 **********************************

 First, we prove the case $\phi\in BUC(\jj,X)$. If
$\jj=\r_+$ we extend $\phi$ to  $\tilde{\phi} \in BUC(\r,X)$. By
Remark 1.1, we have $ sp_{C_0 (\jj,X)}(\tilde{\phi})= sp_{C_0
(\jj,X)}(\phi)$. Let $g_n(t)= (1/n) \psi (t/n)$. Choose $n_0\in
\N$ such that supp $\widehat {g_n}_0 \cap
sp_{C_0(\jj,X)}(\tilde{\phi})=\emptyset$. It follows
$\tilde{\phi}* g_{n_0}|\,\jj \in C_0 (\jj,X)$ and so
$\tilde{\phi}* g_{n_0}|\,\jj \in \E_{u,0} (\jj,X)$. One has
$\tilde{\phi} (t)-\tilde{\phi}* g_{n_0}(t)=
\int_{-\infty}^{\infty} [\tilde{\phi} (t)-\tilde{\phi} (t-n_0
s)]\psi(s)\,ds$ and $[\tilde{\phi} (\cdot)-\tilde{\phi} (\cdot-n_0
s)]\in \E_{u,0} (\r,X)$. The function  $s\to [\tilde{\phi}
(\cdot)-\tilde{\phi} (\cdot-ns)]\psi(s)$ is in $BUC(\r,
\E_{u,0}(\r,X))$ and is therefore  Lebesgue-Bochner integrable. It
follows $[\tilde{\phi} -\tilde{\phi}* g_{n_0}]\in \E_{0} (\r,X)$.
By the above we conclude $\phi = [\tilde{\phi} -\tilde{\phi}*
g_{n_0}+\tilde{\phi}* g_{n_0}]\,|\jj \in  \E_{0}(\jj,X) $ and
$m([\tilde{\phi} -\tilde{\phi}* g_{n_0}]\,|\jj)=m({\phi})
-m(\phi)\, \widehat{g_n}_0(0) =0$. It follows $\phi\in
\E_{0}(\jj,X)$. This finishes the proof of this case.

Secondly, let $\phi\in L^{\infty}(\jj,X)$. By (1.14), we conclude
$0\not\in sp_{C_0(\jj,X)}(M_h\phi)$ for each $h
>0$. Since $M_h \phi\in
BUC(\jj,X)$, one gets by the above $M_h \phi\in  \E_{u,0} (\jj,X)$
for each $h>0$. One has $ (\phi-M_h \phi)\in \E_0 (\jj,X)$
 by (1.12). It follows $\phi \in \E_0 (\jj,X)$. \P

    \proclaim{Definition 1} For S-bi-invariant subspace $ \A \st
    BUC(\r_+,X)$ denote the quotient space by $BUC(\r_+,X)/\A$,
    the quotient map $\Pi : BUC(\r_+,X) \to BUC(\r_+,X) /\A$ and
    the quotient semigroup defined by $S/\A (t)\Pi \phi= \Pi
    S(t)\phi$ by $S/\A $.
       \endproclaim

 \proclaim{Definition 2} Let $ \A \st
    BUC(\r_+,X)$ be a S-bi-invariant subspace and

    \noindent $\phi \in L^{\infty}
    (\r,X)$. We define $sp^{\r_+}_{\A} (\phi)=\{\eta\in \r:$ for every $\e >0$ there exists

 \qquad\qquad  \qquad $f\in L^1(\r)$ such that supp $\hat{f} \st ]\eta-\e,\eta +\e[$ and $ f\star \phi| \r_+ \not \in\A  \}$
    \endproclaim

here $f\star \phi (t): =\int_0^{\infty} f(t-s)\phi (s)\, ds  $ and
$\hat {f} (\la)= \int_{-\infty}^{\infty} e^{-i\la s} f(s)\, ds$.

By a direct verification one can show

\proclaim{Remark 3} (i) Let  $ \A \st
    BUC(\r_+,X)$ be a S-bi-invariant subspace

(i) If  $\phi_1, \phi_2 \in L^{\infty}
    (\r,X)$, $\phi_1 (t)- \phi_2 (t)=0$, $t \ge 0$, then

$sp^{\r_+}_{\A} (\phi_1)= sp^{\r_+}_{\A} (\phi_2)$.

(ii)  If $\phi \in L^{\infty}
    (\r,X)$. Then
    $ f\star \phi| \r_+  \in\A $ if and only if  $ (f* \phi)| \r_+  \in\A
    $.

(iii) If $\phi \in  BUC(\r_+,X)$, then $\Pi \phi =0$ if and only
if $\phi \in \A$.

(iv) By Lemma 2 of [Chill-Fasangova], the $C_0$-semigroup $S/\A
(\cdot)$ extends to a $C_0$-group of isometries  $S/\A(t)$, $t\in
\r$. We denote its generator by $D/\A$.  For $s<0$ define $S(s):
BUC(\r_+,X)\to BUC(\r_+,X)$  by  $S(s)v(t)= v(0)$ if $s+t\le 0$
and $S(s)v(t)= v(t+s)$.  Then one can verify

$S/\A (t)\Pi \phi= \Pi
    S(t)\phi$, $t\in \r$, $\phi \in BUC(\r_+,X)$.
(See [Chill-Fasangova], proof of  Theorem 5).
\endproclaim

 \proclaim{Definition 4} Let $T(\cdot) $ be a bounded $C_0$-group
 on a Banach space $Y$. The local Averson spectrum of some $y\in
 Y$with respect to $T(\cdot) $ is defined by

 Sp $^A (T(\cdot) y):=\{\eta\in \r:$ for every $\e >0$ there exists $f\in L^1(\r)$
    such that

   \qquad\qquad \qquad\qquad \qquad  supp $\hat{f} \st ]-\eta-\e,-\eta +\e[$ and $ \int_{-\infty}^{\infty}  f(s)\,T(s) y ds\not =0 \}$
    \endproclaim

\proclaim{Proposition 5} Sp $^A (T(\cdot) y)= sp_{0, L^1}(T(\cdot)
y)$.
    \endproclaim

\demo{Proof} Since $T(\cdot) $ is a group, one gets  $
\int_{-\infty}^{\infty}  f(s)\,T(s) y ds =0$ if and only if

$ T(t)\int_{-\infty}^{\infty}  f(s) \,T(s) y ds= T(\cdot)y *
\check {f} (t) =0$, $t\in \r$. Here $\check {f} (t)= f(-t)$,
$t\in\r$. Using the definition of $sp_{0, L^1}(T(\cdot) y)$ one
gets the claim.
\enddemo

\proclaim{Theorem 6} (Chill-Fasangova) Let $ \A \st
    BUC(\r_+,X)$ be a S-bi-invariant subspace and $\phi \in BUC
    (\r_+,X)$. Then \,\,\,
    Sp $^{\r_+}_{\A} (\phi)= $ Sp $^A (S/\A (\cdot)\Pi \phi)$.

\endproclaim
\demo{Proof} Assume $f\in L^1 (\r)$. One has   $f*(S/\A (t)\Pi
\phi)= \Pi  (f* S(t) {\phi}$, $t\in\r$, by Remark 3 (iv) . This
implies $f*(S/\A (t)\Pi \phi)=0$, $t\in\r$ if and only if $ \Pi
(f* S(t) \phi=0$, $t\in\r$. Since $S/\A (\cdot)$ is a group, one
gets $f*(S/\A (\cdot)\Pi \phi)=0$ if and only if  $\Pi
S(0)\int_{-\infty}^{\infty}f (s) S(-s)\phi\, ds= \Pi (f*\overline
{\phi}|\r_+ )=0$. Here  $\overline {\phi}|\r_+ =\phi$ and
$\overline {\phi}(t) =\phi(0)$, $t \le 0$. By Remark 3 (iii), one
gets $ \Pi (f* \overline {\phi}|\r_+ )=0$ if and only if $f*
\overline {\phi}|\r_+ \in \A$. By remark 3 (ii), $f* \overline
{\phi}|\r_+ \in \A$ if and only if $f\star \phi|\r_+ \in \A$. This
proves the claim.
\enddemo

\proclaim{Proposition 1.4} Let $\phi \in L^{\infty}(\r_+,X)$, let
$i\,\om \in i\r\setminus sp^{\Cal{L}} (\phi)$ and let $\overline
{\Cal {L}(\phi)}: V\to X$ be a holomorphic extension of $\CAL
{L}(\phi)$ to a connected open neighbourhood  $V $ of $\cc_+ \cup
i\r\setminus sp^{\Cal{L}} (\phi)$. For each $s\in \r_+$, $\CAL
{L}(\phi_s)$ has a holomorphic extension $\overline {\CAL
{L}(\phi_s)}: V\to X$. Moreover, for each $\la \in V$, there exist
a neighbourhood $U$ of $\la$ in $V$ and a constant $c$ such that
$||\overline {\CAL {L}(\phi_s)} (z)||\le c$, $s\in \r_+$, $z\in
U$.
\endproclaim

 \demo{Proof} See [Batty, Proposition 2.1]. \P
\enddemo

Let $ \phi\in
    BC(\r,X)$ and $ {L}(\phi) := \overline{span\{\phi_t, \phi*f: t\in\r, f\in L^1(\r)\}}$. Then  $L(\phi)$ is
    a translation
      invariant subspace of $BC
    (\r,X)$. Denote by $B_{\phi}$ the generator of the semigroup
    $S(t)\psi$, $\psi \in L(\phi)$. If $\phi \not \in BUC(\r,X)$, then  domain $B_{\phi}=
    D(B_{\phi})$ is not dense in $L(\phi)$. Indeed, assuming the
    contrary,
    $S(\cdot)$ is a $C_0$-group, a contradiction by [Pazy, p.].

\proclaim{Theorem 7}  Let $ \phi\in
    BC(\r,X)$ and $L(\phi) := \overline{span\{\phi_t, \phi*f t\in\r, f\in L^1(\r)\}}$. Then  $L(\phi)$ is
    a translation
      invariant subspace of $BC
    (\r,X)$.
\endproclaim

\enddocument

 R/ASST 3.    BASIT, Bolis        C1  2004000268 Ergodicity
and stability of orbits of unbounded semigroup representations
Journal of Australian Mathematical Society 1446-7887

 Basit, B.,
Pryde, A. J  230000 (100) (230107 (40 $\% $)230108 (60 $\% $)).

R/ASST 3.
 BASIT, Bolis        C1  2004000269  Generalized
Esclangon-Landau results and applications to linear
difference-differential systems in Banach spaces    Journal of
Difference Equations and Applications    1023-6198   Basit, B.,
Gunzler, H.  230000 (100) (230107 (50 $\% $ )230108 (60 $\% $)).

 R/ASST 3.    BASIT, Bolis
C1  2005000825  Difference property for perturbations of vector
valued Levitan almost periodic functions and their analogs Russian
Journal of Mathematical Physics 1061-9208   Basit, B., Gunzler, H.
230107 (100) (230107 (50 $\% $ )230108 (50 $\% $ )).

R/ASST 3.    BASIT, Bolis E1  2003000174  Almost periodic behavior
of unbounded solutions of differential equations Proceedings of
the Centre for Mathematics and its applications Australian
National University 07 3155 2058  Basit, B., Pryde, A. J 230000
(100) (230107 (40 $\% $ )230108 (60 $\% $ )).

 R/ASST 3           . BASIT, Bolis        E1
2003002645  A Layered Approach to Extracting Programs from Proofs
with an Application in Graph Theory   Proceedings of the 7th $\&$
8th Asian Logic Conferences 981-238-261-5   Jeavons, J. S, Basit,
B., Poernomo, I. H, Crossley, J. N    230000 (50),280000
(50),280302 (0).

 R/ASST 3 .    BASIT, Bolis        Others  2000040201  A
Layered Approach to Extracting Programs from Proofs with an
Application in Graph Theory           Jeavons, J. S, Poernomo, I.,
Basit, B., Crossley, J. N  230000 (50),280000 (50),280302 (0).

 R/ASST 3           .
BASIT, Bolis Others  2000040203  Fred: An Implementation of a
Layered Approach to Extracting Programs from Proofs.  Part I:  An
Application in Graph Theory          Jeavons, J. S, Poernomo, I.,
Crossley, J. N, Basit, B.  230000 (50),280000 (50),280302 (0).

 R/ASST 3           .
BASIT, Bolis Others  2000040223  A Case Study for Reliable,
Reusable Software Crossley, J. N, Poernomo, I., Basit, B.,
Jeavons, J. S  230000 (50),280000 (50),280302 (0).

\pagebreak

 *******************

Given an $f \in PS(\r,X)$ and $g = e^{i \omega\, t}$, $\omega$ fixed $\in \r$;
to  $\e_n = 2^{-n}$ one can then find  inductively sequences of
intervals $I_n=[-a_n,a_n]$ and $r_n >0$ ($a_n >0$) with , for $n \in \N$ ,
$ |f_{r_n}  - f | < \e_n$  on $I_n$, $I_n + r_n  \st  I_{n+1}$, $a_n\to \infty$,
$r_n\to \infty$; by taking a subsequence one can further assume , with
complex $c$, $|c|=1$, $| g(r_n) - c | < \e_n$  for  $n \in \N$. Since
$I_n+r_n  \st  I_{n+1}$, the above gives for $n+1$,  $| f_{r_{n+1} + r_n} -
f_{r_n} | < \e_{n+1}$  on $I_n$, and inductively
$|f_{r_n+...+r_{n+k}} - f |  <  \e_n +...+ \e_{n+k}$  on $I_n$, $n,k
\in  \N$, $|g_{r_n+...+r_{n+k}}  -  c^{k+1} g | < \e_n+...+\e_{n+k}$
on $\r$. Now with the Kronecker approximation theorem there exist $k$
($\to \infty$ even) with  $| c^{k+1}  -  1 |  <$ given pos. $\e$, so for
such a $k$ with  $ | g_{{r_n}+...+r_{n+k}} - g | < \e _n +...+\e_{n+k}+\e$ on
$\r$, $n \in \N$. Since $a_n \to \infty$, $r_n \to \infty$  $\e_n+...+\e_{n+k}
\to 0 $ as  $n\to \infty$, the above shows that for $f$, $g$ as above, $(f,g)
\in PS(\r,X  \times \cc$) ; inductively $(f,g_1,...,g;m) \in PS$; since PS
is uniformly closed one gets the

Theorem. $f \in PS(\r,X)$, $ g \in AP(\r,Y)$ implies  $(f,g)
\in PS(\r,X\times Y)$.

Corollary.  For any $X$, $PS_b(\r,X)$  satisfies $(\Gamma)$.

With this, analogues of §§3-5 of RSD follow, even for $PS_b$ instead
of $REC_{urc}$ , resp $rec_{urc}$.

\enddocument

\proclaim {Example 1}  A Levitan-ap $\phi:\Bbb{Z}\to l^1$ with
$R(\phi)=(e_n)$, $(e_n)$ bases of $l^1$.

 Note that $(e_n)$ is
weakly complete and bounded because it has no weakly Cauchy
subsequence. But $(e_n)$ is not weekly relatively compact.

Now Define $\phi (n)=x_k$, \qquad $n\in A_k$, $n\in \Bbb{Z}$ as in
(5.1) of [Basit, Bull sci. math. 1977]. Then $\phi$ is Levitan-ap.
\endproclaim

\proclaim {Example 2} Let  $f:\Bbb {Z}\to \cc $, $f(n)=
\frac{1+e^{in}}{|1+e^{in}|}$. Then $f\in AA(\Bbb {Z})$ [BJM, p.
192, Example 7.16] but $f\not\in AP(\Bbb {Z})$. Denote by $O(f)$
the orbit of $f$. Show that

 $h_{\pm}:\Bbb {Z}\to \cc $, $h_{\pm}(n)= \frac{1-e^{in}}{|1-e^{in}|}$, $n\not=0$, $h_+ (0)=i$, $h_- (0)=-i$

\noindent are elements of $O(f)$ (So, $h_{\pm}\in  RC(\Bbb {Z})$)
but $h_++h_- \not \in RC(\Bbb {Z})$.
\endproclaim
\demo {Proof} Let $a=(m_k)\st \Bbb {Z}$ be such that $ \lim_{k\to
\infty} e^{im_k} = e^{i\theta}$. Then $ \lim_{k\to \infty} f(n
+m_k) = \frac{1+e^{in}e^{i\theta}}{|1+e^{in}e^{i\theta}|} = g(n)$,
$e^{in}e^{i\theta}\not = -1$. If $e^{i\theta} = -1$ we can choose
$a=(m_k)$ such that $ \lim_{k\to \infty} f(n +m_k)= h_+$ or $h_-$.
One can easily show that $\lim_{k\to \infty} g(n -m_k) = f(n)$.
This proves $f$ is almost automorphic. Since $h_+ - h_- (n)=0$,
$n\not =0$ and $h_+ - h_- (0)=2 i$, one gets $h_, h_-$ are not
almost automorphic. This implies $f$ is not almost periodic.

We prove $sp^A (T(\cdot))$ is closed. Let $(\la_n)\st sp^A
(T(\cdot))$ and $\la_n \to \la$ as $n\to \infty$. Let $\e >0$.
There is   $\la_{n(\e)} \in ]\la  -\e,\la +\e[$ and $\delta >0$
such that  $ ]\la_{n(\e)}- \delta, \la_{n(\e)}+\delta [ \st
]\la-\e,\la+\e[$. There is $x_{n(\e)}\in X$ and $f_{n(\e)}\in L^1
(\r)$ and supp $\hat{f}_{n(\e)}\st ]\la_{n(\e)}- \delta, \la
_{n(\e)}+\delta [ $.  This implies supp $\hat{f}_{n(\e)}\st
]\la-\e,\la+\e[$ and hence  $\la \in sp^A (T(\cdot))$.

That $h_+ + h_-$ is not recurrent follows with $h_+ + h_- = 2h_+ +
(h_- -h_+)$. Assuming $h_+ + h_-$ is recurrent, we get $O(h_+ +
h_-) = 2O (h_+)= 2 O(f) $. This implies  $h_+ + h_-=2h$ for some
$h \in O(f)$. Therefore  $2|h(0)|=2$ contradicting $|h_+ (0) + h_-
(0)|=0$.
\enddemo
Let $f\in C_{urc}(\r,X)$. Then $O(f)$ will denote the closure in
the topology of local uniform convergence on $\r$ of the set
$\{f_t: t\in \r\}$. Denote by

(1.1) \qquad $d(g,h):= \sum_{k=1}^{\infty} (||g-h||_k /2^k)$,
where

 \qquad \qquad $||g-h||_k = max_{t\in [-k,k]} ||g(t)-h(t)||$, $g,h \in
C_{urc}(\r,X)$.

One can check that $(O(f), d)$ is a complete metric space.

\proclaim {Definition 1.0} Let $f\in C_{rc}(\r,X)$. Then for each
$(t_{n'}) \in\r$ there is a  subnet  $(t_{i})_{i\in \Lambda}$ such
that $lim_{i\in \Lambda} f(t+t_i)$ exists point wise $t\in\r$.
Denote by $O(f)$ the set of all such limits.

If $O(f) \st C_{rc}(\r,X)$, then  $f\in C_{urc}(\r,X)$.
\endproclaim
\demo{Proof} Let $g\in O(f)\st C_{rc}(\r,X)$ and $g(t)=lim_{i\in
\Lambda} f(t+t_i)$, $t\in \r$. By [DS, Theorem] valid also for
Banach valued function, one gets $(f_{t_i}|[-n,n])_{i\in \Lambda}$
is weakly compact for each $n\in \N$ and  hence weakly
sequentially compact. It follows  using the diagonal process that
there is a subsequence $(t_n)\st (t_{n'})$ such that
$g(t)=lim_{n\to \infty} f(t+t_n)$, $t\in \r$. Now using the same
method of the proof  of [Veech, Lemma 4.1.1 ], one gets $f\in
C_{urc}(\r,X)$
\enddemo
\proclaim {Definition 1.1} A function $f\in C_{urc}(\r,X)$ is
called recurrent if  the set

$\{\tau : d(f, f_{\tau})< \e\}$ is relatively dense in $\r$ for
each $\e >0$.
\endproclaim

\proclaim {Proposition 1.2} If $f\in C_{urc}(\r,X)$,  then the set

$\{\tau : d(f, f_{\tau})< \e$ is relatively dense in $\r$  for
each $\e > 0$,

if and only if the  set

$\{\tau : ||f(t) -f(t+\tau)||< \e , |t|\le N \}$ is relatively
dense for each $\e >0$ and $N > 0$.
\endproclaim

A sequence $\alpha=(a_n)\st \r$ is regular for $f$ if $(f_{a_n})$
converges in $(O(f), d)$. We denote by $f_{\alpha}$ the limit of
$(f_{a_n})$.

A set $F\st C_{urc}(\r,X)$ is called minimal if and only if $F$ is
non-empty, closed and invariant  in $(C_{urc}, d)$ and has no
proper subset having these three properties.

 \proclaim {Definition 1.3} A function $f\in C_{urc}(\r,X)$ is
called minimal if and only if  for each regular for $f$ sequence
$\alpha=(a_n)\st \r$ there exists a sequence $\beta=(b_n)\st \r$
regular for $f_{\alpha}$ such that $(f_{\alpha})_{\beta}=f$.
\endproclaim

\proclaim {Proposition 1.4} A function $f\in C_{urc}(\r,X)$ is
minimal if and only if $O(f)$ is a minimal set in $(C_{urc}, d)$.
\endproclaim

\demo{Proof} Let $O(f)$ be a minimal set in $(C_{ucr}, d)$ and
$\alpha= (a_n)$ be a regular sequence for $f$. Then
$O(f)=O(f_{\alpha})$. This implies $f= (f_{\alpha})$ for some
 regular sequence $\beta$  for $(f_{\alpha})$ and proves that
 $f$ is minimal.

 Conversely, let $f$ be minimal and assume that $\sigma \st O(f)$
 is a proper minimal subset.  Then  $d(f,\sigma) = min_{g\in \sigma} d(f,g) > 0$.
Assume $g= f_{\alpha}\in \sigma$. Then $f\not = g_{\beta}$ for
each regular sequence for $g$. This is a contradiction which
proves $O(f)$ is minimal.
\enddemo

Using Propositions 1.2, 1.4 and [NS, p. 373-379], one can show
\proclaim {Definition 1} Let $f:\r \to X$, $f (\r)$ relatively
compact ($f\in B_{rc} (\r,X)$). Then $f$ is called recurrent if
$E(f,\e, F)$ is relatively dense in $\r$,  for all $\e > 0$ and
all $F$ finite subset of $\r$.

$E(f,\e, F)$ relatively dense means: there is a finite subset
$H\st\r$ such that $\r= H+E(f,\e, F)$.
\endproclaim
 \proclaim {Theorem 2} Let $f\in  B_{rc}(\r,X)$.
Then $f$ is recurrent if and only if $f$ is minimal.
\endproclaim

\proclaim {Theorem 1.6} (a) Let $\gamma_{\om}$ be any continuous
character  of $\r$ and $g \in C_{rc} (\r,X) \cap MIN (\r, X)$.
Then $(f,\gamma_{\om})$ is minimal.

\endproclaim
\demo{Proof} (a) By Birkhoff's Theorem  for  $(f, \gamma_{\om})$
[LZ, p.9], [BJM,  Proposition 1.6.12, p. 53] , there is a net
$\alpha =(a_i)\st \r$ such that
$(f_{\alpha},(\gamma_{\om})_{\alpha})$  is minimal. As $f$ is
minimal, $f_{\alpha}$ is minimal. There exists a net $\beta \st
\r$ such that $f=(f_{\alpha}) _{\beta}$ and $h
=((\gamma_{\om})_{\alpha}) _{\beta}= c\, \gamma_{\om}$  with
$|c|=1$. One has $(f,c \,\gamma_{\om})$ is minimal.  It follows
$(f,\gamma_{\om})$ is minimal.
\enddemo

\proclaim{Proposition} Let $f\in BC (\r,X)$ and for each $\e >0$
and each finite set $F\st \r$, there is a finite set $H \st \r$
such that $\r=E(f,\e,F)+H$. Then for each $\e >0$ and each compact
set $K\st\r$, the set $E(f,\e,K)$ is also relatively dense in
$\r$.
\endproclaim
\demo{Proof} Let $\e >0$,   $H$ finite, $\st \r$, $\r= H+
E(f,\e,F)$ and  $E(f,\e,F)\not = \r$ . We can assume $H=\{h_1,
\cdots, h_m\}$ with $0=h_1< \cdots < h_m=L$. Let $E(f,\e,F)=
\cup_{k\in \z} (\tau_k +[-\delta_k,\delta_k])$, where $\tau_0 \in
[-\delta_0,\delta_0]$ and $\tau_k < \tau_{k+1}$, $k\in \z$. We
show that for any $a\in\r$, the interval $[a, a+L]$ contains at
least one element of $E(f,\e,F)$. Indeed, let $a+L= h_{k(a)}+
\tau_a$ with $h_{k(a)}\in H$, $ \tau_a \in E(f,\e,F)$, then $ a
\le\tau_a = a+L-h_{k(x)}\le a+L$ implying $\tau_a\in  [a, a+L]
\cap E(f,\e,F)$. Now we show each interval $[a, a+L]$ contains at
least one element $\tau_k +[-\delta_k,\delta_k]$ with $\delta =$
min $_{1\le k \le m-1}\,\, (h_{k+1}-h_k)$.
 Assuming the contrary, let $[a,a+L]$

 $\delta_{k_0} <
\delta/2 $. Then $ \tau_{k_0}+ (h_k+h_{k+1})/2 \not \in
\tau_{k_0}+ [-\delta_{k_0}, \delta_{k_0}]+H$, $k=1, \cdots, m-1$.
$\square$
\enddemo
 \proclaim {Theorem 1.6} (a) Let $f\in AP(\r,X)$ and $g
\in RC(\r,X)\cap C_{rc}$. Then $(f,g)$ is recurrent.

(b) Let $f\in AA_u(\r,X)$ and $g \in RC(\r,X)\cap C_{urc}$. Then
$(f,g)$ is recurrent.
\endproclaim
\demo{Proof} (a) By Birkhoff's Theorem [LZ, p.9], there is a
sequence $\alpha =(a_n)\st \r$ such that $(f_{\alpha},g_{\alpha})$
is recurrent (= minimal). As $g$ is minimal, $g_{\alpha}$ is
minimal. There exists $\beta \st \r$ such that $g=(g_{\alpha})
_{\beta}$ and $h =(f_{\alpha}) _{\beta}$ is almost periodic.  One
has $(h,g)$ is recurrent. It follows that $E(h,\e,\r) \cap
E(g,\e,N)$ is relatively dense in $\r$. But
$E(h,\e,\r)=E(f,\e,\r)$.
 It follows that $E(f,\e,\r) \cap E(g,\e,N)$ is relatively
dense in $\r$ and hence $(f,g)$ is recurrent.

(b) See [Flor, Satz 4].
\enddemo

\proclaim {Lemma 1.7} (a) Let $f, g\in  C_{urc}(\r,X) \cap
RC(\r,X)$ and $(f,g)  \in RC(\r,X\times X)$. Then $ O(f+g) =\{F:
F= \alpha
 f+\alpha g, \alpha \in E_{X} (\r) \}$, where $E_{X} (\r)$ is the envelope
compact Ellis  semigroup  with $X= O((f,g))$.

(b) Let $f, g\in  C_{urc}(\r,X) \cap RC(\r,X)$ and $(f,g) \not \in
RC(\r,X\times X)$. Then $f+g \not \in
 RC(\r,X)$.

 (c)Let $f\in
C_{urc}(\r,X) \cap RC(\r,X)$. Let $g\in C_{urc}(\r,X)$ but $g\not
\in  RC(\r,X)$. Then $f+g$ or $f-g$  is not recurrent.
\endproclaim
\demo{Proof} (a) Since $\alpha \in E_{X} (\r)$, $\alpha
 (f+ g) =\alpha
 f+\alpha g$ implying $ O^*(f+g)=\{F:
F= \alpha
 f+\alpha g, \alpha \in E_{X} (\r) \} \st O(f+g)$. Since $E_{X} (\r)$ is compact, $
 O^*(f+g)$ is closed.  Since $ \{\alpha:  \alpha =r, r \in \r\}$
 is dense in $E_{X} (\r)$, $ O^*(f+g)$ is dense in  $O(f+g)$. It
 follows   $O^*(f+g)=O(f+g)$.

(b) By Birkhoff theorem there is $a=(a_n)\st \r$ such that $(g_a,
f_a, f_a+g_a)$ is recurrent. There is $b=(b_n)\st \r$ such that
$((g_a)_b, (f_a)_b, (f_a+g_a)_b)=(g, (f_a)_b, (f_a)_b
 +g)$ and  we can choose $b$  such that $((
 (f_a)_b)_a)_b = (f_a)_b$. Indeed, consider $\Cal {R}= \{b: (f_a ,g_a)_b=((f_a)_b ,g)
 \}$. Then $\Cal {R}\not = \emptyset$ and  $a\Cal {R} $  is a
 compact right topological
 subsemigroup  (see  [Flor, proof of Satz 1 (1 implies 3)]) of the compact envelope semigroup $E_{X} (\r)$ with $X= O(F)$, $F=
 (f,g)$:
  If $ab, ab'\in a\Cal {R}$, then  $ab ab'\in a\Cal {R}$ because $((((f_a)_b ,g)_a)_{b'}= ((((f_a)_b)_a)_{b'} ,g)$.
   By [BJM, Theorem 3.11, p. 31],
 there is  an idempotent $ab\in a\Cal {R}$, i.e.  $abab= ab$. It
 follows  $ ((f-(f_a)_b)_a)_b = (f_a)_b -(f_a)_b=0$. This implies
 $ f-(f_a)_b$  is not
 recurrent for $O(f-(f_a)_b)\not =\{0\}$ but $ O(((f-(f_a)_b)_a)_b)=\{0\}$ (would $f-(f_a)_b$ be recurrent,
 $O(f-(f
 _a)_b)=O((f-(f_a)_b)_a)=
 O(((f-(f_a)_b)_a)_b)=\{0\}$). We have $f+g= f+g+ (f_a)_b+g -((f_a)_b+g)=(f_a)_b+g+
 f-(f_a)_b$ and $(((f_a)_b+g)_a)_b+
  ((f-(f_a)_b)_a)_b =(f_a)_b+g$. Assuming $f+g$ is recurrent, one gets $O(f+g)=
  O((f_a)_b+g)$. One has  $((f_a)_b)_{\Cal{C}}+g_{\Cal{C}} \in O((f_a)_b+g)$ for
  some $c\st \r$ if and only if $(((f_a)_b)_{\Cal{C}}, g_{\Cal{C}}) \in O(((f_a)_b,g))$.
  Therefore  $f+g\not \in O((f_a)_b+g)$, a contradiction which proves
  that
  $f+g$ is not recurrent.

(c) Assume $f+g$ is recurrent. Since $g\not \in RC(\r,X)$, $O(g)$
is not recurrent (=minimal). By Birkhoff theorem there is
$a=(a_n)\st \r$ such that $(g_a, f_a, f_a+g_a)$ is recurrent.
There is $b=(b_n)\st \r$ such that $((g_a)_b, (f_a)_b,
(f_a+g_a)_b)$ is recurrent and $(f_a+g_a)_b)=f+g$. This implies
$g=(f_a+g_a)_b)-f$. Since $g$ is not recurrent, $((f_a+g_a)_b),
f)$ is not recurrent too. It follows $(-(f_a+g_a)_b), 2f)$ is not
recurrent. By (a),   $f-g= 2f-(f_a+g_a)_b)$ is not recurrent.
\enddemo
\proclaim {Theorem 1.8} Let $f,g\in  C_{urc}(\r,X) \cap RC(\r,X)$
and  $g\in L(f) $. Then

 (a) $(f,g)$ is  recurrent.

 (b) $\phi g \in L(f)$ for each $\phi \in AA_{u} (\r,\cc)$
\endproclaim
\demo{Proof} (a) Assuming $ F=(f,g)$ is not recurrent,  we get
$(f,-g)$ is not recurrent and hence  $f-g$ is not recurrent  by
Theorem 1.7 (a). This implies $f-g \not \in L(f)$ contradicting
$g\in L(f)$.

(b) By Theorem 1.6 (b), $((f,g), \phi)$ is recurrent. Hence for
each $a\st \r$ regular for $ F=((f,g), \phi)$ there is $b\st \r$
regular for $ F_a=((f_a,g_a), \phi_a)$ such that $
(F_a)_b=((f_a,g_a), \phi_a)_b = F$. This implies $G=(f, \phi g)$
is recurrent. Hence $f +\phi g $ is recurrent and hence $\in
L(f)$. It follows $\phi g \in L(f)$.

\enddemo
\proclaim{Theorem} Let $T(\cdot)\st L(X)$ be a $C_0$-group of
bounded operator such that $T(\cdot)x \in AP(\r,X)$ for each $x\in
X$. Let $f: \r \to X$ be weakly continuous and weakly relatively
compact. Then $\phi (t)= T(t)f(t)$ is weakly relatively compact
\endproclaim
\demo{Proof} Since $T(\cdot)x \in AP(\r,X)$, the mean $A_{\la}x
=\lim _{T\to \infty} \int_0^{T} T(t)x e^{i\la t} \, dt$ exists for
each $\la \in \r$.  Since $f$ is weakly continuous and weakly
relatively compact, we can consider $X$ is separable. Let $(x_n)
\st X$  be a dense countable subset. Then $\{y : y= A_{\la} x_n
\not =0, n \in \N, \la \in \r\}$ is countable. It follows $\{y :
y= A_{\la} x \not =0, x \in X, \la \in \r\}$ is countable too. It
follows that there exists $(\la_k) \st \r$ such that

(1) \qquad $T(t)x = \sum _{k=1}^{\infty} A_k x e^{i\,\la_k} t$,
where $A_k x = A_{\la_k} x$, $x\in X$.

\noindent It follows that for each $(t_{k'} \st \r$, there exists
a subsequence $(t_{k} \st (t_{k'} $ such that $T(t_n)x \to Ax$,
where $A\in L(X)$, $e^{i, t_n} \to e^{i\, \theta}$ and $f(t_n) \to
y$  weakly as $n \to \infty$.

Using (1) and Fourier series of almost periodic functions, one
gets $T(t)Ax = \sum _{k=1}^{\infty} A_k x e^{i\,\la_k}
(t+\theta)$.

\noindent Now, we show that $T(t_n) f(t_n)\to Ay$ weakly as $n\to
\infty$. Indeed,  one has

 $T(t)T(t_n) f(t_n) = \sum
_{k=1}^{\infty} A_k f(t_n) e^{i\,\la_k} (t+t_n)$. Since $A_k
f(t_n) \to A_k y$ as $n\to \infty$ for each $k \in \N$, one gets
$\sum _{k=1}^{\infty} A_k f(t_n) e^{i\,\la_k} (t+t_n)\to \sum
_{k=1}^{\infty} A_k y e^{i\,\la_k} (t+\theta)= T(t)Ay$ weakly as
$n \to \infty$.
\enddemo
\pagebreak

 inequality  for estimating $| \int_0^T \sin t^{\alpha}
\, dt|$, $1>\alpha >0 $:

Let  $y= t^{\alpha}$.  Then  $t = y^{1/\alpha}$ and

$I(T)= \alpha \int_0^T \sin t^{\alpha} \, dt=$

$  \int_0^{T^{\alpha}} [\sin y] {y}^{(1/\alpha)-1} \, dy
=\sum_{j=0}^{k-1} \int_0^{2\pi} [\sin y] [y+2j\pi]^{(1/\alpha)-1}
\, dy+ I_{\theta}$,

Where $0\le \theta < 2\pi$ and  $\theta + 2k\pi = T^{\alpha}$ and
$I_{\theta}= \int_{2k\pi} ^{T^{\alpha}}\, [\sin y]
(y+2k\pi)^{(1/\alpha)-1} \, dy$.

One has
 $|I(T)|\le  |\int_0^{\pi}
\sin y ( y^{(1/\alpha) -1}- (y+\pi)^{(1/\alpha) -1}+\cdots+$

$(y+(2k-3)\pi)^{(1/\alpha) -1}- (y+ 2(k-1)\pi)^{(1/\alpha) -1})\,
dy| + |I_{\theta}|\le $

 $\int_0^{\pi} [y^{(1/\alpha) -1}+(y+ (2k-1)\pi)^{(1/\alpha) -1}]\, dy + |I_{\theta}|\le (\pi)^{1/\alpha}+ T^{\alpha} +|I_{\theta}|$.

% \head{L-ap Levin,  Levitan-Zykov  Definitions of N-ap }\endhead

 Choose $\e
>0$,   non-negative $f\in \f(\r)$  such that $ sp^{w\Cal{L}}(\phi) \cap [-\e,\e]
=\emptyset$, supp $\hat{f} \st [-\e,\e]$  and $\hat{f} (0)=1$. By
part (i), we conclude $sp^{w\Cal{L}}(\phi*f)=\emptyset$. By [Y,
Proposition 1,  p. 156] we conclude $\phi*f\in BUC(\r,X)$ and
hence by Theorem 2.3 (iii) $\phi*f\in C_0(\r,X)$. Consider the
sequence $g_n (\cdot)= (1/n) \psi(\cdot/n)$ of Proposition 1.2
(ii). We prove that $\phi* g_n\to 0$ as $n\to \infty$ in
$BUC(\r,X)$. Indeed, one has $\phi* g_n
(t)=\int_{-\infty}^{t/n}\phi(t-sn)\psi (s)\, ds=\int_0
^{\infty}\psi (\frac{t}{n} -s)\phi (ns)\, ds $. Therefore,
$||\phi* g_n||_{\infty}\le ||\phi||_{\infty} $ and $||\phi* g_n
(t+h)-\phi* g_n (t)||\le ||\phi||_{\infty} ||\psi_{h/n}
-\psi||_{L^1} $.

\proclaim {Theorem} Let $f\in AP(r,X)$ and $g \in RC(r,X)\cap
C_{rc}$. Then $(f,g)$ is recurrent.
\endproclaim
\demo{Proof} By Birkhoff's Theorem [LZ, p.9], there is a sequence
$\alpha =(a_n)\st \r$ such that $(f_{\alpha},g_{\alpha})$ is
recurrent (= minimal). As $g$ is minimal, $g_{\alpha}$ is minimal.
There exists $\beta \st \r$ such that $g=(g_{\alpha}) _{\beta}$
and $h =(f_{\alpha}) _{\beta}$ is almost periodic.  One has
$(h,g)$ is recurrent. It follows that $E(h,\e,\r) \cap E(g,\e,N)$
is relatively dense in $\r$. But $E(h,\e,\r)=E(f,\e,\r)$. This
implies. It follows that $E(f,\e,\r) \cap E(g,\e,N)$ is relatively
dense in $\r$ and hance $(f,g)$ is recurrent.
\enddemo

(5)\qquad  $m_{\la}(T)= (1/2T)\int _{-T}^T f(t) e^{-i\la \,t}\,
dt$, $c(\la)=\lim_{T\to \infty} m_{\la}(T)$

\qquad \qquad  and $\e_{\la}(T)= m_{\la}(T)-c(\la)$.

One has

$\e_{\la}(T)=(1/2T)\int _{-T}^T [f(t) e^{-i\la\,t}- c(\la)]\, dt=
(1/2)(\e^{+}_{\la}(T)+\e^{-}_{\la}(T))$

with

\qquad $\e^{+}_{\la}(T)= (1/T)\int _0^T [f(t) e^{-i\la\,t}-
c(\la)]\, dt$ \,\,\,and

\qquad  $\e^{-}_{\la}(T)= (1/T)\int _0^T [f(-t) e^{i\la\,t}-
c(\la)]\, dt$.

Since $f(t) e^{-i\la\,t}- c(\la)=\frac{d}{dt} (t\e_{\la}^+ (t))$
and $f(-t) e^{i\la\,t}- c(\la)=\frac{d}{dt} (t\e_{\la}^{-} (t))$,
one has

$m_{\la\pm \mu}(T)$

$:= (1/2T)\int _{-T}^T c(\la) [e^{-i\mu t}+e^{i\mu t}]\, dt +
(1/2T)\int _{-T}^T [f(t) e^{-i\la\,t}- c(\la)][e^{-i\mu t}+e^{i\mu
t}]\,dt $

$=  c(\la) (2\sin {\mu T})/(\mu
 T)
  +
(1/T)\int _{0}^T [\frac{d}{dt} (t\e_{\la}^+ (t))+\frac{d}{dt}
(t\e_{\la}^{-} (t))]\cos {\mu t}\,dt $

$= c(\la) (2\sin {\mu T})/(\mu
 T)
 +(1/2T)\int _{0}^T  \frac{d}{dt}
(t\e_{\la} (t))]\cos {\mu t}\,dt$.

$= c(\la) (2\sin {\mu T})/(\mu
 T)
 + (1/2)\e_{\la} (T)) \cos {\mu T} -(\mu /2T)\int _{0}^T (t\e_{\la} (t))]\sin {\mu
 t}\, dt$.

Now with $m_{\la\pm \mu}(T)$ one can prove (6) of Kahane and may
be can finish the proof.

Instead of (11) one gets

(11*) \qquad sup$_{T  \in  U_j} ||m_{u,v}(T)||>5\e/6$,

          inf $_{T  \in  U_j} ||m_{u,v}(T) || < \e/2$

      if  $u = \lambda_j, v  \in  \r$  with $u+v  \in  I_j  \cap F_{\e}$,

    where  $m_{u,v} : =  m_{u+v} + m_{u-v}$    (your $m_{u\pm v}$).

With unique  $y  \in  \cap ^{\infty}_1  K_j$, $ K_j : =$ closure
$\overline{I_j \cap F_{\e}}$, and $y = \lambda_j + v_j $ with
unique $v_j$, one gets as in my email Kahane5

(12) \qquad  sup $_{U_j} ||m_{\lambda_j,v_j}|| \ge 5\e /6$
         inf $_{U_j} ||m_{\lambda_j, v_j} ||  \le  \e /2$ ,
            $j  \in  \N$, with $\lambda_y + v_j = y $.

So this unfortunately does  not  give a contradiction, because we
do not know how  $m_{\la _j - v_j}$ behaves ; it is true that
$\lambda_j - v_j $ converges also to $y$, but this does not seem
to help

Do you have any ideas here

 To our best knowledge
examples of $BAA(\r,\cc)$ given in (see [V], [T], [M]) are uniformly continuous. This led  to the need
of the following:

\proclaim {Example 1} If $p(t)=2+\cos t + \sqrt{2} t$ then $g=
\sin(1/p)\in BAA$ and $g\in S^1$-$AP $.
\endproclaim
\demo {Proof} For $(t_{n'})\st \r$  there is a subsequence
$(t_{n})$ such that $p(\cdot+t_n) \to  q(\cdot)$.  The  set $
C=\{s\in \r: q(s)=0\}$ is at most  countable. Passing to a
subsequence if necessary, we can assume
  $\sin(1/p(t+t_n)) \to a(t)$ for all $t \in \r$ with
 $a(t)= sin (1/q(t))$, $t\not \in C$ and $|a(t)|\le 1$.
One can check that $a(t-t_n)\to g(t)$, $t\in \r$. This proves $g\in BAA$.

Now consider $ f_n (t)= \sin\, \frac {1}{2+ \max\, \{\cos\, t,
-1+1/n\} + \sqrt {2} t}$. Then $||g-f_n||_{S^1}\le \mu\{t\in
[0,\pi]: |1-\cos\, t |\le 1/n \}\to 0$ as $n\to \infty$. Since $2+
\max \{\cos t, -1+1/n\} + \sqrt {2} t \ge 1/n$, then $f_n \in
AP\st S^1$-$AP$. This implies $g\in S^1$-$AP$.
\enddemo

\proclaim {Example 2} An $f\in \m AP(\r,\r) \cap C(\r,\r)$, $-1
\le f \le 1$,
 with $Pf\in  AP$ but $f\not \in S^1$-$AP$ and $f\not \in  AA$, is given by

$f:= \sum_{n=2}^{\infty} \, h_n$, with $h_n \in C(\r,[-1,1])$,
$h_n \equiv 0$ on $[- 2 ^{n-1}, 2 ^{n-1} -1 ]$,

 $h_n =\sin (2 ^{n-1} \pi (t-2 ^{n-1}))$  on $[ 2 ^{n-1}-1, 2
^{n-1}]$,

$h_n$ period $2 ^{n}$.

Since $\tau =2^{n+1}$ is not $\e$-period for $h_{n+k}$, $k \ge 2$,
$0 < \e\le 1/2$, it follows $E(f,\r,\e)$ is not r.d. Also, one can
prove that $E(f,N,\e)$ is not r.d. Hence $f$ is not aa or
recurrent function.
\endproclaim
 \proclaim {Proposition} Let $f\in C_u  (\r,X)\cap
Rec(\r,X)$. Then $f\in BUC  (\r,X)\cap Rec(\r,X)$.
\endproclaim
\demo{Proof} Let $\e >0$. Then $E(f,\e, [-1,1])$ is relatively
dense in $\r$. Choose $L$, $1 >\delta >0$ such that

(1.1)\,\,$\r= [L,L]+E(f,\e, [-1,1])$ and $||f(t+h)-f(t)||\le \e$
if $|h|\le \delta$.

\noindent Set $K_s=[-s,s]+E(f,\e, [-1,1])$, $a=$ sup$_{t\in
[-1,1]} ||f(t)||$.
 One has sup $_ {t\in K_1} ||f(t)||\le a+\e$.
  By (1.1), one can verify sup $_ {t\in K_{n\delta}} ||f(t)||\le (1+ n)\e +a$.
Let $t\in \r$. $t= x+\tau$ with $x\in [-L,L]$ and $\tau\in
E(f,\e, [-1,1])$. Let $|x|= n\delta +h$, $h < \delta $. Then $n
\le  L/\delta$.
 This implies
$ ||f(x+\tau) ||\le (1+ (L/\delta))\e+ a+\e$. This proves the
statement
\enddemo

\proclaim {Definition (LZ)} For $f  \in  C(R,X)$, $f$ is called
N-ap (Levitan-Zhikov book) if to each $\e >0$ and $ N>0$ there
exist a set $ U(\e,N)  \st  \r$, such that the collection of these
sets satisfies :
                all $U$ are r.d.(Bohr),  $U(\e,N)  \st  E(\epsilon;N)$ for all $\e, N$;

to each $\e > 0$ and $N$ there is $\delta > 0$ such that

                                $ U(\delta,N) \pm U(\delta,N)  \st  U(\e,N)$,

     where  $E(\e,N) : = \{ p \in R : ||f(t+p)-f(t)||<\epsilon for |t|\le N \}$.
\endproclaim
\smallskip

Question :   where is a proof that such N-ap f satisfy

(*)     For each $\e > 0$ and $N > 0$  there are  $\delta > 0$  and  $L > 0$
         such that     $E(\delta,L)-E(\delta,L)  \st  E(\e,N)$

.....................................

\smallskip

From Levin's Paper [Levin, 1949]

\proclaim{Definition 1} $f\in BC(\r,\cc)$ is Levitan's ap (L-ap)
if

(i) For  each $\e > o$, $N >0$ there is  r. d. (Bohr) $E_{\e,N}$, where  $E_{\e,N}$ is a set of intervals of $\r$ and
$E_{\e,N}\st E(f, N,\e)=\{p: |f(t\pm p)-f(t)| <\e, |t| <N\}$;

(ii)For each $\e$ and fixed $N>0$ there is $\delta_{\e}>0$ such that

$E_{\delta,N}\pm E_{\delta,N}\st E_{\e,N}\,$ for all $\,\delta <\delta_{\e}$.
\endproclaim

The remark says $ E_{\e,N}$ does not necessary  contain all $E(f, N,\e)$.
See the example of Levin below.

After that Levin proves:

Levin Definition with conditionally convergent sequences as in Reich or [Basit, 77] implies (i), (ii) with

(i')$\qquad  E_{\e,N}=E(f, N,\e)$, $N\ge 1/\e$; \qquad
$ E_{\e,N}=E(f, 1/\e,\e)$, $N<  1/\e$.

The proof that (i), (ii) implies Levin's ap:

First Levin showed that each  $E_{\e,N}$ contains an interval $(-r, r)$. This implies $E_{\e,N}$ is relatively
dense in the sense of Veech because  there is $\delta> 0, r>0$ such that  $(-r,r)+ E_{\delta,N}\st E_{\e,N}$.  This answers your questions.

\smallskip

The proof of Levin :

Translation: We claim that $E_{\e,N}$ contains an interval with center $0$. Indeed, there is $\rho >0$ such

$E_{\rho,N}+E_{\rho,N} \st E_{\e,N}$,

but $E_{\rho,N}$ consists itself of difference of intervals and $I-I$ with $I$ interval is an interval with center $0$.

.......

\smallskip

To generalize this proof  to $U(\e,N)$ it is enough to Show that $U{\e,N}$ has at lest an interval. But I can not prove that from your
(*). After thinking a long time I read again Def.1  and I discovered that Def. 1 states that $E_{\e,N}$ states that it consists
of intervals from $R$ as translated above. But (*) does not require that. I believe that (*) is equiv to Def. 1 but
 the proof
is not direct:

$f$ with (*) implies $f$  satisfies Def.3 of [LZ, p. 54] implies  $f$ satisfies Levins  implies $f$ satisfies (i), (ii) with
$E_{\e,N}$ given by (i') above.

\smallskip

Example of Levin: (the following is a translation but I think one should add $\phi, \phi_1$ are periodically extended)

Let $\phi (t)=1$ if $|t-n\omega|<\delta$, ($n\in \N_0$, $\delta =\omega/10$), $\phi=0$ else;

$\phi_1 (t)=1$ if $|t-n p\omega-\omega/2|< \delta$, $n\in \N$, $p$ some natural number.

 $f =\phi+\phi_1$ is periodic with period $p\omega$:

$\tau_n =n\omega$, $\tau'_n =np\omega+\omega/2$ are $0,\delta$-periods ($\in E(f,\delta,0)$) but

$\tau'_n-\tau_{np}=\omega/2\not \in E(f,\delta,0)$.

One can also change the above slightly to make f continuous.

\bigskip

................................

Let $X$ be a separable Banach space and  $ \{a_n: n\in \N\}$  a dense $\Q$ subspace of $X$.
 Denote by $\n=\{F\st  \{a_n: n\in \N\}:  F$ is finite $\}$. Then $\n =\{F_n: n\in \N\}$.
If $F\in \n$, $\e >0$, $F_{\e}= \cup_{a\in F} \{x\in X: ||x-a||\le \e\}$.

Let $(x^*_n)$ be a dense subset of  $X^*$ endowed with its $X$-topology (see [DSch, p. 425]).  Set (see [DSch, p. 434], Theorem 3)

(1.1)\qquad \qquad $d(x,y)= \Sigma_{k=1}^{\infty} \frac{|x_k^*(x-y)|}{2^k(1+|x_k^*(x-y)|)}$,

(1.2)\qquad \qquad $B_{\e} (x)=\{y\in X: \,\, d(y,x)\le \e\}$.

Then $B_{\e} (x)$ is weakly closed in $X$. For $F\st X$, we define $F(d,\e)= \cup_{a\in F} B_{\e} (a)$.

Let $K(n)= \cap _{j=1}^n ker\, {x_j^*}$,  $ker\, {x_j^*}=\{x\in X: x_j^* (x)=0\}$. Then $X/K(n)$ is finite dimensional.
Therefore,
there exists a finite dimensional subspace $X(n)\st X$ with $X(n)\cap K(n) =\{0\}$ and $X=X(n)\oplus K(n)$. It follows

(1.3) \qquad \qquad $B_{\e,n} (x) : =B_{\e} (x)\cap X(n) $ is
compact for each $n\in \N$, $x\in X$, $\e >0$.

(1.5) \qquad \qquad  $V_{\e,n} (x) : =B_{\e} (x)\cap K(n) $

One has $B_{\e} (x)=B_{\e,n} (x) + V_{\e,n} (x) $ and $co (B_{\e} (x))= co (B_{\e,n} (x)) + co (V_{\e,n} (x)) $.

For $\e >$ choose $n _{\e}\in \N$ such that $\sum_{> n_{\e}} (1/2^k)<\e$. It follows

(1.5) \qquad \qquad $co (V_{\e,n_{\e}} (x))\st K(n_{\e}) $ and
$d(x,0)< \e$, $x\in K(n_{\e})$.

\proclaim {Lemma 1}  To $F\in \n$, $\e >0$ and $\delta >0$ there is
$F', F''\in \n$ with $\overline{co\, F{d,\e}}\st
F'(d,\e+\delta)$ and $\overline{co\, (F(d,\e)- F(d,\e))}\st F''(d,2\e+\delta)$.
\endproclaim
\demo{Proof} Choose $n= n(\e)$ such that $\sum_{k=n(\e)+1}^{\infty} (1/2^k) \le \e/2$. Consider $X_n=\{x\in X: x^*_k(x)$ for all $k > n(\e)\}$.
Then $X_n$ is finite dimensional  and
hence  $ F{d,\e}= F_1 +F_2$ with $ F_1 \st X_n$ and $F_2\st Y$, where $X=X_n\oplus Y$.

\enddemo
\demo{Proof} Let $F={a}$. Then $I(a) =\{ta$, $0\le t
 \le 1\}$ is compact. Choose $F'\in \n$ such that $I(a) \st F'_{\delta}$. We claim  $ co B_{\e} (a) \st F'{d, \e +\delta}$. Indeed, let $ x=\sum _{i=1}^m  u_i h_i \in co B_{\e} (a)$,
$h_i \in  B_{\e}$, $0\le u_i \le
 1$, $1\le i\le m$, $\sum_{i=1}^m  u_i =1$.
\enddemo
\proclaim {Theorem 2}  For any $\jj$, $X$,

(i) $C_{wrc} (\jj,X)$ satisfies  $C_{wrc}\st \m C_{wrc}$, $(\Delta)$,

(ii)$L^{\infty}_{rc} (\jj,X)$ satisfies  $L^{\infty}_{rc}\st \m L^{\infty}_{rc}$, $(\Delta)$.

\endproclaim
\demo{Proof}(i) Assume $f\in L^1_{loc}$, $\Delta_h f \in  C_{rc}
$, $h >0$. By [], $f\in C$. One can assume that $X$ is generated
by $f(\jj)$ and therefore separable. Choose $\n=\{F_n: n\in \N\}$
as above. Denote by $A_m^{\e}=\{\phi \in BC: \phi (\jj)\st
(F_m)_{\e}\} $. Define $V_m= \{h\in [0,1]: \Delta_h f\in
A_m^{\e}\}$. Then $V_m$ are closed and $[0,1]=\cup_{i=1}^{\infty}
V_m$. By Baire category theorem, there is $m_0 \in \N$, $0 <
u+\rho <1$ with $[u, u+\rho] \st V_{m_0}$. This means that
$\Delta_ h f  (\jj)\st (F_{m_0})_{\e}$, $u \le h \le u+\rho$.
Since  $(\Delta_ k f)_u (t)= \Delta_ {u+k} f (t)- \Delta_u f (t)$,
$(\Delta_ k f)_u (\jj)\st (F_{m_0})_{\e}-(F_{m_0})_{\e}$, $0 <k <
\rho$.
 By
$(1/k) \int_0^{k} [f(t+u+s)-f(t+u)]\, ds = lim_{n \to \infty}(1/nk)\sum_{j=1}^n [f(t+u+ j k/n)-f(t+u)]\in
\overline {co ((F_{m_0})_{\e}- co ((F_{m_0})}$, $0<k <\rho$.
By Lemma above and since $\e$, $\delta$ are arbitrary,  one gets $(1/k) \int_0^{k} [f(\cdot+s+u)-f(\cdot+u)]\, ds$
$ \in C_{rc}$, $0< k< \rho$. By $f-M_ kf -(f-M_k f)_u = (f-f_u)- M_k (f-f_u)$ and $M_k (f-f_u)\in C_{rc}$, one concludes
$f-M_k f -(f-M_k f)_u\in C_{rc}$ implying
  $f-M_{k} f \in C_{rc}$, $0< k< \rho$.
This implies $f- M_h f\in C_{rc}$ for all $h>0$ and proves $(\Delta)$.

(ii) $L^infty_{\Cal{C}}r : = \{f \in L^{\infty} :$ to $f$ exists a
$L$-nullset $P$ with $f(\jj-P)$ relatively compact $\}$.

In the def. of $A^{\epsilon}_m$  we replace the $phi(\jj) \st (F_m)_{\epsilon}$
by  $phi(\jj) \st (F_m)_{\epsilon}$ a.e., meaning to $\phi$ exists a $L$-nullset $P st \jj$
with
$\phi(\jj/P)  \st  (F_m)_{\epsilon} $.
One gets then :  For each  $k  \in  [0,ro]$ there exists a $L$-nullset $P_k$ such
that
$(\Delta_k f)_u (\jj-P_k)  \st  (F_{m_0})_{\e} - (F_{m_0})_{\e}$.
Now  $F(k,t) : = ((\Delta_k f)_u) (t)$ defines a $(k,t)$ measurable
$F:[0,ro] \times \jj \to X$.
Then for any open $G \st $X,  $F^{-1}(G)$ is a measurable subset of $[0,ro] \times \jj$
by [Gun, p. 192,
Aufgabe 126, "Beispiel"].
With Fubini-Tonelli for measurable sets one gets , with     $ G : =  X  -
\overline {co}((F_{m_0})_{\e} - (F_{m_0})_{\e})$, that  $F^{-1}(G)$ is
($2$-dimensional) nullset and so there exists a nullset $P \st \jj$ such that for $t
\in \jj/P$
one has  $F(k,t)  \in  X / G$ for $k \in [0,ro] /U_k$,   $U_k$ a nullset.
Then one can proceed as in the proof of case (i), replacing the $C_{rc}$
everywhere
by $L^infty_{rc}$  as defined above .
 \P
\enddemo

An $f\in \m AP(\r,\r) \cap C(\r,\r)$, $-1 \le f \le 1$,
 with $Pf\in \m AP$ but $f\not \in  AP$ and $f\not \in  AA$, $f \not \in S^1 AP (\r,\r)$, is given by

$f:= \sum_{n=2}^{\infty} \, h_n$, with $h_n \in
C(\r,[-1,1])$, $h_n \equiv 0$ on $[- 2 ^{n-1}, 2 ^{n-1}
-1 ]$,

 $h_n =\sin (2 ^{n-1} \pi (t-2 ^{n-1}))$  on $[ 2 ^{n-1}-1, 2
^{n-1}]$ ,

$h_n$ period $2 ^{n}$.

Since $\tau =2^{n+1}$ is not $\e$-period for $h_{n+k}$, $k \ge 2$, $0 < \e\le 1/2$, it follows $E(f,\r,\e)$ is not r.d.
Also, one can prove that $E(f,N,\e)$ is not r.d. Hence $f$ is not aa or recurrent function.

\proclaim{Theorem 3.1} If $\A  \st  C_u(\jj,X)$  is $\r$-linear,
positive-invariant and uniformly
           closed, then
           $\A$
           satisfies  $(\Delta)$.
\endproclaim

\demo{Proof} Denote by $\A_{ub} =\A \cap BUC(\jj,X)$. Then
$\A_{ub}$ satisfies $(\Delta)$ by [7, Proposition 3.1]. Let $f\in
\A$. Then

(3.5) \qquad $f= (f-M_h f) +M_h f =u+x +Pv =u+w,\,\,\,\,$  $h> 0$

\noindent with $u =f-M_h f$, $v =(\Delta_h f)/h \in \A_{ub}$ by
[7, Proposition 2.2], $x\in X$ and $w = x+ Pv \in \A$. Now we show
that $\A$ satisfies

$(\Delta P)$ \qquad $f\in \A$ implies $Pf-M_h Pf \in \A$ for all
$h> 0$.

\noindent Since $(Pu)_h-Pu= h M_h u \in \A_{ub}$ for all $h> 0$
and $\A_{ub}$ has $(\Delta)$, it  follows $\A_{ub}$ has $(\Delta
P)$. We prove $Pw-M_h Pw\in \A$. We have

$M_h Pw- Pw= hM_h w - (1/h) \int_0^h s\, w(\cdot+s)\, ds= hM_h w -
(h/2) w(\cdot+h) $

$+ (1/h) \int_0^h (s^2/2)v(\cdot+s)\,ds=w_1+w_2+w_3$.

\noindent One has $w_1\in \A$ by [7, Proposition 2.2], $w_2\in \A
$ because $w\in \A $ and the properties of $\A$ and $w_3 \in
\A_{ub}$ because it is a uniform limit of a sequence from
$\A_{ub}$. This proves $(\Delta P)$ for $\A$.

Assume now $f\in L^1_{loc} (\jj,X)$ and $\Delta_h f \in\A$ for all
$h>0$. Then $f= (f-M_h f) + P((\Delta_h f)/h) +x = u+Pv$ with $u=
f-M_h f+x \in UC$ because $UC(\jj,X)$ has $(\Delta)$ ([7, ]),
$\Delta_k u = \Delta_k f - M_h (\Delta_k f) \in \A$ by [7,
Proposition 2.2]. By the same method used in the proof of [7,
Proposition 3.1], one gets $u-M_h u \in\A$ for all $h>0$. Since
$v\in \A$ and $\A$ has $(\Delta P)$, one gets $Pv- M_h Pv \in \A$
for all $h>0$. Together, $f-M_k f\in \A$ for all $k>0$.
\enddemo

\proclaim {Examples 3.8}  $C_{w_N,ub}$, $C_{w_N,0}$, $AP_{w_N}$
and
     $AP_{[N]}$ satisfy    $(\Delta)$, for any $\jj$, $X$, $n\in \N$.
\endproclaim

 \demo{Proof} Case $C_{w_N,ub}$:
Let $f\in L^1_{loc} (\jj,X)$ with $F(h):=\Delta_h f\in
C_{w_N,ub}=Y$.
  With  $h(f (t) -M_h f(t)) = -  \int_0^h F(s)(t)\, ds = u(t)$ we have to show $
  u\in Y$. Now by Proposition 1.5 (ii), $F\in C(\r_+,Y)$, so by
 [41, Theorem 1, p. 133] the $L^1_{loc} (\r_+,Y)$-Bochner integral
 $v:= \int_0^h F(s)\, ds$ exists and belongs to $Y$.

For fixed $t  \in \jj$, define $Ty : =y(t) $; then $T : Y \to X$
is a bounded linear operator, so by a special case of Hille's
theorem [41, p. 134, Corollary 2] one has
    $T(v)  =  L^1_{loc}(\r_+,Y)$-Bochner integral  $\int^h_0 F(s)(t) ds, =
u(t)$. This gives   $u = v  \in  Y $ as desired.

In the same way one can show $(\Delta)$ for $ C_{w_N,0}$
respectively $AP_{w_N}$   with $Y = C_{w_N,0}$ respectively $
Y=AP_{w_N}$.

Since $(\Delta)$ for the $AP_{[N]}$ is not used here, we omit the
proof.

First we show $(\Delta)$ for $ AP_{[1]} (\jj,X)$. So,
  let $\Delta_h f \in AP_{[1]}$, $h>0$. As $AP_{[1]} \st
AP_{w_1}$, one gets $u=f- M_h f \in AP_{w_1}= t\cdot AP \oplus
C_{w_1,0}$. Since $AP_{[1]}$ is invariant under convolution,
$\Delta_k u =\Delta_k f+ M_h (\Delta_k f) \in AP_{[1]}$, $k>0$. We
claim that $u- M_h u \in AP_{[1]}$, $h>0$. Indeed, one has $u =t
F+
 G$ for  a unique $ F\in AP$ and a unique $G$ satisfying $G/w_1
\in C_0$. Moreover with  $\Delta_k u \in AP_{[1]} $, one gets
$\Delta_k G \in AP$, $k
>0$. It follows $G-M_h G\in AP $ by $(\Delta)$ for $AP$. Also,
 $tF -M_h (t F)=t (F -M_h ( F)) -
(1/h)\int_{0}^{h} s F(\cdot +s)\, ds  \in AP_{[1]} $  by
$(\Delta)$ for $AP$ and application of theorems on Bochner
integral  cited above to $\int_{0}^{h} s F(\cdot +s)\, ds$. This
proves our claim.

As in Theorem 3.1, it remains only to prove $(\Delta P)$ for
$AP_{[1]}$: Let $\phi =t p+ q \in AP_{[1]}$. Then $P\phi -M_h
P\phi = -[(1/h)\int _0^h \int_0^s q (\cdot+x) \,dx\,ds + t\cdot
(1/h)\int _0^h \int_0^s p(\cdot+x) \,dx\,ds + (1/h)\int _0^h
\int_0^s x\,p(\cdot+x) \,dx\,ds] \in AP_{[1]}$.

Now assume that the statement is true for $n=N$ and let $\Delta_h
f \in AP_{[N+1]}$, $h>0$. Then $u=f- M_h f \in AP_{w_{N+1}}=
t^{N+1}\cdot AP \oplus C_{w_{N+1},0}$ by $\Delta$ for $
AP_{w_{N+1}}$. It follows $u = t^{N+1} F+ G$ with $F \in AP $, $G
\in C_{w_{N+1},0}$ and $\Delta_k G \in  AP_{[N]}$, $k
>0$. By the assumptions of the induction it follows $G-M_h G \in
 AP_{[N]}$. Directly  as in the case $n=1$ one can conclude $
t^{N+1} F- M_h (t^{N+1} F)\in AP_{[N+1]}$ and $(\Delta P)$ for $
AP_{[N+1]}$.
 \P
\enddemo

\bf{Remark}. The Examples 3.8 hold for any continuous $w : \jj \to
(0,\infty)$ with \,\,\,
     sup $\{||w_h/w||_{\infty} : 0<h<=1\} < \infty $.

Case $\jj=\r$ is similar

***By (1.10) we need to prove $\om_0 \not \in sp^{\Cal{C}}(\phi)$
implies $\om_0 \not \in sp^{\Cal{C}_u} (\phi)$. By Proposition 1.1
(ii) and Corollary 2.3, one has $sp^{\Cal{C}} (\phi)=\cup_{h>0}
sp^{\Cal{C}} (M_h \Phi) =\cup_{h>0} sp^{\Cal{C}_u} (M_h \phi)$.
So, if
 $\om_0\not  \in sp^{\Cal{C}}(\phi)$, one gets  $\om_0\not  \in
sp^{\Cal{C}_u}(M_h\phi) $ for all $h>0$. By (1.7), one has  $\Cal
{C}(S(s)M_h \phi) (\la) =$
  $g(\la\, h) \Cal {C}\phi_s (\la) - (1/h)\int_0^h (e^{\la \, \tau} \int_0 ^u e^{-\la \, t}\phi(s+t+\tau)
 dt)\,d\tau$. Choose
$h_0
>0$ such that $g(i\om_0\, h_0)\not =0$. Since $i\,\om_0$ is a
regular point for $\Cal {C} ( S(\cdot)M_{h_0}\phi)$, it follows
$i\,\om_0$ is a regular point for $ \Cal {C}_u \phi$.

Given an $f \in PS(\r,X)$ and $g = e^{i \omega\, t}$, $\omega$
fixed $\in \r$; to  $\e_n = 2^{-n}$ one can then find  inductively
sequences of intervals $I_n=[-a_n,a_n]$ and $r_n >0$ ($a_n >0$)
with , for $n \in \N$ , $ |f_{r_n}  - f | < \e_n$  on $I_n$, $I_n
+ r_n  \st  I_{n+1}$, $a_n\to \infty$, $r_n\to \infty$; by taking
a subsequence one can further assume , with complex $c$, $|c|=1$,
$| g(r_n) - c | < \e_n$  for  $n \in \N$. Since $I_n+r_n  \st
I_{n+1}$, the above gives for $n+1$,  $| f_{r_{n+1} + r_n} -
f_{r_n} | < \e_{n+1}$  on $I_n$, and inductively
$|f_{r_n+...+r_{n+k}} - f |  <  \e_n +...+ \e_{n+k}$  on $I_n$,
$n,k \in  \N$, $|g_{r_n+...+r_{n+k}}  -  c^{k+1} g | <
\e_n+...+\e_{n+k}$ on $\r$. Now with the Kronecker approximation
theorem there exist $k$ ($\to \infty$ even) with  $| c^{k+1}  -  1
|  <$ given pos. $\e$, so for such a $k$ with  $ |
g_{{r_n}+...+r_{n+k}} - g | < \e _n +...+\e_{n+k}+\e$ on $\r$, $n
\in \N$. Since $a_n \to \infty$, $r_n \to \infty$
$\e_n+...+\e_{n+k} \to 0 $ as  $n\to \infty$, the above shows that
for $f$, $g$ as above, $(f,g) \in PS(\r,X  \times \cc$) ;
inductively $(f,g_1,...,g;m) \in PS$; since PS is uniformly closed
one gets the

Theorem. $f \in PS(\r,X)$, $ g \in AP(\r,Y)$ implies  $(f,g) \in
PS(\r,X\times Y)$.

Corollary.  For any $X$, $PS_b(\r,X)$  satisfies $(\Gamma)$.

With this, analogues of §§3-5 of RSD follow, even for $PS_b$
instead of $REC_{urc}$ , resp $rec_{urc}$.

\enddocument

In the following we give alternative equivalent definition of
Laplace spectrum. As is customary used, we define convolution  in
$L^1(\r_+)$ by $f*g (t) =\int_0^t f(t-s)g(s)\, ds$ for  $f, g \in
L^1(\r^+)$. Then for $\phi\in L^{\infty} (\r_+,X)$ and $f \in
L^1(\r^+)$, the convolution $\phi*f (t)=\int_0^t \phi(t-s)f(s)\,
ds$ is uniformly continuous and  bounded.  Extending $f\in
L^1(\r_+)$ by $0$ on $(-\infty,0)$ and $\phi\in L^{\infty}
(\r_+,X)$ by $0$ on $(-\infty,0)$, one gets $\int_0^t
\phi(t-s)f(s)\, ds=\in_{-\infty}^{\infty}\phi(t-s)f(s)\, ds$. So,
for $\phi\in L^{\infty} (\r_+,X)$ we may define

\qquad $I^{\jj}(\phi)=\{f\in L^1 (\jj): \phi *f=0\}$.

Define $sp^{\jj}(\phi)= \{\om\in \r: \phi*f =0:\text {for\,} f\in
I^{\jj} \}$.

\proclaim{Proposition 1.4}  Let $\phi\in L^{\infty}(\jj,X)$. Then
for $\jj=\r_+$, one has $sp^{\r_+} (\phi)= sp^{\Cal{L}} (\phi)$
and for $\jj=\r$, one has $sp^{\r} (\phi)= sp_B (\phi)$.
\endproclaim
 \demo{Proof} Consider first, the case $\jj=\r_+$. Assume $\om\in sp^{\Cal{L}}
(\phi)$. Let $f\in I^{\r_+} (\phi)$. Then $\phi*f=0$ and so
$\Cal{L}\phi*f (\la)=\Cal{L}\phi (\la)\Cal{L}f (\la)=0$ for
$\la\in\cc_+$. This implies $\lim_{\la\to \om}\Cal{L}\phi
(\la)\Cal{L}f (\la)=0$. Since $\lim_{\la\to \om}\Cal{L}\phi (\la)$
does not exist  but $\lim_{\la\to \om}\Cal{L}f (\la)= \hat{f}(\om)
$, it follows $\hat{f}(\om)=0$. This implies $\om \in sp^{\r_+}
(\phi)$. Now,  assume $\om\not \in sp^{\r_+} (\phi)$. There is
$f\in I^{\r_+} (\phi)$ such that $\hat{f}(\om)\not =0$. We have

$\lim_{\la\to \om}\Cal{L}\phi*f (\la)=\lim_{\la\to \om}\Cal{L}\phi
(\la)\Cal{L}f (\la)= \lim_{\la\to \om}\Cal{L}\phi (\la)\hat{f}
(\om)$.

\noindent This implies $\lim_{\la\to \om}\Cal{L}\phi (\la)=0$. \P

\enddemo

\proclaim{Proposition 1.4} (i) Let $\phi\in E(\Lambda_{\jj})$ and
$\beta \in \r\setminus \Lambda$. Then

(1.11)\qquad $\psi'-i\beta \psi = \phi$

\noindent has a unique solution $\psi =
\gamma_{\beta}[P(\gamma_{-\beta}\phi)- \overline {\Cal
{L}(\gamma_{-\beta}\phi)}(0)]$ where $ \psi\in \Lambda_{\jj}$.

(ii) $\sigma(A_{\Cal{L}})=i\, \Lambda$ and
$\sigma(A_{\Cal{C}})=i\, \Lambda$.
\endproclaim
 \demo{Proof}(i) Directly, one can check that $\psi$ satisfies (1.11). For $\jj=\r_+$,  Ingham's
 Theorem gives
  $\psi \in BUC(\r_+,X)$ and for $\jj=\r$,  [BP] gives $\psi \in BUC(\r,X)$. Simple calculations and
   Proposition 1.1(iii) show that $sp^{\Cal{L}}(\psi)\st sp^{\Cal{L}}(\phi)$
  and $sp^{\Cal{C}}(\psi)\st sp^{\Cal{C}}(\phi)$. This proves $ \psi\in \Lambda_{\jj}$.
   Uniqueness follows because if $\psi'-i\beta \psi
  =0$ then $\psi=e^{i\beta t}x$ which gives
  $sp_{L}(\psi)=\{\beta\}$ (respectively  $sp_{C}(\psi)=\{\beta\}$) or $x=0$.

(ii) Since $e^{i\beta t}x|\,\jj$ are eigen-vectors for
$A_{\Cal{L}}$ and $A_{\Cal{C}}$ respectively for all $\beta\in
\Lambda$, $0\not = x\in X$, one gets $ i\, \Lambda \st\sigma
(A_{\Cal{L}})$ and $ i\, \Lambda \st\sigma (A_{\Cal{C}})$
respectively. That $ \sigma (A_{\Cal{L}}) \st i\, \Lambda$  and  $
\sigma(A_{\Cal{C}}) \st i\, \Lambda$ follows by (i) and the closed
graph theorem. \P
\enddemo

Let $\Lambda $ be a closed subset of $\r$.  Then

\qquad $E(\Lambda_{\r})=\{\phi \in BUC(\r,X): sp_{C}(\phi) \st
\Lambda\}$.

\noindent is a closed subset of $BUC(\r,X)$. But
  $E(\Lambda_{\r_+})=\{\phi \in BUC(\r_+,X): sp_{L}(\phi) \st
\Lambda\}$
 is not closed in  $BUC(\r_+,X)$. Denote by

\qquad  $\overline {E(\Lambda_{\r_+})}$ the closure  of
$E(\Lambda_{\r_+})$ in $BUC(\r_+,X)$.

\noindent Then $E(\Lambda_{\r})$ respectively
$\overline{E(\Lambda_{\r_+})}$ is a closed subspace of $BUC(\r,X)$
(respectively $BUC(\r_+,X)$) which is invariant under translations
from $\jj$, by Proposition 1.1 (i). Moreover, $\gamma_{\om}x|\,\jj
\in E(\Lambda_{\jj})$ for all $x\in X$ and $\om \in \Lambda$. When
$\Lambda=\emptyset$, we set

 $\overline
{E(\Lambda_{\r_+})}=E(\emptyset)$ and $ E_{ E(\emptyset)}
(\Lambda_{\r_+})=\overline {E(\Lambda_{\r_+})}/E(\emptyset) $.

Denote by $S_{\Lambda}^{+}(t)$  and $S_{\Lambda}(t)$ the
restrictions of $S^+(t)$ and $S(t)$ on
$\overline{E(\Lambda_{\r_+})}$ and $E (\Lambda_{\r})$
respectively. It follows $S_{\Lambda}^+(t)$ and $S(t)$, $t\in
\r^+$ are $C_0$-semigroups whose generators we denote by
$A_{\Cal{L}}$ and $A_{\Cal{C}}$ respectively.

\proclaim{Corollary 2.3} For $\phi\in BUC(\jj,X)$, one has
$sp^{\Cal{L}}(\phi)= sp^{\Cal{L}_u} (\phi)$ for $\jj=\r_+$ and
$sp^{\Cal{C}}(\phi)= sp^{\Cal{C}_u} (\phi)$ for $\jj=\r$.
\endproclaim
\demo{Proof} We consider first the case $\jj=\r_+$. Take $\Lambda=
sp^{\Cal{L}}(\phi)$. By Proposition 1.3 (ii) and Proposition 2.1,
one has $\sigma (A_{\Cal{L}})=i sp^{\Cal{L}} (\phi)=i \cup_{\xi\in
\Lambda_{\Cal{L}}}S^+(\cdot)\xi $. This implies $ sp^{\Cal{L}}
(S^+ (\cdot)\phi)\st sp^{\Cal{L}}(\phi)$.  By Lemma  1.3 (i),
$\Cal{L}_u\phi= \Cal{L}S^+(\cdot)\phi$ and so $ sp^{\Cal{L}} (S^+
(\cdot)\phi)= sp_{L_u}(\phi)$. By (1.10), concludes
$sp^{\Cal{L}}(\phi) =sp^{\Cal{L}_u}(\phi)$.

The case $\jj=\r$ is similar using in the above Proposition 2.2
instead of Proposition 2.1. \P
\enddemo

\demo{Proof} First, we prove the case $\phi \in LAP_{ub}(\r,X)$.
Set
 $L(\phi)=\overline{\text{\, span}}\{\phi_t: t\in\r\}$ and
$L^+(\phi)=\overline{\text{\, span}}\{(\phi|\r_+)_t: t\in\r_+\}$.
By [bb,diss], $m: L(\phi) \to L^+(\phi)$ defined by  $m\psi
=\psi|\, \r_+$ is an isomorphic isometric map. since
$S(\cdot)\phi$ is Bochner integrable in $BUC(\r,X)$, $\Cal{L}
(S(\cdot)\phi (\la)\in L(\phi)$ for each $\la\in \cc_+$ (see also
[bb, diss Lemma 1.2.2]) and hence its extension
$\overline{\Cal{L}} (S(\cdot)\phi (\la)\in L(\phi)$. This implies
$m(\overline{\Cal {L}} S(\cdot)\phi(\la)) = \overline{\Cal {L}}
S^+(\cdot)\phi|\,\r_+ (\la)$ and hence $sp^{\Cal{L}}
(S^+(\cdot)\phi|\,\r_+)=sp^{\Cal{L}} (S(\cdot)\phi)$. By (2.7) and
Corollary 2.4, we conclude $sp^{\Cal{L}} (\phi)=sp^{\Cal{C}}
(\phi)$.

 Now, assume  $\phi \in
LAP_{b}(\r,X)$. Then $M_h\phi \in LAP_{ub}(\r,X)$ for each $h>0$.
So, by Proposition 1.1 (ii), we conclude $sp^{\Cal{L}}
(\phi)=\cup_{h>0} sp^{\Cal{L}} (M_h\phi)=\cup_{h>0} sp^{\Cal{C}}
(M_h\phi)=sp^{\Cal{C}} (\phi)$. \P
\enddemo

Problem 21.07. Let $S$ be a finite  set of non-negative real
numbers containing $0$ with the property that if $a, b \in S$,
then $a+b$ or $|a-b|\in S$. Prove that

(i) Either $S$ consists of exactly four different real numbers;

(ii) Or there is $n\in \{0,1,2,3,\cdots \}$ and positive real
number $r$ such that $S=\{0,r, 2r,  \cdots, nr\}$

Solution. Assume $ 0\not =a,b\in S$ and $b\not = 2a$. Then
$S=\{0,a, b-a,a\}$ satisfies (i). If $|S|\le 3$, then

$S=\{0\}$ or $S=\{0,a\}$ or $S=\{0,a, 2a\}$.

So, assume $n > 3$ and $S=\{0,s_1,s_2, s_3,\cdots , s_n\}$ and
$0<s_1<\cdots <s_n$. We prove

\qquad (1) $s_n-s_j=s_{n-j}$, \qquad $j=1,2,\cdots, n-1$,

\qquad (2) $s_{n-1}-s_j=s_{n-j-1}$, \qquad $j=1,2,\cdots, n-2$,

\noindent Indeed, (1) follows because $ s_n < s_n+s_j\not\in S$
for $j=1,2,\cdots, n-1$. So, $ s_n > s_n -s_j \in S$ for
$j=1,2,\cdots, n-1$. This implies (1) because $s_1 \le
s_n-s_{n-1}$. Now, we prove (2) for $n\ge 4$. We have
$s_{n-1}+s_j\not \in S$ for $j=2,\cdot,n-1$ because $s_{n-1}+s_j
> s_{n-1}+s_1 =s_n$. This implies $s_{n-1}-s_j \in S$ for
$j=2,\cdot,n-1$.  It follows $0 <s_{n-1}-s_{n-2}=s_n -
(s_{n-2}+s_1)=s_n - (s_{n-2}+s_2) +s_2 -s_1= s_2 -s_1\ge s_1$.
Since $s_2 -s_1\le s_2$, one gets $s_{n-1}-s_{n-2}=s_1$. It
follows also, $s_{n-1}-s_1=s_{n-2}\in S$.  So, $s_{n-1}-s_j\in S$,
\qquad $j=1,2,\cdots, n-2$. This implies (2) because $ s_1\ge
s_{n-1}-s_j \le s_{n-2}$.

From (1), (2), one gets $s_k= k s_1$, $k=1,2,\cdot,n$.

 .........................................

 Let $ABC$ be a triangle
with acute angels and $M$ a point on $AC$. Let $S_1$ be the circle
passing through $ABM$ and $S_2$ the circle passing through $MBC$.
Show that the area of the intersection $S_1\cap S_2$ is minimum if
$BM$ is perpendicular to $AC$.

Solution. Let $O_1$ be the center of $S_1$ and $O_2$ be  the
center of $S_2$. One has  $\angle BO_1M = 2\angle BAM $ and
$\angle BO_12M = 2\angle BCM $. So, $\angle BO_1M $ and $\angle
BO_12M$ do not depend on the position of $M$.

 The area of the intersection $S_1\cap S_2$ is minimum when the
 radii  $BO_1$ and $BO_2$ are minimum. Since $BO_1= (1/2)BM\sin\angle
 BAM
 $ and  $BO_2= (1/2)BM\sin\angle BCM$ one gets the
 radii  $BO_1$ and $BO_2$ are minimum when $BM$ is minimum. But
$BM$ is minimum when  it is perpendicular to $AC$.

.

.

.

.

 ...................................................

Let $ABCD$ be a square and let $E$ be a point on its diagonal
$BD$. Suppose that $O_1$ is the center of the circle passing
through $ABE$ and $O_2$ of the circle passing through $ADE$. Show
that $AO_1 E O_2$ is a square.

 Solution. One has

 $AO_1=O_1E$, $AO_2=O_2 E$,

$\angle AO_1 E= 2 \angle ABE=90 ^{\circ}$ and $\angle AO_2 E= 2
\angle ADE=90 ^{\circ}$.

This implies $\angle O_1 A E= \angle O_1  E A= \angle O_2 A E=
 \angle A E O_2= 45^{\circ} $. Therefore,
triangle  $AO_1E$ is congruent to triangle  $AO_2E$. This implies
$AO_1 E O_2$ is a square.

.

.

.

.

.

.

.

.

.

.
 \pagebreak

.........................................................

 22.06 Let $a,b
$ and $c$ be positive real numbers. Determine all non-negative
real numbers $x$ such that

(*) $\frac{a}{b+x} + \frac{b}{c+x}+\frac{c}{a+x}\ge
\frac{3(a+b+c)}{a+b+c+3x}$

Solution. Without loss of generality we can assume $a \le b\le c$.
Set $a+b+c =A$. One has

$\frac{a}{b+x}-\frac{a+b+c}{a+b+c+3x}=
\frac{(a-b)(A+x)}{(b+x)(A+3x)}+ \frac{a-c}{(b+x)(A+3x)}= I_1+
I_2$,

$\frac{b}{c+x}-\frac{a+b+c}{a+b+c+3x}=
\frac{(b-c)(A+x)}{(c+x)(A+3x)}+ \frac{b-a}{(c+x)(A+3x)}= J_1+
J_2$,

$\frac{c}{a+x}-\frac{a+b+c}{a+b+c+3x}=
\frac{(c-a)(A+x)}{(a+x)(A+3x)}+ \frac{c-b}{(a+x)(A+3x)}= K_1+
K_2$.

We show that $I_1+J_1+K_1 \ge 0$ and $I_2+J_2+K_2 \ge 0$ for each
$x\ge 0$.

We have

$I_1\ge \frac{(a-b)(A+x)}{(a+x)(A+3x)}$, $J_1\ge
\frac{(b-c)(A+x)}{(a+x)(A+3x)}$,
  $I_1+J_1+K_1
\ge\frac{(a-b+b-c+c-a) (A+x)}{(a+x)(A+3x)} = 0$,

$I_2+J_2 \ge \frac{a-c}{(c+x)(A+3x)}+\frac{b-a}{(c+x)(A+3x)}=
\frac{b-c}{(c+x)(A+3x)}$ and hence

  $I_2+J_2+K_2 \ge
\frac{(c-b)}{A+3x} [\frac {1}{a+x}-\frac{1}{c+x}] \ge 0$.

This shows that (*) holds for each $x\ge 0$.

\pagebreak ............

 Given three circles $S_1$, $S_2$, $S$. $S_1$ and
$S_2$ passes through the center  $O$ of $S$.  $S_1$  intersects
$S_2$ in $M$. $S_1$  intersects  $S$ in $C$, $D$ and $S_2$
intersects $S$ in $A$, $E$. The segment $AD$ extended  intersects
the segment $CE$ extended  in $B$ and $B\not = M$. Prove that
$\angle BMO= \pi/2$.

\noindent Solution. Denote by $K$ the middle point of $CD$. Then
$OK$ is perpendicular to $CD$. Since $\angle BMO= \angle BMD+
\angle DMO$,
 it is enough to prove

(1) $\angle DMO=\angle KCO=\alpha$.

(2) $\angle BED=\angle COK$.

(3) $\angle BMD=\angle BED$.

(1) follows  because $BDOC$ lies on the circle $S_1$.

\noindent Since $\angle DEC =\pi-\angle BED$, $\angle DEC +\angle
DAC=\pi$, it follows $ \angle DAC=\angle DEM$. But $2 \angle DAC=
\angle DOC$. Hence (2) follows.

\noindent To prove (3), we need to show $BMED$ lies on a circle
which is equivalent to prove $\angle DME =\angle DBE =
\alpha+\beta$. We have

$\angle ADC=\angle ABC+\angle BCD= \angle ABC+\gamma$, $\angle
ADC=\alpha+ \angle ODA=\alpha+\gamma+\beta$ because $OD=OA$. This
gives $\angle ABC= \alpha+\beta$ and proves (3).

\enddocument
Solve the system of equations

(1)\qquad $(1+x)(1+x^2)(1+x^4)=1+y^7$,

(2)\qquad $(1+y)(1+y^2)(1+y^4)=1+x^7$.

Solution. Assume $x=y$. Then $x=y=1$ is not a solution and hence
$(1+x)(1+x^2)(1+x^4)= \frac{1-x^8}{(1-x} =1+x^7$. This implies
$x^7 =x$. Hence $x=0$ or $x=-1$. Therefore, $x=y=0$,  $x=y=-1$ are
two solutions.

We show that if $x\not = y$, then $(x,y)$ is not a solution.

If $x=0$, then $y=0$. If $x>0$ and $y<0$, then
$(1+x)(1+x^2)(1+x^4)> 1$ and $1+y^7<1$ so (1) is not satisfied by
(x,y). By symmetry, $(x,y)$ with $y> 0$, $x< 0$ is not a solution.

Assume $ x> y >0$. Then $(1+x)(1+x^2)(1+x^4)> 1+x^7  >1+y^7$, so
$(x,y)$ is not a solution. By symmetry $(x,y)$ is not a solution
if $y> x >0$.

Assume $ x <0$,  $y < 0$ and $x< y$. Multiply (1) by $1-x$ and (2)
by $1-y$ and subtract, we get

$y^8- x^8= y-x + y^7 -x^7+ xy (x^6- y^6)$.  One has $y^8- x^8<0$
but $ y-x + y^7 -x^7+ xy (x^6- y^6) >0$, implying $(x,y)$ is not a
solution. By symmetry if  $ x <0$,  $y < 0$ and $y< x$, then
$(x,y)$ is not a solution.

\noindent Determine all pairs $(x,y)$ of real numbers such that

(1)\qquad $(x+1)(x^2+1)(x^4+1)=y^7+1$ and
$(y+1)(y^2+1)(y^4+1)=x^7+1$.

\noindent Solution. From (1) one gets

$(x+1)(x^2+1)(x^4+1)(y+1)(y^2+1)(y^4+1)=(y^7+1) (x^7+1)$.

\noindent Obviously, $(-1,-1)$ is a solution to (1). Now, assume
$x\not =-1$ and $y\not=-1$. Divide both sides by $(x+1)(y+1)$ and
simplify, one gets

(2)\qquad $(x^2+1)(x^4+1)(y^2+1)(y^4+1)=[(x^2+1)(x^4+1)-x(1+x^2+
x^4)]$

\qquad\qquad \qquad\qquad \qquad\qquad $[(y^2+1)(y^4+1)-y(1+y^2
+y^4)]$

\noindent From (2) one can check directly that if $x< y$ or $x> y$
then $(x,y)$ is not a solution for (2). (One can check only the
cases $x>y >0$, $x<y <0$ and $ 1 > x > -y >0$. The last case is
easily seen from the first part  of (1)).

 So, the possible solutions
are the pairs $(x,y)$ with $x=y$. But this implies

(3)\qquad $(x^2+1)(x^4+1)=(x^2+1)(x^4+1)-x(1+x^2+x^4)$;

\noindent and the only real solution for (3) is $x=0$. Hence
$(0,0)$ is the only solution to (2).

\noindent So, the  solutions to (1) are the pairs $(-1,-1)$,
$(0,0)$.

\enddocument
**************************************

\pagebreak \pagebreak ............

 Given three circles $S_1$, $S_2$, $S$. $S_1$ and
$S_2$ passes through the center  $O$ of $S$.  $S_1$  intersects
$S_2$ in $M$. $S_1$  intersects  $S$ in $C$, $D$ and $S_2$
intersects $S$ in $A$, $E$. The segment $AD$ extended  intersects
the segment $CE$ extended  in $B$ and $B\not = M$. Prove that
$\angle BMO= \pi/2$.

\noindent Solution. Denote by $K$ the middle point of $CD$. Then
$OK$ is perpendicular to $CD$. Since $\angle BMO= \angle BMD+
\angle DMO$,
 it is enough to prove

(1) $\angle DMO=\angle KCO=\alpha$.

(2) $\angle BED=\angle COK$.

(3) $\angle BMD=\angle BED$.

(1) follows  because $BDOC$ lies on the circle $S_1$.

\noindent Since $\angle DEC =\pi-\angle BED$, $\angle DEC +\angle
DAC=\pi$, it follows $ \angle DAC=\angle DEM$. But $2 \angle DAC=
\angle DOC$. Hence (2) follows.

\noindent To prove (3), we need to show $BMED$ lies on a circle
which is equivalent to prove $\angle DME =\angle DBE =
\alpha+\beta$. We have

$\angle ADC=\angle ABC+\angle BCD= \angle ABC+\gamma$, $\angle
ADC=\alpha+ \angle ODA=\alpha+\gamma+\beta$ because $OD=OA$. This
gives $\angle ABC= \alpha+\beta$ and proves (3).

\enddocument